\documentclass[reqno]{amsart}
\usepackage{amsmath,amsthm, amsbsy, amsfonts, amssymb}
\usepackage{euscript, mathrsfs,stmaryrd} 
\usepackage{mathbbol} 
\usepackage{tabularx, booktabs, supertabular} 
\usepackage{mathtools} 
\usepackage{graphicx, amscd, tikz-cd, tikz} 
\usetikzlibrary[matrix, arrows, arrows.meta, calc,
  decorations.pathmorphing, backgrounds, positioning, fit, petri, cd]
\usepackage{setspace} \usepackage{hyperref, url}

\graphicspath{
  {./}
  {figures/}
}

\newtheorem{theorem}{Theorem}[section]
\newtheorem{thm}[theorem]{Theorem}
\newtheorem{lemma}[theorem]{Lemma}
\newtheorem{proposition}[theorem]{Proposition}
\newtheorem{prop}[theorem]{Proposition}

\newtheorem{cor}[theorem]{Corollary}

\newtheorem{mainthm}{Theorem}

\newtheorem{maincor}[mainthm]{Corollary}

\theoremstyle{definition}
\newtheorem{definition}[theorem]{Definition}
\newtheorem{defn}[theorem]{Definition}
\newtheorem{example}[theorem]{Example}

\newtheorem{remark}[theorem]{Remark}
\newtheorem{rem}[theorem]{Remark}
\newtheorem{notation}[theorem]{Notation}
\theoremstyle{remark}
\newcommand{\C}{\mathbb{C}}  
\newcommand{\F}{\mathbb{F}}  
\newcommand{\N}{\mathbb{N}}  
\newcommand{\Q}{\mathbb{Q}}  
\newcommand{\R}{\mathbb{R}}  
\newcommand{\Z}{\mathbb{Z}}  
\newcommand{\K}{\mathbb{K}}  

\newcommand{\abs}[1]{\lvert#1\rvert}
\newcommand{\eu}[1]{\EuScript{#1}}

\newcommand{\mbb}[1]{\mathbb{#1}}
\newcommand{\mc}[1]{\mathcal{#1}}

\newcommand{\msf}[1]{\mathsf{#1}}
\newcommand{\op}[1]{\operatorname{#1}}

\newcommand{\cO}{\mathcal{O}}

\newcommand{\ip}[2]{\langle #1 , #2 \rangle }
\newcommand{\tr}{\msf{T}}
\newcommand{\htr}{\msf{H}}

\newcommand{\hyp}{\msf{H}}
\newcommand{\proj}{\msf{P}}  
  
\newcommand{\clatt}{\mbb{L}}
\newcommand{\ball}{\msf{B}}  

\newcommand{\ver}{\msf{V}}
\newcommand{\edg}{\msf{E}}
\newcommand{\size}[1]{|#1|}

\newcommand{\points}{\mc{P}}
\newcommand{\lines}{\mc{L}}
\newcommand{\I}{I}
\newcommand{\II}{I\!I}

\newcommand{\sgn}{\epsilon}

\newcommand{\into}{\hookrightarrow}
\newcommand{\onto}{\twoheadrightarrow}

\newcommand{\bx}{\mathbf{x}}
\newcommand{\by}{\mathbf{y}}
\newcommand{\ba}{\mathbf{a}}
\newcommand{\one}{\mathbf{1}}
\newcommand{\spanning}{\mathcal{S}}
\newcommand{\basis}{\mathcal{B}}
\newcommand{\dual}{\vee}
\newcommand{\eisen}{\mc{E}}
\newcommand{\gauss}{\mc{G}}
\newcommand{\hurwitz}{\mc{H}}
\newcommand{\petersen}{\eu{P}}

\newcommand{\kneser}{KG}
\newcommand{\McL}{\textsc{McL}}
\newcommand{\graph}{\textsc{Graph}}
\newcommand{\bigraph}{\textsc{BiGraph}}
\newcommand{\gram}{\textsc{Gram}}

\newcommand{\lattice}{\textsc{Lattice}}

\newcommand{\vecDelta}{\mathbb{\Delta}}

\DeclareMathOperator{\aut}{Aut}
\DeclareMathOperator{\inc}{Inc}
\DeclareMathOperator{\Refl}{\textsc{Refl}}
\DeclareMathOperator{\cox}{\textsc{Cox}}

\DeclareMathOperator{\spec}{\textsc{Spec}}
\DeclareMathOperator{\rad}{rad}

\newcounter{tbcomments}

\newcounter{jmcomments}

\begin{document}

\title[Lattices from Graphs]{Lorentzian Lattices from Graphs}

\author{Tathagata Basak}
\email{tathagat@iastate.edu}
\address{Department of Mathematics, Iowa State University
  Ames, IA 50011} 

\author{Jon McCammond}
\email{jon.mccammond@math.ucsb.edu}
\address{Department of Mathematics, UC Santa Barbara, 
  Santa Barbara, CA 93106} 
\thanks{TB was partially supported by Simons Foundation Collaboration Grant 637005.}






\date{\today}

\begin{abstract}
  In a classic 1983 article, Neumaier and Seidel showed how to
  construct a real Lorentzian lattice from any graph whose second
  eigenvalue is an integer.  We extend their ``discrete hyperbolic
  geometry'' construction to more general rings of integers and
  highlight its functorial nature.  We also introduce a variation that
  produces complex Lorentzian lattices in complex Lorentzian spaces.
  The complex version produces uniform constructions of several
  complex Lorentzian lattices in the literature that are preserved by
  arithmetic complex reflection groups which leads to interesting ball
  quotients.  Examples include (1) a complex reflection group in
  $U(4,1)$ that corresponds to the moduli space of cubic surfaces in
  $\proj^3_\C$, (2) a complex reflection group in $U(9,1)$, found by
  Deligne-Mostow and Thurston, that corresponds to the moduli space of
  $12$ points on $\proj^1_\C$, and (3) a complex reflection group in
  $U(13, 1)$ that plays a central role in the Daniel Allcock's
  monstrous proposal.
\end{abstract}

\maketitle

\section{Introduction}\label{sec:intro}

This paper builds real Lorentzian lattices and complex Lorentzian
lattices from graphs, with a particular emphasis on strongly regular
graphs when working over $\R$, and symmetric $2$-designs when working
over $\C$.  These types of constructions are not completely new,
particularly over the reals; they are inspired by the straight-forward
way one builds root lattices from Dynkin diagrams.  More immediately,
our main result over the reals is a slight extension of observations
in \cite{NeumaierSeidel83} where Neumeier and Seidel construct a real
Lorentzian lattice from any graph whose second eigenvalue is an
integer.  We extend their ``discrete hyperbolic geometry''
construction to more general rings of integers and we highlight its
functorial nature, showing that some induced subgraphs lead to lattice
embeddings.  Given $a, b, \dotsb \in \C$, let $\cO(a, b, \dotsb)$
denote the ring of integers in the number field $\Q(a, b, \dotsb)$.
Roughly speaking, we prove the following.

\begin{mainthm}[Theorem~\ref{thm:latt-functor}]\label{main:latt-functor}
  If $r$ is a real algebraic integer and $\cO(r)$ is a PID, there is
  a functor from graphs with second eigenvalue $r$ to Lorentzian
  $\cO(r)$-lattices.
\end{mainthm}

\noindent 
For a more precise statement see \S\ref{sec:graph-lattice}.  The
variation over $\C$, which we believe is new, produces complex
Lorentzian lattices and interesting ball quotients with nontrivial
moduli interpretations.  Roughly speaking, we prove the following.

\begin{mainthm}[Theorem~\ref{thm:clatt-functor}]\label{main:clatt-functor}
  If $r$ is a real algebraic integer and $\zeta$ is a root of unity so
  that $\cO = \cO(r \zeta, \overline{r \zeta})$ is a PID, then there
  is a functor from bipartite graphs with second eigenvalue $r$ to
  complex Hermitian Lorentzian $\cO$-lattices.
\end{mainthm}

\noindent
For a more precise statement see \S\ref{sec:bipartite}.  The proof
idea in both cases is very straightforward.  For the real version, one
creates an $(r,0,-1)$ representation
\cite[p.88]{BrouwerCohenNeumaier89} of the graph $\Delta$ by taking a
basis $\{e_i\}$ indexed by the vertices of $\Delta$ and declaring that
$\ip{e_i}{e_i}=r$ for all $i$, $\ip{e_i}{e_j}=-1$ when $\{i,j\}$ is an
edge and $\ip{e_i}{e_j}=0$ when $\{i,j\}$ is not an edge.  Because $r$
is an eigenvalue, this form has a non-trivial radical.  In the
quotient that kills the radical, the basis becomes a spanning set
$\{\alpha_i\}$ and the form becomes non-degenerate.  For the complex
version we start with a bipartite graph $\Delta$ and then its real
Lorentzian lattice is complexified, stretched, and one class of
vertices is rotated by a unit complex number $\zeta$.

The most striking illustrations of Theorem~\ref{main:latt-functor}, as
Neumaier and Seidel were aware, are Lorentzian $\Z$-lattices derived
from strongly regular graphs with second eigenvalue $1$ or $2$.
Section~\ref{sec:kneser} discusses lattices derived from an infinite
family of Kneser graphs $\kneser(n,2)$ with $r=1$ and
Section~\ref{sec:mclaughlin} discusses lattices from a dozen
interesting graphs inside the $276$-vertex McLaughlin $2$-graph with $r=2$,
with special emphasis on the $275$-vertex McLaughlin graph.

These sections highlight the benefits of the integer-valued Lorentzian
form.  The McLaughlin $2$-graph, for example, contains several
strongly regular graphs $\Delta$ whose automorphism group is,
essentially, a sporadic finite simple group and the group
$\aut(\Delta)$ acts on the Lorentzian representation fixing a
$1$-dimensional subspace. Projecting the spanning vectors onto its
orthogonal complement produces the Euclidean representation of
$\Delta$ \cite{BrouwerVanMaldeghem22} and the action of $\aut(\Delta)$
on this space is often the smallest irreducible representation of the
sporadic finite simple group.  Concretely, the graphs inside the
McLaughlin $2$-graph have Lorentzian $\Z$-lattices inside
$\II_{25,1}$, which include and extend their Euclidean representations
in the Leech lattice.  By adding one extra dimension to the smallest
irreducible representation, one obtains a representation of the
sporadic group defined over the integers preserving an integral
Lorentzian form.  In fact, the Lorentzian lattice provides a
representation of the much larger infinite group $\cox(\Delta) \rtimes
\aut(\Delta)$ where $\cox(\Delta)$ is the simply-laced Coxeter group
defined by $\Delta$.  See Remark~\ref{rem:large-radical}.

Using Lorentzian representations we have discovered a result about
equipartitions among disjoint $r$-regular induced subgraphs inside a
connected graph with integral second eigenvalue $r$.  We explore this
result and some of its graph theoretic consequences in a separate
short preprint \cite{BaMc-graphs}.  The integrality of the Lorentzian
forms makes some of the calculations in recent articles look more
natural.  See \cite{BonadrenkoRadchenko13, BPR17, BMPRV18}.  And the
integrality of the form makes it possible to consider reduction modulo
a prime.  These observations might prove useful in answering existence
questions about strongly regular graphs.

The complex variants of our graph-to-lattice construction give uniform
definitions of several remarkable complex Lorentzian lattices that
occur in the literature starting with the seminal work of Deligne and
Mostow \cite{DeligneMostow86}.  To explain the context for our
results, we discuss two pairs of such lattices named $\eisen_3
\subseteq \eisen_4$ and $\eisen_9 \subseteq \eisen_{13}$.  These
lattices are defined over the ring of Eisenstein integers $\eisen =
\Z[\omega]$ where $\omega$ is a cube-root of unity.  Our naming
convention is that if the Lorentzian Eisenstein lattice is in
$\C^{n,1}$, we call the graph $\Delta_n$, the lattice $\eisen_n$, and
the corresponding complex hyperbolic reflection group $\Gamma_n$.  So
$\eisen_n \otimes \C \simeq \C^{n, 1}$ and $\Gamma_n = \proj
\Refl(\eisen_n)$ is a discrete subgroup of $\proj U(n, 1)$ acting on
the complex $n$-ball $\ball^n$ preserving its complex hyperbolic
metric.

\subsection*{Example 1 and 2}  Let $\Delta_4 = \petersen = \kneser(5, 2)$
be the Petersen graph.  Let $\Z_\petersen$ be the real Lorentzian
$\Z$-lattice in $\R^{4,1}$ constructed in \S\ref{sec:kneser} from
$\petersen$ and let $\eisen_4 = \Z_\petersen \otimes \eisen$ be the
complex Lorentzian lattice of signature $(4, 1)$ obtained by tensoring
$\Z_\petersen$ with the Eisenstein integers and promoting the bilinear
form to a Hermitian form.  Let $\Gamma_4 = \proj \Refl(\eisen_4)
\subseteq \proj U(4, 1)$.  The largest induced cycle in $\petersen$ is
a $6$-cycle or the affine Dynkin diagram $ \widetilde{A}_5$. According
to theorem \ref{main:latt-functor}, the subgraph $\Delta_3 =
\widetilde{A}_5$ yields a signature $(3, 1)$ sub-lattice $\eisen_3
\subseteq \eisen_4$.  Allcock, Carleson and Toledo proved in
\cite{ACT02-cubic-surfaces} that the ball quotient $X_4 = \Gamma_4
\backslash \ball^4$ is the moduli space of stable cubic surfaces in
$\proj^3_{\C}$. The sub-ball quotient $X_3 = \Gamma_3 \backslash
\ball^3$ is the moduli space of six points in $\proj^1_{\C}$.

The most striking uses of Theorem~\ref{main:clatt-functor} are special
cases that produce lattices over the Eisenstein integers $\eisen =
\Z[\omega]$ and lattices over the Gaussian integers $\gauss = \Z[i]$.

\begin{maincor}[Eisenstein lattices]\label{maincor:eisen}
  There is a functor from bipartite graphs with second eigenvalue
  $r=\sqrt{3}$ to Lorentzian $\eisen$-lattices.
\end{maincor}

\begin{maincor}[Gaussian lattices]\label{maincor:gauss}
  There is a functor from bipartite graphs with second eigenvalue
  $r=\sqrt{2}$ to Lorentzian $\gauss$-lattices.
\end{maincor}

Although we are only aware of finitely many bipartite graphs with
these properties, the Eisenstein and Gaussian lattices they construct
are quite exceptional.  Below we briefly discuss a pair of Eisenstein
examples $\eisen_{9} \subseteq \eisen_{13}$ constructed using
Corollary~\ref{maincor:eisen}.  There is a similar pair of Gaussian
examples $\gauss_5 \subseteq \gauss_7$ constructed using
Corollary~\ref{maincor:gauss}; these are discussed in \ref{ex:8-gon}
and \ref{ex:gauss-heawood}.

\subsection*{Example 3}  
The simplest bipartite graphs to which Corollary~\ref{maincor:eisen}
can be applied is the $12$-cycle or the extended Dynkin diagram
$\widetilde{A}_{11}$.  Since the Lorentzian Eisenstein lattice
Corollary~\ref{maincor:eisen} constructs from this graph is in
$\C^{9,1}$, we call the graph $\Delta_{9}$, the lattice $\eisen_{9}$,
and its complex hyperbolic reflection group $\Gamma_{9}$.  The
classification of arithmetic complex reflection groups in $\proj
U(n,1)$ is far from complete and very few examples are known in high
dimensions.  Deligne and Mostow constructed a family of examples using
monodromy of hypergeometric functions \cite{DeligneMostow86}, and
Thurston independently found many of these examples while studying the
moduli space of certain triangulations of a $2$-sphere
\cite{Thurston98}. There is a single complex reflection group in
$\proj U(9,1)$ which is the ``largest'' among both these classes of
examples, and it is the reflection group $\Gamma_9 = \proj
\Refl(\eisen_9)$.  The ball quotient $X_9 = \Gamma_9 \backslash
\ball^9$ turns out to be the moduli space of $12$ points on
$\proj^1_\C$.  See Example~\ref{ex:12-gon}.

\subsection*{Example 4}  
At the other extreme, the largest bipartite graph we know with second
eigenvalue $\sqrt{3}$ is the $26$ vertex incidence graph of the finite
projective plane over $\F_3$.  Since the Lorentzian Eisenstein lattice
Corollary~\ref{maincor:eisen} constructs from this graph is in
$\C^{13,1}$, we call the graph $\Delta_{13}$, the lattice
$\eisen_{13}$, and its projectivized automorphism group $\Gamma_{13} =
\proj \Refl(\eisen_{13}) = \proj \aut(\eisen_{13}) \subset \proj
U(13,1)$.  Since $\Delta_9$ is a maximal circuit inside $\Delta_{13}$,
the functoriality in Corollary~\ref{maincor:eisen} implies $\eisen_9
\subseteq \eisen_{13}$.  The largest dimension in which an arithmetic
complex hyperbolic group is known to exist is thirteen, and the group
$\Gamma_{13}$ is the unique known example in this dimension
\cite{Allcock00-leech}.  See Example~\ref{ex:monstrous}.  \par The
group $\Gamma_{13}$ has an intriguing conjectured connection with the
Monster simple group $M$.  The group $\Gamma_{13}$ acts on the complex
hyperbolic $13$-ball $\ball^{13}$ preserving its $U(13,1)$-invariant
metric.  Let $\eu{H}_{13} \subseteq \ball^{13}$ be the union of the
mirrors of the complex reflections in $\Gamma_{13}$.  Daniel Allcock's
``monstrous proposal'' \cite{Allcock09-monstrous} conjectures that the
group $M \wr 2$ (known as the bimonster) is a natural explicit
quotient of the orbifold fundamental group $G$ of $X_{13}^{\circ} =
\Gamma_{13} \backslash (\ball^{13} - \eu{H}_{13})$, i.e.  the
complement of the discriminant divisor in the ball quotient $X_{13} =
\Gamma_{13} \backslash \ball^{13}$.

There are many similarities between the two pairs of examples
$\eisen_3 \subseteq \eisen_4$ and $\eisen_9 \subseteq \eisen_{13}$
discussed above.  Consider one of these four examples $(\Delta,
\Gamma, X, \ball) = (\Delta_n, \Gamma_n, X_n, \ball^n)$ where $n = 3,
4, 9$ or $13$.  The principal eigenvector of the graph $\Delta$ always
determines a special point $\tau$ in the complex ball $\ball$ and our
graph-to-lattice construction produces a set of `roots' $\lbrace
\alpha_i \colon i \in \Delta \rbrace$ indexed by the vertices of the
graph $\Delta$ such that $\Gamma$ contains complex reflections
$\lbrace R_i \colon i \in \Delta \rbrace$ of appropriate order in the
mirrors $\lbrace \alpha_i^{\bot} \colon i \in \Delta \rbrace$.  Let
$X^{\circ} \subseteq X$ be the complement of the mirrors of $\Gamma$
(i.e. fixed points of the complex reflections in $\Gamma$).  For each
$i \in \Delta$, there is a natural loop $g_i \in G$ (called \emph{a
generator of monodromy}) based at $\tau$ and going around the mirror
$\alpha_i^{\bot}$ once.  The following statements hold in each case (a
more detailed statement is contained in
\ref{rem:similarities-summary}) \emph{
\begin{enumerate}
\item $\tau$ is the unique fixed point in $\ball$ under the
  automorphism group of $\Delta$.
\item The mirrors of $\Gamma$ that are closest to $\tau$ are precisely
  $\lbrace \alpha_i^{\bot} \colon i \in \Delta \rbrace$.
\item The reflection group $\Gamma$ is generated by complex
  reflections $\lbrace R_i \colon i \in \Delta \rbrace$.
\item The fundamental group $\pi_1^{\op{orb}}( X^{\circ}, \tau)$ is
  generated by the loops $\lbrace g_i \colon i \in \Delta \rbrace$.
\item The above generators of $\Gamma$ and $\pi_1^{\op{orb}}(
  X^{\circ}, \tau)$ realizes these groups as quotients of the Artin
  group of the graph $\Delta$.
\end{enumerate}
}

These statements are reminiscent of the properties of simple roots
$\lbrace \alpha_i \rbrace$ and Weyl vector $\tau$ in Weyl groups. The
real hyperbolic versions of statements (1) and (2) are also true for
all the real examples studied in \S\ref{sec:kneser} and
\S\ref{sec:mclaughlin}.  Indeed for some of the real examples, one can
build a Kac-Moody Lie algebra that makes this analogy with Weyl groups
precise.  Even in the complex examples, where the analogy is vague and
imprecise, it has still been useful for studying the complex
hyperbolic reflection groups and associated ball quotients.  The
transition from the reflection group $\Gamma$ in (3) to the
fundamental group $\pi_1^{\op{orb}}(X^{\circ}, \tau)$ in (4) is like
transitioning from the symmetric group to the braid group.

We conclude with a few remarks about the role of $\tau$ in the four
parallel examples.  When $n=4$, $\Delta_4$ is the Petersen graph
$\petersen$, $X_4^{\circ}$ is the moduli space of smooth cubic
surfaces, the point $\tau$ corresponds to the Clebsch cubic surface,
and studying the geometry of the moduli space near $\tau$ has led to a
new presentation of $\pi_1^{\op{orb}}( X_4^{\circ}, \tau)$ as a
quotient of the Artin group of $\eu{P}$.  See \cite{ABL}.
When $n=3$, $\Delta_3$ is a $6$-cycle, the quotient $X_3$ is the
moduli space of $6$ points in $\proj^1_{\C}$, and the point $\tau$
corresponds to a symmetric configuration of $6$ points around the
equator in $\proj^1_{\C} \simeq S^2$.
Similarly, when $n=9$, $\Delta_9$ is a $12$-cycle, the ball quotient
$X_9$ is the moduli space of $12$ points in $\proj^1_{\C}$, and the
point $\tau$ corresponds to a symmetric configuration of $12$ points
around the equator in $\proj^1_{\C} \simeq S^2$.
In both cases, there are presentations of
$\pi_1^{\op{orb}}(X^{\circ},\tau)$ as quotients of the corresponding
affine Artin groups, which are very close to spherical braid groups.
For $n=13$, a complete presentation of $G = \pi_1^{\op{orb}}(
X^{\circ}_{13}, \tau)$ is not known, but it is a quotient of the Artin
group of $\Delta_{13}$ generated by twenty-six loops $\lbrace g_i
\colon i \in \Delta_{13} \rbrace$ based at $\tau$.  See
\cite{AB16-braid-like, AB18-generators}.  Daniel Allcock's monstrous
proposal predicts that the quotient of $G$ by the normal subgroup
generated by $\lbrace g_i^2 \colon i \in \Delta_{13} \rbrace$ is $M
\wr 2$. In \cite{AB25-deflation}, it was shown that this quotient is
either $M \wr 2$ or $\Z/2$.  The proof crucially uses the fact that
$X_{9}$ is a (complex) codimension $4$ sub-orbifold of $X_{13}$ and
the Deligne-Mostow description of $X_9$ as the moduli space of twelve
points. As noted above, the inclusion of $X_9$ in $X_{13}$ is a
consequence of $\eisen_9 \subseteq \eisen_{13}$, which follows
immediately from the functorial nature of
Corollary~\ref{maincor:eisen}.

\subsection*{Organization}
The article is structured as follows.  After establishing notation for
spaces with forms in \S\ref{sec:spaces} and Lorentzian lattices in
\S\ref{sec:lattices}, we present the Neumaier-Seidel construction of
Lorentzian lattices from graphs in \S\ref{sec:graph-lattice} with a
focus on functoriality and on more general rings of integers.  The
lattice symmetries that are obvious from the construction are
discussed in \S\ref{sec:groups} along with a standard decomposition of
the Lorentzian space.  Section~\ref{sec:srg} recalls the structure of
strongly regular graphs, including a statement of Neumaier's near
classification of strongly regular graphs with a fixed smallest (or a
fixed second largest) eigenvalue, and \S\ref{sec:srg-lattices}
explores the consequences for the corresponding lattices.  A family of
lattices from strongly regular graphs with $\lambda_2=1$ are described
in \S\ref{sec:kneser} and examples of lattices from strongly regular
graphs with $\lambda_2=2$ are discussed in \S\ref{sec:mclaughlin}.
The remaining sections shift attention to complex lattices.  Our
complex variation of the Neumaier-Seidel construction is introduced in
\S\ref{sec:bipartite} where it is used to construct Eisenstein and
Gaussian Lorentzian lattices from bipartite graphs.  Lattice
symmetries and a complex version of the standard decomposition are
discussed in \S\ref{sec:std-decomposition}.
Section~\ref{sec:sym-2-designs} specializes to the case of symmetric
$2$-designs, which are bipartite analogs of strongly regular graphs.
Section~\ref{sec:proj-planes} discusses the infinite family of
incidence graphs of finite projective planes, and
\S\ref{sec:ball-quoients} describes some examples of complex
Lorentzian lattices (and two quaternionic Lorentzian lattices) that
have appeared in the literature, which can now be constructed in a
more systematic fashion.  Finally, Appendix~\ref{app:algebra}
establishes the basic algebraic facts which are used to uniformly
prove both main theorems.

\subsection*{AI Use Statement}
The main body of the paper contains no AI generated mathematics. Most
of the mathematics was done in 2024 and 2025. AI was only used in
Appendix ~\ref{app:algebra}.  All the results in the Appendix are also
due to the authors; but ChatGPT was used to clean up and fine-tune
some statements and proofs in the Appendix and to proof-read the
Appendix.

\section{Spaces with Forms}\label{sec:spaces}

This section establishes notation for real vector spaces with
symmetric bilinear forms, complex vector spaces with Hermitian
sesquilinear forms, and the projective models of real and complex
hyperbolic spaces.  Let $\K = \R$ or $\C$.  If $M$ is a matrix with
entries in $\K$, then $M^\tr$ is the \emph{transpose} of $M$ and
$M^\htr$ is the \emph{Hermitian transpose} or \emph{conjugate
transpose} of $M$.  If $\K=\R$, then $M^\htr = M^\tr$.  A matrix $M$
is \emph{symmetric} when $M=M^\tr$ and \emph{Hermitian symmetric} when
$M=M^\htr$.

\begin{defn}[Hermitian forms]\label{def:forms}
  Let $V$ be an $n$-dimensional vector space over $\K$ with a fixed
  ordered basis.  Then vectors in $V$ are represented as an $n\times
  1$ column vector of coordinates.  An $n\times n$ Hermitian symmetric
  matrix $M$ defines a \emph{Hermitian form} on $V$ by setting
  $\ip{x}{y} = x^\htr M y$.  So our Hermitian forms are linear in the
  second variable.  When $\K = \R$, Hermitian form just means
  symmetric bilinear form.  A \emph{Hermitian space} means a
  $\K$-vector space with a Hermitian form and a map between them means
  a form preserving linear transformation.  Two Hermitian spaces are
  \emph{isometric} if there is a form preserving linear isomorphism
  between them.  For any $x \in V$, the \emph{norm of $x$} is
  $\ip{x}{x}$, abbreviated as $x^2$.  For any $S \subseteq V$, let
  $S(n) = \{ x \in S\colon x^2 = n \}$ and let $S_< = \{ x \in S
  \colon x^2 < 0 \}$.
\end{defn}

\begin{defn}[Gram matrix]\label{rem:gram-matrix}
  Let $\mathcal{L} = (u_1,\cdots, u_m)$ be an $m$-tuple of vectors in 
  a Hermitian space $V$. The \textit{Gram matrix} of
  $\mathcal{L}$ is the Hermitian symmetric matrix
  \[
  \gram(\mathcal{L})
  = \gram(u_1, \dotsb, u_m) := \llbracket \ip{u_i}{u_j} \rrbracket.
  \]
  When $V$ is identified with $\K^n$ via coordinates with respect to
  an ordered basis $\basis = (u_1, \dotsb, u_n)$, the Hermitian form
  on $V$ is defined using the matrix $M = \gram(\basis)$.  We say that
  $M$ and $V$ has signature $(a, b, c)$ if $M$ has $a$ positive
  eigenvalues, $b$ negative eigenvalues, and $c$ zero eigenvalues.  Note
  that $a + b + c = n$.
  \end{defn}
  \begin{defn}[Standard forms]\label{def:std-forms}
  For integers $a,b,c \geq 0$ with $a+b+c=n$, let $V=\K^{a,b,c}$
  denote the $n$-dimensional $\K$-vector space $V$ with the Hermitian
  form 
  \begin{equation*}
    \ip{x}{y} = \overline{x}_1 y_1 + \dotsb +
    \overline{x}_a y_a - \overline{x}_{a + 1} y_{a + 1} - \dotsb -
    \overline{x}_{a+b} y_{a+b}.
  \end{equation*}
  The space $\K^{a,b,c}$ with this form is the \emph{standard
  Hermitian space of signature $(a,b,c)$}. Any finite dimensional
  Hermitian space $V$ is isomorphic to $\K^{a, b, c}$ for some tuple
  $(a, b, c)$ that is uniquely determined by $V$.  A Hermitian space
  $V \simeq \K^{a,b,c}$ is \emph{non-degenerate} when $c=0$,
  \emph{positive semi-definite} when $b=0$, \emph{positive
  definite} when $b=c=0$ and \emph{Lorentzian} when $b=1$ and $c=0$.
  Trailing $0$'s are removed from the notation.  Thus $\K^{a,b}$ is
  non-degenerate, $\K^a$ is positive definite, and $\K^{a,1}$ is
  Lorentzian.
 \end{defn}

\begin{defn}[Radical] \label{def:orthogonal}
  Let $V$ be a Hermitian space. Given $x \in V$ (resp. $U \subseteq
  V$), let $x^{\bot}$ (resp. $U^{\bot}$) denote the orthogonal
  complement of $x$ (resp. $U$).  The \emph{radical} of $V$ is the
  subspace $\rad(V) = V^{\bot}$. Note that $V$ is
  \emph{non-degenerate} when $\rad(V)$ is trivial.  If $M$ is the
  Hermitian matrix defining the form on $V$ with respect to some
  ordered basis then note that $\rad{V}$ is the kernel of the linear
  map $x \mapsto M x$.  A Hermitian space $V$ is \emph{degenerate}
  when $V \simeq \K^{a, b, c}$ with $c > 0$ or, equivalently, when
  $\rad(V)$ is not trivial and has dimension $c > 0$. In this situation
  the Hermitian form on $V$ naturally induces a nondegenerate form on
  $V/\rad (V) \simeq K^{a, b}$ such that the quotient map $V \to
  V/\rad(V)$ is form preserving.  Let $\aut(V)$ be the automorphism
  group of $V$.  For a nondegenerate Hermitian space $V \simeq \K^{a,
    b}$, this is the (generalized) \emph{orthogonal group} $O(a,b)$
  when $\K=\R$ and the (generalized) \emph{unitary group} $U(a,b)$
  when $\K=\C$.
\end{defn}

The following observation is elementary but crucial.

\begin{remark}[Negative Norm]\label{rem:negative-norm}
  Let $U$ be a subspace of a Lorentzian space with $u \in U$ of
  negative norm.  If $z \in \rad(U)$, then $z$ has norm $0$ and it
  belongs to the positive definite subspace $u^{\bot}$, so $z = 0$.
  In other words, \emph{any subspace of a Lorentzian space that
  contains a negative norm vector is non-degenerate.}
\end{remark}

\begin{defn}[Projective space]\label{def:proj-space}
  Let $V$ be a $\K$-vector space. The \emph{projective space} of $V$
  is the set $\proj (V) = \proj V$ of $1$-dimensional subspaces (or
  \emph{lines}) in $V$.  If $\dim V = n+1$, then $\proj V \simeq
  \proj^n_\K$.  We write $\proj : V\setminus \{0\} \to \proj V$ for
  the projection map, and $\bx = \proj x = \K x$ for the projective
  point determined by a non-zero vector $x$.
  If $V \simeq \K^{a,b,c}$ and $\bx \in \proj V$, then note that there
  always exist a vector $x$ representing $\bx$ with norm $1$, $-1$, or
  $0$.
\end{defn}
  
\begin{defn}[Hyperbolic space]\label{def:hyperbolic} 
  Let $V \simeq \K^{n,1}$ be a Lorentzian space and let $V_<$ be the
  set of vectors in $V$ of negative norm.  We write $\ball V \simeq
  \hyp^n_\K$ for the open subset $\proj V_< \simeq \proj\K^{n,1}_<$ of
  $\proj V$ consisting of $1$-dimensional subspaces spanned by
  negative norm vectors.  The notation $\ball V$, or \emph{ball of
  $V$}, refers to the fact that $\ball V \simeq \hyp^n_\K$ is an open
  topological ball. The space $\hyp^n_\R$ is \emph{real hyperbolic
  space} of real dimension~$n$, and $\hyp^n_\C$ is \emph{complex
  hyperbolic space} of complex dimension~$n$.
  The space $\hyp^n_\K$ has an essentially unique $\aut(V)$-invariant
  metric.  The distance between $\bx,\by \in \ball V$ is
  \[
  d(\bx,\by) = \cosh^{-1} \abs{\ip{x}{y}}
  \]
  where $x$ and $y$ are vector representatives of $\bx$ and $\by$ of
  norm $-1$.  We have chosen the scaled
  distance function with the simplest formula.  The space $\hyp^n_\K$
  is a complete geodesic Riemannian manifold of 
  constant negative curvature when $\K=\R$, and variable negative
  curvature when $\K = \C$.
\end{defn}

\begin{defn}[Hyperplanes]
  Let $V \simeq \K^{n,1}$ be a Lorentzian space, let $\ba \in \proj V$
  be a positive projective point. Choose $a \in V$ such that $\ba = \K
  a$. The hyperplane $a^{\bot} \simeq \K^{n-1,1}$ restricts to
  $a^\bot_<$ in $V_<$ and this projects to a totally geodesic
  hyperplane $\proj a^\bot_< \simeq H^{n-1}_\K$ inside $\ball V \simeq
  \hyp^n_\K$. We write $H_\ba = \proj a^\bot_<$, or simply $a^\bot =
  \proj a^\bot_<$ when there is no chance of confusion. Let $\hat{a} =
  a/\sqrt{a^2}$, so $\hat{a}^2 = 1$. The map $V \to V$ that sends $u
  \mapsto \hat{a} \ip{\hat{a} }{u}$ is \emph{orthogonal projection}
  onto the line $\ba = \K a$ since multiples of $a$ are fixed and
  vectors orthogonal to $a$ are in the kernel. So $u \mapsto u -
  \hat{a} \ip{\hat{a}}{u}$ is orthogonal projection onto
  $a^\bot$.
\end{defn}

\begin{defn}[Reflections and Mirrors]\label{def:refl} 
  Let $V = \K^{n,1}$ be a Lorentzian space, let $a \in V$ be a vector
  of positive norm, and let $\zeta$ be a root of unity in
  $\K$. Because the form on $V$ is non-degenerate and $a^2 \neq 0$,
  there is an orthogonal line-hyperplane decomposition $V = \K a
  \oplus a^\bot$, and an automorphism $R=R_{a,\zeta}$ of $V$ called a
  \emph{$\zeta$-reflection} which maps $u \mapsto \zeta u$ for all $u$
  in $\ba = \K a$ and point-wise fixes the $\K$-hyperplane $a^\bot$
  called its \emph{mirror}. Letting $\hat{a} = a/\sqrt{a^2}$, the
  formula for $R$ is \[ R(u) = (u - \hat{a}\ip{\hat{a}}{u}) + \zeta
  \hat{a} \ip{\hat{a}}{u} \] for all $u \in V$. When $\K=\R$, the only
  available root of unity is $\zeta = -1$, and $R = R_a$ is the real
  reflection that sends $a$ to $-a$ and fixes $a^\bot$. The
  reflections $R_{a, \zeta}$ act on the hyperbolic space $\ball V$
  point-wise fixing the totally geodesic hyperplane (determined by)
  $a^{\bot}$, called its \emph{mirror}. These reflections generate the
  group of isometries of the hyperbolic space $\ball V$. This group is
  isomorphic to $\proj O(n,1)$ if $\K = \R$ and is isomorphic to
  $\proj U(n, 1)$ if $\K = \C$.
\end{defn}

\begin{rem}[Distance to a hyperplane]\label{rem:dist-hyperplane}
  The distance between a (negative) point $\bx \in \hyp^n_\K$
  (represented by a vector $x \in \K^{n,1}$ of norm $-1$) and a
  hyperplane $H_\ba$ determined by a (positive) point $\ba$
  (represented by a vector $a \in \K^{n,1}$ of norm $1$) is given by
  the formula
  \[
  d(\bx,H_\ba) = \sinh^{-1} \abs{\ip{x}{a}}.
  \]
\end{rem}

\section{Lattices}\label{sec:lattices}

An integral \emph{$\Z$-lattice} is a free abelian group $\Z^n$ with a
non-degenerate symmetric bilinear form $\Z^n \times \Z^n \to \Z$.
This section reviews basic facts about $\Z$-lattices and
$\cO$-lattices where $\cO$ is the ring of integers in an algebraic
number field. For simplicity we almost always assume that $\cO$ is a
PID; the only exception being Example~\ref{ex:proj-plane-lattices}.

\begin{notation}[$\cO$ and $\K$]\label{not:cO-F-K}
  Given finitely many algebraic integers $r_1, r_2, \dotsb \in \C$,
  let $\Q(r_1, r_2, \dotsb)$ be the number field generated by $r_1,
  r_2, \dotsb$ and let $\cO = \cO(r_1, r_2, \dotsb)$ be its ring of
  integers.  Set $\K = \R$ when each $r_i \in \R$, and set $\K = \C$
  otherwise.  We always assume $\cO$ is stable under complex
  conjugation.  We write $\zeta$ for a root of unity, with specific
  values $\zeta_n = e^{2 \pi i/n}$, $\omega = \zeta_3$, and $i =
  \zeta_4$.  In our examples $\cO$ is often $\Z$ or the Eisenstein
  integers $\eisen=\Z[\omega]$ or the Gaussian integers
  $\gauss=\Z[i]$.
 \end{notation}

\begin{defn}[Extending Scalars]\label{def:extend}
  Let $L$ be a free $\cO$-module of finite rank.  Given any
  commutative ring $R \supseteq \cO$, we write
  \[
  R L = R  \otimes_{\cO} L,
  \] 
  and $R C$ for the $R$-span of $C$ in $R L$ for any subset $C \subseteq
  L$. The natural map $L \to RL$ is injective and we usually identify
  $L$ inside $RL$ using this map.  Note that if $L \simeq \cO^n$, then
  $R L \simeq R^n$.  When $R$ is a field such as $\K$, we say that $\K
  \otimes_{\cO} L$ is the \emph{underlying $\K$-vector space} of $L$.
\end{defn}

\begin{defn}[Lattices]\label{def:lattice}
  A \emph{Hermitian $\cO$-module} $L$ is a free $\cO$-module of finite
  rank, together with a Hermitian form $\ip{\;}{\;}: L \times L \to
  \Q(r)$.  A map of Hermitian $\cO$-modules means an inner product
  preserving $\cO$-module map and an \textit{isometry} means a
  bijective map of Hermitian $\cO$-modules.  We say $L$ is
  \emph{integral} if the Hermitian form takes values in $\cO$.  Given
  $S \subseteq L$, let $S^{\bot}$ be its orthogonal complement.  The
  \emph{radical} of $L$ is $\op{rad}(L) := L^{\bot}$.  We view $L$ as
  a subset of the vector space $V = \K L$.  The Hermitian form on $L$
  naturally extends to define a Hermitian form on $V$.  The
  \emph{signature} of $L$ is the signature of $V$, say $(a,b,c)$. We
  say that $L$ is \emph{degenerate} when $V$ is degenerate ($c > 0$)
  and $L$ is \emph{non-degenerate} when $V$ is non-degenerate ($c=0$).
  A \emph{$\cO$-lattice} means a nondegenerate Hermitian $\cO$-module.
  A map of $\cO$-lattices is just a map of Hermitian $\cO$-modules.
  Degenerate Hermitian $\cO$-modules as sometimes referred to as
  degenerate $\cO$-lattices.
\end{defn}

\begin{rem}[Sublattices]\label{rem:sublattices}
  Let $L$ be an $\cO$-lattice and let $L'$ be a free $\cO$-submodule
  of $L$ together with the form of $L$ restricted to $L'$.  Even
  though the form on $L$ is non-degenerate, the restricted form on
  $L'$ may be degenerate. We only say $L'$ is a \emph{sublattice of
  $L$} when the restricted form remains non-degenerate.
\end{rem}

\begin{defn}[Automorphisms]\label{def:aut}
  Let $L$ be an $\cO$-lattice and let $V = \K L \simeq \K^{a,b}$ be
  the corresponding Hermitian space.  The \emph{automorphism group of
  $L$} is the group $\aut(L) \subset \aut(V)$ of all symmetries of $V$
  (Definition~\ref{def:orthogonal}) that send $L$ to $L$.  When $V
  \simeq \K^{g,1}$ is Lorentzian, the \emph{ball of $L$} is the ball
  of $V$, a real or complex hyperbolic space.
  (Definition~\ref{def:hyperbolic}) In symbols $\ball L = \ball V
  \simeq \hyp^g_\K$. Since the action of $\aut(L)$ on $L$ preserves
  norms and sends lines to lines, $\aut(L)$ also acts on $\ball L
  \simeq \hyp^g_\K$.
\end{defn}

We often assume $\cO$ is a PID since this eliminates
technicalities. For example, when $\cO$ is a PID, every Hermitian
$\cO$-module has a canonical lattice quotient.

\begin{prop}[Lattice Quotients]\label{prop:degenerate}
  Let $\cO$ be a PID and let $L$ be a Hermitian $\cO$-module.
  \begin{enumerate}
  \item $L$ contains a Hermitian $\cO$-submodule $\bar{L}$ such
    that $L = \op{rad}(L) \oplus \bar{L}$ and the Hermitian form on
    $\bar{L}$ is nondegenerate; i.e. $\bar{L}$ is a lattice.
  \item The quotient map $L \to L/\op{rad}(L)$ induces an isomorphism
    $\bar{L} \simeq L/\op{rad}(L)$. In particular, $L/\rad(L)$ is a
    $\cO$-lattice, uniquely determined up to isometry.
   \end{enumerate}
\end{prop}

\begin{proof}
  Part (1) is Proposition~6.2.26 in \cite{AdkinsWeintraub92}.  As a
  result, the quotient map $L \to L/\op{rad}(L)$ induces a
  $\cO$-linear inner product preserving bijection $\bar{L} \to
  L/\op{rad}(L)$.  Since $\bar{L}$ is a lattice, so is
  $L/\op{rad}(L)$.
\end{proof}

We record a few useful consequences. The first one is nearly
immediate.

\begin{prop}[Onto Maps]\label{prop:maps-onto}
  Let $\cO$ be a PID. If $T: L \to L'$ is an onto map of Hermitian
  $\cO$-modules and $L'$ is an $\cO$-lattice, then $\ker(T) = \rad(L)$
  and the induced map $\overline{T} \colon L/\rad(L) \to L'$ is an
  isometry of $\cO$-lattices.
\end{prop}

\begin{prop}[Spanning Sets]\label{prop:gram-spanning}
  Let $L_1, L_2$ be $\cO$-lattices with finite ordered spanning sets
  $\spanning_1 = \{x_1,\cdots,x_n\}$ and $\spanning_2 =
  \{y_1,\cdots,y_n\}$.  If $\gram(\spanning_1) = \gram(\spanning_2)$
  then there is a $\cO$-lattice isometry $L _1\to L_2$ sending $x_i$
  to $y_i$ for all $i$.
\end{prop}

\begin{proof}
  Let $L = \cO^n$ be the Hermitian $\cO$-module with basis $\basis =
  \{e_1,\cdots,e_n\}$ and Hermitian form defined by the matrix
  $\gram(\spanning_1) = \gram(\spanning_2)$.  Let $T_j \colon L \to
  L_j$ be the form preserving surjection of $\cO$-modules sending
  $e_i$ to $x_i$ (if $j = 1$) and $e_i$ to $y_i$ (if $j= 2$).  By
  Proposition~\ref{prop:maps-onto}, the $T_j$'s induce isometries
  $\overline{T}_j \colon L/\rad(L) \to L_j$.  The composition
  $\overline{T}_2 \circ (\overline{T}_1)^{-1} \colon L_1 \to L_2$ is
  the desired isometry.
\end{proof}

The \emph{symmetry} of a subset $\spanning \subset L$ is a permutation
of $\spanning$ that preserves the inner product.  The following is an
immediate corollary of Proposition~\ref{prop:gram-spanning}.

\begin{cor}[Symmetries]\label{cor:symmetries}
  If $L$ is an $\cO$-lattice with finite spanning set $\spanning$,
  then the symmetries of $\spanning$ are automorphisms of $L$.
\end{cor}

Next recall the standard self-dual even and odd $\Z$-lattices.

\begin{defn}[Self-Dual]\label{def:self-dual}
  The \textit{dual} of an $\cO$-lattice $L$ is 
  \begin{equation*}
    L^\dual := \{ v \in \K L \colon \ip{v}{L} \subseteq \cO \}. 
  \end{equation*}
  Notice that a lattice $L$ is \emph{integral} if and only if
  $\ip{L}{L} \subseteq \cO$ if and only if $L \subseteq L^\dual$.  A
  lattice is \emph{self-dual} when $L = L^\dual$.  Self-dual lattices
  are also known as \emph{unimodular} lattices because their
  discriminant ideal is the unit ideal.  The \emph{discriminant ideal}
  of an $\cO$-lattice $L$ is the principal ideal $(d)$ where $d =
  \det(\gram(\basis))$ for an $\cO$-basis $\basis$ of $L$.  When $\cO
  = \Z$, $\gauss=\Z[i]$, or $\eisen=\Z[\omega]$, all units in $\cO$
  are roots of unity and the \emph{discriminant} is the positive
  integer $\size{d}$.
\end{defn}

\begin{defn}[A Few Important Lattices]\label{def:important-lattices}
  An integral $\Z$-lattice is \emph{even} when all lattice vectors
  have even norm, and it is \emph{odd} if it is not even.  The
  \emph{standard odd lattice} $I_{a,b} := \Z^{a,b} \subset \R^{a,b}$
  is the set of all vectors with integer coordinates inside the
  standard space $\R^{a,b}$.  This $I_{a, b}$ is an odd self-dual
  lattice of signature $(a, b)$ and is the unique such lattice if both
  $a, b$ are strictly positive.  \par Even self dual $\Z$-lattices of
  signature $(a, b)$ exists if and only if $a \equiv b \bmod 8$;
  assume this.  The \emph{standard even lattice} $\II_{a,b}$ is the
  set of all $(x_1, \dotsb, x_{a+b}) \in \R^{a,b}$ such that $x_i \in
  \Z$ for all $i$ or $(x_i + \tfrac{1}{2}) \in \Z$ for all $i$, with
  the restriction that $\sum x_i \equiv 0 \bmod 2$.  This $\II_{a, b}$
  is an even self-dual lattice of signature $(a, b)$ and is the unique
  such lattice if if both $a, b$ are strictly positive.  The lattice
  $\II_{8,0}$ is the $E_8$ root lattice in $\R^8$.  In
  \S\ref{sec:ball-quoients} we use a rank~$4$ Hermitian
  $\eisen$-lattice $E_8^{\eisen} = (E_8^{\eisen}, \ip{\;}{\;})$ whose
  underlying $\Z$-module with the bilinear form
  $\tfrac{2}{3}\op{Re}(\ip{x}{y})$ is $E_8$; we call $E_8^{\eisen}$
  the Eisenstein $E_8$ lattice.  We also use the complex and
  quaternionic forms of some other lattices like $D_4$ or the Leech
  lattice which are defined in similar manner.
\end{defn}

\begin{defn}[Hyperbolic Cells]\label{def:hyp-cell}
  For each $c \in \cO - \lbrace 0 \rbrace$, the \emph{hyperbolic cell}
  $\II^{\cO,c}_{1,1}$ is the rank two $\cO$-lattice of signature
  $(1,1)$ with Hermitian form defined by the matrix
  \[
  M = \begin{bmatrix} 0 & c \\ \bar{c} & 0 \end{bmatrix} 
  \]
  One can show that $\II^{\Z, 1}_{1,1}$ is isometric to $\II_{1,1}$
  and if $k > 0$, then $\II_{8 m + k, k}$ is isometric to the
  orthogonal direct sum of $m$ copies of $E_8 = \II_{8,0}$ and $k$
  copies of $\II_{1,1}$.
\end{defn}

\section{Lattices from Graphs}\label{sec:graph-lattice}

This section recalls and extends the main construction in
\cite{NeumaierSeidel83}.  We highlight its functorial nature
(Theorem~\ref{thm:latt-functor}) and the failure of functoriality when
one tries to extend the functor beyond the Lorentzian case
(Example~\ref{ex:affine-fail}).

\begin{defn}[Graphs]\label{def:graphs}
  In this article, $\Delta$ is a simplicial graph with \emph{vertices}
  $\ver_\Delta$ and \emph{edges} $\edg_\Delta$. We assume
  $\ver_\Delta$ has a fixed linear order, i.e. a bijection
  $\ver_\Delta \simeq \{1,2,\cdots,v\}$, where $v =
  \size{\ver_\Delta}$.  When there is no chance of confusion, we write
  $\Delta$ for $\ver_\Delta$. Edges are labeled by their unordered
  pair of (distinct) endpoints.  The \emph{adjacency matrix of
  $\Delta$} is a $v\times v$ matrix $A = A_\Delta = \llbracket a_{i j}
  \rrbracket$ with $a_{i j} = 1$ when $\{i,j\}$ is edge of $\Delta$
  and $a_{i j} = 0$ otherwise.  The (adjacency) \emph{spectrum} of
  $\Delta$ is the multiset of eigenvalues of $A_\Delta$.  Since
  $A_{\Delta}$ is real and symmetric, all of its eigenvalues are real.
  We write $\spec(\Delta)$ or $\spec(A_\Delta)$ for the multiset
  $\{\lambda_1^{m_1}\lambda_2^{m_2} \cdots \lambda_\ell^{m_\ell}\}$
  where $\lambda_1 > \lambda_2 \cdots > \lambda_\ell$ are the distinct
  eigenvalues of $A_\Delta$ and $m_i$ is the multiplicity of
  $\lambda_i$.  The largest eigenvalue $\lambda_1$ is the
  \emph{principal eigenvalue}. When $\Delta$ is connected, $\lambda_1$
  has multiplicity $1$ and it has a corresponding \emph{principal
  eigenvector} with all positive coordinates.  When $\Delta$ is
  $k$-regular, $\lambda_1=k$ and $\one = (1,1,\ldots,1)^\tr$ is a
  $\lambda_1$-eigenvector.
\end{defn}

\begin{defn}[Modules with Forms from Graphs]\label{def:shift-module}
  Let $\Delta$ be a graph with an ordering on its $v$ vertices and with
  adjacency matrix $A$.  If $a$ and $b$ are real constants and $b \neq
  0$, note that the matrix $aI + bA$ has eigenvalues $\mu_i = a +
  b\lambda_i$ with the same multiplicities and eigenvectors.  In
  particular, the matrix $r I - A$ with $r
  \in \R$ has spectrum $\{(r - \lambda_\ell)^{m_\ell} \cdots (r -
  \lambda_2)^{m_2} (r - \lambda_1)^{m_1} \}$.

  Let $R$ be a commutative ring, let $R^{\Delta} $ denote the free
  $R$-module of the set of all functions from $\Delta$ to $R$, and fix $r
  \in R$.  Use the ordering on the vertices of $\Delta$ to identify
  $R^{\Delta} = R^v$.  Let $R^{\Delta}_r$ denote the free module
  $R^{\Delta}$ together with the $R$-valued symmetric bilinear form
  defined by
  \[
  \ip{u}{w} = u^T (r I - A) w \text{\; for all \;} u, w \in R^{\Delta}.
 \]
 Let
  \[
 R_{\Delta,r} = R^{\Delta}_r/ \op{rad}( R^{\Delta}_r)
 \]
 with the nondegenerate bilinear form induced from $R^{\Delta}_r$.  In
 the general notation of \ref{definition-Hermitian-R-modules} (with
 the involution $\sigma_R = \op{id}_R$) one has
 \[
 R^{\Delta}_r = R^v_B 
 \text{\; and \;}
 R_{\Delta,r} = R_{v,B}
 \text{\; where \;} B = (r I - A).
 \]
 Let $\lbrace e_i \colon i \in \Delta \rbrace$ be the standard basis
 of the free module $R^{\Delta}_r$ and let
 \[
 \spanning_{\Delta, r} = 
 \lbrace \alpha_i \colon i \in \Delta \rbrace
 \]
 be the spanning set in $R_{\Delta, r}$ formed by the images of these
 basis vectors\footnote{We use $\alpha_i$ for the spanning vectors in
 $\spanning_{\Delta, r}$ to echo the way Bourbaki used $\alpha_i$ for
 the simple roots of classical root lattices.  The symbol $\spanning$
 in $\spanning_{\Delta, r}$ is meant to remind the reader of spanning sets
 and simple roots.}.  For future calculations, we reiterate that the
 dot products among the basis vectors $\lbrace e_i \colon i \in \Delta
 \rbrace$ and the spanning vectors $\lbrace \alpha_i \colon i \in
 \Delta \rbrace$ are given by
 \begin{equation}\label{eq:rI-A}
  \ip{e_i}{e_j} = 
  \ip{\alpha_i}{\alpha_j}
  =
  \begin{cases}
  \; \; r & \text{\; if $i = j$,} \\
  -1 & \text{\; if $\lbrace i, j \rbrace$ is an edge of $\Delta$,}\\
  \; \; 0 & \text{\; otherwise.}
  \end{cases}
 \end{equation}

\end{defn}

\begin{remark}\label{remark-real-change-of-rings}
  Recall from Lemma \ref{lemma-change-of-rings} that if $R \subseteq
  Q$ is an extension of commutative rings and $r \in R$, then we have
  natural maps
 \[
 \begin{tikzcd}[column sep=large]
   R^{\Delta}_r
   \arrow[r, ->>]
   \arrow[d, hook] &
   R_{\Delta, r}
   \arrow[d, hook]
   \\
   Q^{\Delta}_r
   \arrow[r, ->>] &
   Q_{\Delta, r}.
 \end{tikzcd}
 \]
 and these maps are compatible with triples of extensions. Using this,
 we'll identify $R^{\Delta}_r \subseteq Q^{\Delta}_r$ and $R_{\Delta,
   r} \subseteq Q_{\Delta, r}$ whenever we have commutative ring
 extension $R \subseteq Q$.  Note that the standard basis $\lbrace e_i
 \rbrace$ of $R^{\Delta}_r$ (resp. the standard spanning set $\lbrace
 \alpha_i \rbrace$ of $R_{\Delta, r}$) are identified with the
 corresponding vectors with the same names in $Q^{\Delta}_r$
 (resp. $Q_{\Delta, r}$).
\end{remark}

\begin{definition}[Important Special Case]\label{def-R_Delta}
  We are specifically interested in the situation where $\Delta$ is a
  connected incomplete graph and $r = \lambda_2$ is the second largest
  eigenvalue of $A$.  So for a commutative ring $R$ containing
  $\lambda_2$, we write
  \[
  R_{\Delta} = R_{\Delta, \lambda_2} \text{\; and \;}
  \spanning_{\Delta} = \spanning_{\Delta, \lambda_2}.
  \]
\end{definition}

\begin{remark}[Rank and Signature]\label{rem:rank-signature}
  Assume $R$ is a PID and $\Z[r] \subseteq R$.  Recall from
  \ref{lemma-Lorentzian-rep-is-image-of-Gram} that $R_{\Delta, r}$ is
  then a free $R$-module whose $R$-rank is the rank of the matrix $(r
  I - A)$. So
  \[
  \op{rank}_R( R_{\Delta, r}) = 
  \begin{cases} 
    v & \text{\; if $r$ is not an eigenvalue of $A$,} \\
    v - m_i  & \text{\; if $r = \lambda_i$ is an eigenvalue of $A$
      with multiplicity $m_i$.} 
  \end{cases}
  \]
  In the second case, it is assumed that $\Z[\lambda_i] \subseteq R$.
  Specializing to $R = \R$, this gives us the dimensions of the vector
  spaces $\R_{\Delta, r}$.  The form on $\R^{\Delta}_r$ has signature
  $(a, b, c)$ where $a$ is the sum of the multiplicities of the
  $A_\Delta$-eigenvalues strictly \emph{less} than $r$ and $b$ is the
  sum of the multiplicities of the $A_\Delta$-eigenvalues strictly
  \emph{greater} than $r$ and $c = v - a - b$.  So the nondegenerate
  form on $\R_{\Delta, r}$ has signature $(a, b)$.
\end{remark}

\begin{definition}\label{def-O_r}
  Let $r \in \R$ be an algebraic integer.  Let $\cO=\cO(r)$ be the
  ring of integers of the number field $\Q(r)$.  Often we work over $R
  = \cO$ or a similar ring.  For simplicity, we assume that $\cO$ is a
  PID unless otherwise stated.
\end{definition}

\begin{lemma}\label{lemma-O_Delta-in-R_Delta}
  Let $r \in \R$ and $\cO=\cO(r)$ be as in Definition~\ref{def-O_r}. Then
  $\cO_{\Delta, r}$ is the $\cO$-span of a $\R$-basis of the vector
  space $\R_{\Delta, r}$. Thus 
  $\R_{\Delta, r} = \R \otimes_{\cO} \cO_{\Delta, r}$.
\end{lemma}

\begin{proof}
  This follows from Lemma \ref{lemma-base-change} in the
  Appendix. Here is another argument.  Let $F$ be the fraction field
  of $\cO$. Identify $\cO_{\Delta, r} \subseteq F_{\Delta, r}
  \subseteq \R_{\Delta, r}$.  Let $u_1, \dotsb, u_n$ be a $\cO$-basis
  of $\cO_{\Delta, r}$.  Then $u_1, \dotsb, u_n$ are still linearly
  independent in $F_{\Delta, r}$ because a dependence relation over
  $F$ would give a dependence relation over $\cO$ by clearing
  denominators.  As we already noted
  \[
  \op{rank}_{\cO}( \cO_{\Delta, r} ) 
  =  \op{dim}_{F}( F_{\Delta, r} )  
  = \op{dim}_{\R}( \R_{\Delta, r} ).
  \]
  It follows that $u_1,\dotsc,u_n$ is an $F$-basis of $F_{\Delta, r}$
  and furthermore that $u_1, \dotsb, u_n$ form an $\R$-basis of
  $\R_{\Delta, r}$\footnote{The natural inclusion $j:F_{\Delta, r}
  \hookrightarrow \R_{\Delta, r}$ induces an $\R$-linear map $ J:\R
  \otimes_F F_{\Delta, r}\longrightarrow \R_{\Delta, r}, $ given by
  $a\otimes x\mapsto a\,j(x)$. Since $u_1,\dots,u_n$ is an $F$-basis
  of $F_{\Delta, r}$, the vectors $1\otimes u_1, \dotsb, 1\otimes u_n$
  form an $\R$-basis of $\R \otimes_F F_{\Delta, r}$
  \cite[p.308]{Roman08}.  Moreover, $J$ is onto, since the vectors
  $J(1 \otimes \alpha_i) = j(\alpha_i)$ span $\R_{\Delta, r}$.  Both
  source and target of $J$ have $\R$-dimension $n$, so $J$ is an
  isomorphism. Therefore $j(u_1),\dotsb, j(u_n)$ is an $\R$-basis of
  $\R_{\Delta, r}$.}.  So $\cO_{\Delta, r}$ is the $\cO$-span of a
  basis of $\R_{\Delta, r}$.
\end{proof}

It is sometimes convenient to rescale the form used to define
$\R_{\Delta, r}$ and $\Z_{\Delta, r}$.  For the sake of simplicity,
the proofs in this section use the original scaling, but normalized
version plays a key role in the complex variation 
(Definition~\ref{def:cplx-construction}).

\begin{example}[ADE root lattices]\label{ex:ade}
  The only connected graphs $\Delta$ with $\lambda_1 = 2$ are the
  extended simply-laced Dynkin diagrams: $\widetilde{A}_n$ ($n\geq
  2$), $\widetilde{D}_n$ ($n\geq 4$), and $\widetilde{E}_n$
  ($n=6,7,8$).  Note that $v = \size{\Delta} = n+1$. If $\Delta$ is
  one of these, then $\Z^\Delta_2$ is a positive semidefinite
  $\Z$-lattice in $\R^\Delta_2 \simeq \R^{n,0,1}$ with a
  $1$-dimensional radical and $\Z_{\Delta,2}$ is the corresponding
  positive definite ADE root lattice in $\R_{\Delta,2} \simeq \R^n$.
  The spanning set $\spanning_{\Delta, 2}$, in Bourbaki notation, is
  $\{\alpha_1,\ldots,\alpha_n,\widetilde\alpha\}$ where $\alpha_1,
  \dotsb, \alpha_n$ is a set of simple roots (indexed by the vertices
  of the corresponding unextended diagram) and $\widetilde{\alpha}$ is
  the lowest root.  The subset of $\spanning_{\Delta, 2}$ with
  $\widetilde \alpha$ removed is a basis of the lattice $\Z_{\Delta,
    2}$, and these vectors are indexed by the vertices of the
  corresponding (non-extended) Dynkin diagram.
\end{example}

In our context the map from $\lbrace e_i \colon i \in \Delta \rbrace$
to $\spanning_{\Delta, r}$ is almost always injective.

\begin{lemma}[Distinct]\label{lem:distinct}
  Assume $r \in R \subseteq \R$.  When $r$ is positive the projection
  $R^{\Delta}_r \to R_{\Delta, r}$ yields a bijection $\lbrace e_i
  \colon i \in \Delta \rbrace \to \spanning_{\Delta, r}$.  If, in
  addition, $\Delta$ is connected with $\size{\Delta}>2$, then the
  $\alpha_i \in \spanning_{\Delta, r}$ are also projectively distinct
  with distinct hyperplanes $\alpha_i^\bot$.
\end{lemma}

\begin{proof}
  If $\alpha_i = \lambda \alpha_j$ for some $i \neq j$, then $r =
  \ip{\alpha_j}{\alpha_j} = \lambda^2 \ip{\alpha_i}{\alpha_i} =
  \lambda^2 r$, so $\lambda$ is $\pm 1$.  Next $r > 0$ and $0 \geq
  \ip{\alpha_i}{\alpha_j} = \lambda r$, so $\lambda \neq 1$.  This
  already shows that the map $e_i \mapsto \alpha_i$ is $1$-to-$1$ and
  the map is onto by construction.  Suppose $\lambda= -1$. When
  $\abs{\Delta} > 2$, then there is a vertex $k \in \Delta - \{ i , j
  \}$ that has an edge to $i$ or $j$.  Without loss, suppose $\{k , i
  \}$ is an edge.  Then $ 0 > \ip{\alpha_j}{\alpha_k} =
  -\ip{\alpha_i}{\alpha_k} = 1$ is a contradiction.  This rules out
  $\lambda = -1$, So $\alpha_i$ and $\alpha_j$ are linearly
  independent.
\end{proof}

The restriction $r=\lambda_2>0$ in Lemma~\ref{lem:distinct} only
excludes a few graphs, but these graphs show that the restriction is
necessary.

\begin{rem}[Indistinct]\label{rem:indistinct}
  The \emph{complete graph} $\Delta = K_{n+1}$ on $n+1$ vertices has
  $\spec(\Delta)=\{n^1 \; (-1)^n\}$ and $\lambda_2=-1$.  The space
  $\Z^{\Delta}_{\lambda_2}$ has an $n$-dimensional radical, $\Z_\Delta$ is only
  $1$-dimensional, and, in fact, all the $n+1$ basis vectors $\lbrace
  e_i \colon i \in \Delta \rbrace \subseteq \Z^{\Delta}$ have the same
  image in $\Z_\Delta$.  Fortunately, complete graphs are the only
  connected graphs with $\lambda_2 <0$. There is a similar collapse of
  the spanning set $\spanning_\Delta$ when $\Delta$ is a complete
  multipartite graph.  These graphs have $\lambda_2=0$ and they are
  the only connected graphs with $\lambda_2=0$.
\end{rem}

With the restriction $\lambda_2>0$, the resulting spaces and lattices
are Lorentzian.

\begin{prop}[Lorentzian spaces and lattices]\label{prop:lor-spaces}
  Let $\Delta$ be a connected graph with $r = \lambda_2>0$. Then the
  vector space $\R_\Delta$ is Lorentzian of signature $(g,1)$ where $g
  = v-m_2-1 >0$.  Let $\cO=\cO(r)$ be the ring of integers of $\Q(r)$.
  Assume $\cO$ is a PID.  Then $\cO_{\Delta}$ is a Lorentzian
  $\cO$-lattice of signature $(g, 1)$ inside $\R_{\Delta}$, i.e.  it
  is $\cO$-span of an $\R$-basis of $\R_{\Delta}$.
\end{prop}

\begin{proof}
  The signature of $\R_{\Delta}$ is computed in
  Remark~\ref{rem:rank-signature} and $\cO_{\Delta}$ is $\cO$-span of
  an $\R$ basis of $\R_{\Delta}$ by
  Lemma~\ref{lemma-O_Delta-in-R_Delta}, so $\R \otimes_{\cO}
  \cO_{\Delta} \simeq \R_{\Delta}$ and hence $\cO_{\Delta}$ is a
  $\cO$-lattice of the same signature.
\end{proof}
 
\begin{defn}[Lorentzian Representation]\label{def:lor-rep}
  Let $\Delta$ be a connected graph with $r=\lambda_2>0$ and assume
  that the ring of integers of $\Q(r)$, denoted by $\cO$ is a PID. The
  \emph{Lorentzian representation} of $\Delta$ is the pair
  $(\cO_\Delta,\spanning_\Delta)$, the Lorentzian $\cO$-lattice
  $\cO_\Delta$ together with its spanning set $\spanning_\Delta =
  \{\alpha_i\}$ indexed by the vertices of $\Delta$.
\end{defn}

\begin{defn}[Categories]\label{def:categories}
  Let $\graph$ be the category of connected graphs whose arrows are
  inclusions of connected induced subgraphs.  For each positive
  algebraic real $r$ such that the ring of integers $\cO$ of $\Q(r)$
  is a PID, let $\graph(r)$ be the full subcategory of $\graph$ where
  we restrict to graphs where $\lambda_2=r$.  Let $\cO$-$\lattice$ be
  the category of Lorentzian $\cO$-lattices and $\cO$-lattice maps.
\end{defn}

\begin{theorem}[Lorentzian lattice functor]\label{thm:latt-functor}
  Let $r > 0$ be a real algebraic integer such that $\cO$ is a PID.
  Then there is a functor $ \graph(r) \to \cO$-$\lattice$ sending
  $\Delta \mapsto \cO_\Delta$.  In other words, (a) if $\Delta$ is a
  connected graph with $\lambda_2=r$, then $\cO_\Delta$ is a
  Lorentzian $\cO$-lattice, and (b) if $\Delta' \subset \Delta$ is a
  connected induced subgraph with the same second eigenvalue
  $\lambda_2=r$, then the inclusion map $\Delta' \into \Delta$ induces
  an $\cO$-linear embedding $\cO_{\Delta'} \into \cO_\Delta$
  between Lorentzian $\cO$-lattices sending $\spanning_{\Delta'} \into
  \spanning_\Delta$.
\end{theorem}

\begin{proof}
  Part (a) is Proposition~\ref{prop:lor-spaces}.  For Part (b) let
  $\cO^{\Delta'}_r$ have basis $\{ e'_j \colon j \in \Delta'\}$ and let
  the quotient $\cO_{\Delta'}$ have spanning set $\{ \alpha'_j \colon
  j \in \Delta' \}$, where the primes distinguish these from the
  similarly indexed elements inside $\cO^\Delta_r$ and $\cO_\Delta$.
  Next, let $\mathfrak{E}$ be the $\cO$-span of the basis elements $\{
  e_i \colon i \in \Delta'\}$ inside $\cO^{\Delta}_r$ and let
  $\mathfrak{A}$ be the $\cO$-span of the spanning set elements
  $\{\alpha_i \colon i \in \Delta'\}$ inside $\cO_\Delta$. Note that
  the quotient map $q \colon \cO^{\Delta}_r \onto \cO_{\Delta}$
  restricted to $\mathfrak{E}$ maps onto $\mathfrak{A}$ with kernel
  $\mathfrak{E} \cap \rad(\cO^\Delta_r)$.
  Consider the commutative diagram
  \begin{equation*}
  \begin{tikzcd}
    \cO^{\Delta'}_r \arrow[r,leftarrow,"\iota"]
    \arrow[d,twoheadrightarrow,"q"'] & \mathfrak{E}
    \arrow[r,hookrightarrow] \arrow[d,twoheadrightarrow,"q"'] &
    \cO^\Delta_r \arrow[d,twoheadrightarrow,"q"']
    \arrow[r,hookrightarrow] & \R^\Delta_r
    \arrow[d,twoheadrightarrow,"q"']\\ \cO_{\Delta'}
    \arrow[r,twoheadleftarrow,"\bar{\iota}"] & \mathfrak{A}
    \arrow[r,hookrightarrow] & \cO_\Delta \arrow[r,hookrightarrow] &
    \R_\Delta
  \end{tikzcd}
  \end{equation*}
  The map $\iota \colon \mathfrak{E}
  \to \cO^{\Delta'}_r$ sending $e_j \mapsto e_j'$ is an isometry since 
  it takes a basis to a basis.
  Let $q \circ \iota \colon \mathfrak{E} \onto \cO_{\Delta'}$ be the
  form preserving $\cO$-linear map sending $e_j \mapsto e'_j \mapsto
  \alpha'_j$.  If $u \in \rad(\cO^\Delta_r) \cap \mathfrak{E}$, then
  $\ip{u}{e_{j}} = 0$ and $\ip{\iota (u)}{e'_j} = 0$ for all $j \in
  \Delta'$, so $q(\iota(u))=0$.  As a result, there is an induced
  form-preserving $\cO$-linear map $\bar{\iota} \colon \mathfrak{A}
  \onto \cO_{\Delta'}$ sending $\alpha_j \mapsto \alpha'_j$.
  Note that $\mathfrak{A}$ contains vectors of negative norm (since
  $\mathfrak{E}$ contains such vectors and the map to $\mathfrak{A}$
  is onto)
  and the $\R$-span of $\mathfrak{A}$ is a subspace of the Lorentzian
  space $\R_\Delta$.  By Remark~\ref{rem:negative-norm} the form on
  $\mathfrak{A}$ is non-degenerate and by
  Proposition~\ref{prop:maps-onto} the form-preserving onto map
  $\bar{\iota}$ is an isomorphism.  The inverse of $\bar{\iota}$
  composed with the inclusion $\mathfrak{A} \hookrightarrow
  \cO_{\Delta}$ gives the required map $\cO_{\Delta'} \hookrightarrow
  \cO_{\Delta}$.
\end{proof}

\begin{rem}[History]
  Proposition~\ref{prop:lor-spaces} and
  Theorem~\ref{thm:latt-functor} essentially goes back to
  \cite{NeumaierSeidel83}, although these authors only work over the
  integers.  We include a proof because the same argument goes though
  in the much less studied complex Hermitian case.  See
  Section~\ref{sec:bipartite}.  We should note that we have not seen
  the functorial aspect stated explicitly in the literature.
\end{rem}

\begin{rem}[Extending the Functor]\label{rem:extend-functor}
  One can extend Theorem~\ref{thm:latt-functor} to a larger set of
  graphs as follows.  Fix a positive algebraic integer $r \in \R$ such
  that $\cO(r)$ is a PID and let $\graph( \lambda_1 \neq r \geq
  \lambda_2)$ be the category of connected graphs with $ \lambda_1
  \neq r \geq \lambda_2$.  Then $\Delta \mapsto \cO_{\Delta, r}$ is a
  functor to $\cO$-lattices and lattice maps.  If $\lambda_1 > r \geq
  \lambda_2$, then $\cO_{\Delta,r}$ is still a Lorentzian lattice.
  But if $r > \lambda_1$, then $\cO_{\Delta,r}$ is a positive definite
  lattice.  In these additional cases where $r$ is not an eigenvalue
  of $\Delta$, the proof simplifies since the radical is trivial and
  $\cO^\Delta_r = \cO_{\Delta, r}$ (Remark~\ref{rem:rank-signature}).
\end{rem}

The restriction $r\neq \lambda_1$ is necessary because functoriality
can fail when $r=\lambda_1$.

\begin{example}[Failing Functoriality]\label{ex:affine-fail}
  Let $\Delta$ be the triangular graph $T(6)$. Its vertices are
  labeled by the fifteen $2$-elements subsets of $\{1,2,3,4,5,6\}$ and
  it has an edge connecting two vertices if and only if their labels
  are not disjoint.  The spectrum of $\Delta$ is $\{8^1 \; 2^5 \;
  (-2)^9\}$, with $\lambda_2 = 2$, and $\Z_\Delta = \Z_{\Delta, 2}$ is
  a Lorentzian lattice in $\R^{9,1}$ with a $15$ element spanning set
  $\spanning_\Delta$.  Let $\Delta'$ be one of the $80$ triangles in
  $\Delta$ and let $\spanning = \{\alpha_1,\alpha_2,\alpha_3\}$ be the
  vectors in $\spanning_\Delta$ indexed by its vertices.  The
  $\Z$-span of $\spanning$ is $\mathfrak{A} \subset \Z_\Delta$.
  The triangle $\Delta'$ has $\lambda_1=2$, its lattice $\Z_{\Delta',
    2}$ is the $A_2$ root lattice in $\R^2$, and the vectors in its
  spanning set $\{\alpha'_1,\alpha'_2,\alpha'_3\}$ have sum $0$. The
  vectors in $\spanning$, on the other hand, do not have sum $0$ since
  their sum has a non-zero inner product with any vector in
  $\spanning_\Delta$ indexed by a vertex in $\Delta\setminus \Delta'$
  connected to $\Delta'$ by one or more edges.  This shows that the
  function sending $\alpha'_i$ to $\alpha_i$ does not extend to an
  embedding of $\Z_{\Delta', 2}$ in $\Z_{\Delta, 2}$.
\end{example}

\section{Automorphism Groups}\label{sec:groups}

The Lorentzian representation of a graph leads to a group action on
hyperbolic space that preserves a configurations of points and
hyperplanes indexed by the vertices of $\Delta$. Automorphisms of one
of these Lorentzian lattices include symmetries of $\Delta$ and
reflections in its roots.  This section establishes the necessary
notation for these groups.  The definitions are given over the field
$\K$ where $\K=\R$ or $\C$.

\begin{defn}[Roots, Reflections and Mirrors]\label{def:aut-refls}
  Let $\cO$ be as in Definition~\ref{def-O_r}.
  Let $L$ be a Lorentzian $\cO$-lattice with vector space $V = \K L
  \simeq \K^{g,1}$. A positive norm vector $\alpha \in L$ is called a
  \emph{root} of $L$ there is a non-trivial root of unity $\zeta$ such
  that the reflection $R=R_{\alpha,\zeta}$ (Definition~\ref{def:refl})
  sends $L$ to $L$.  The hyperplane $\alpha^{\bot} \subset V = \K L$
  (or the hyperplane $\hyp_\alpha$ it determines in $\ball V \simeq
  \hyp^g_\K$), is called the \textit{mirror of the root} $\alpha$ or
  the \emph{mirror of the reflection} $R$.  The \textit{reflection
    group of $L$} is the group $\Refl(L) \subset \aut(L)$ generated by
  reflections in the roots of $L$.  For any set $\spanning$ of roots
  of $L$, say a spanning set, $\Refl(\spanning) \subset \Refl(L)$ is
  the group generated by reflections in the roots in $\spanning$. In
  many cases all of these inclusions are proper.
\end{defn}

Lattice vectors with very small integer norms are always roots.

\begin{defn}[Short Roots]\label{def:short-roots}
  Let $L$ be an integral $\cO$-lattice.  Lattice elements of norm~$1$
  and norm~$2$ are always roots of reflections of order~$2$ (with
  $\zeta_2 =-1$) since the formula for these reflections sends lattice
  elements to $\Z$-linear combinations of lattice elements
  (Definition~\ref{def:refl}).  Lattice vectors of norm~$1$ are called
  \emph{short roots}.  The mirror of a short root is a \emph{short
  mirror}.  For integral self-dual $\Z$-lattices, the only roots are
  those of norm $1$ and $2$.  The group $\Refl_{*}(L)$ is the subgroup
  of $\Refl(L)$ generated by reflections in the short roots of $L$.
\end{defn}

Recall that classical reflection groups such as Coxeter groups can be
defined by abstract group presentations, and that the presentation of
a simply-laced Coxeter group is completely determined by a simplicial
graph.

\begin{defn}[Coxeter, Artin and Shephard]\label{def:uniform}
  Let $\Delta$ be a graph with vertices $\ver_\Delta$ and edges
  $\edg_\Delta$, and let $p$ be a number in $\{2,3,\ldots\} \cup
  \{\infty\}$.  The group $\cox(\Delta,p)$ is the abstract group with
  generators $\lbrace r_i \colon i \in \Delta \rbrace$ indexed by
  $\Delta$, and with presentation
  \begin{equation*}
    \cox(\Delta,p) = \left< r_i \left|
    \begin{array}{cl}
      r_i^p = 1 & i \in \ver_\Delta\\
      r_ir_jr_i = r_jr_ir_j & \{i,j\} \in \edg_\Delta\\
      r_ir_j = r_jr_i & \{i,j\} \not \in
      \edg_\Delta
    \end{array}
    \right> \right.
  \end{equation*}
  When $p=2$ these are \emph{simply-laced Coxeter groups}.  When $p$
  is $\infty$, there are no vertex relations and the groups are
  \emph{simply-laced Artin groups}.  For general $p$ these are
  examples of the \emph{Shephard groups} recently introduced by
  Katherine Goldman \cite{Goldman25}.  We also define
  $\cox(\Delta,p,\infty)$ to be the group where the braid relations
  are omitted.  In this case $\cox(\Delta,2,\infty)$ is a
  \emph{right-angled Coxeter group} and $\cox(\Delta,\infty,\infty)$
  is a \emph{right-angled Artin group}.
\end{defn}

The reflection groups of Lorentzian lattices constructed from graphs
with $\lambda_2$ equal to $1$ or $2$ contain right-angled and
simply-laced Coxeter groups, respectively.

\begin{prop}[Coxeter Groups]
  Let $\Delta$ be a graph with $\lambda_2>0$ and let
  $(\cO_\Delta,\spanning_\Delta)$ be its Lorentzian representation.
  If $\lambda_2=1$ then $\cO = \Z$ and the reflection group
  $\Refl(\spanning_\Delta) \subset \Refl(\Z_\Delta)$ is the
  right-angled Coxeter group $\cox(\Delta,2,\infty)$.  If
  $\lambda_2=2$ then $\cO = \Z$ and the reflection group
  $\Refl(\spanning_\Delta) \subset \Refl(\Z_\Delta)$ is the
  simply-laced Coxeter group $\cox(\Delta,2)$.
\end{prop}

\begin{proof}
  When $\lambda_2=1$, the $\alpha_i \in \spanning_\Delta$ are short
  roots of $\Z_\Delta$ and the inner products between them show that
  there is a representation of the right-angled Coxeter group
  $\cox(\Delta, 2,\infty)$ acting on $\R_\Delta$.  When $\lambda_2=2$,
  the $\alpha_i \in \spanning_\Delta$ have norm $2$ and they are roots
  of $\Z_\Delta$.  The inner products between them show that there is
  a representation of the simply-laced Coxeter group $\cox(\Delta, 2)$
  acting on $\R_\Delta$.  Concretely, the Coxeter generator $r_i$ acts
  on $\Z_\Delta$ as the (order~$2$) reflection $R_{\alpha_i}$ in the
  root $\alpha_i$.  There is a corresponding action of the Coxeter
  group on the degenerate Lorentzian space $\R^\Delta$ and this is
  precisely the standard faithful Tits representation of this Coxeter
  group \cite[Section~4.1]{BjornerBrenti05}.  The induced action on
  the quotient space $\R_\Delta$ is also faithful by a result of Daan
  Krammer \cite[Proposition~6.1.3]{Krammer09}.
\end{proof}

Graph automorphisms are a second source of lattice automorphisms and
they are compatible with what we call the standard decomposition.

\begin{remark}[Standard Decomposition]\label{rem:decomposition}
  Let $\R_\Delta = \R_{\Delta, \lambda_2}$ be the Lorentzian space of
  a connected graph $\Delta$ with $\lambda_2>0$. The principal
  eigenvector of $A_\Delta$ spans a $1$-dimensional negative definite
  subspace in $\R^\Delta_{\lambda_2}$ (Definition~\ref{def:graphs})
  which projects to the $1$-dimensional negative definite subspace
  $\R_{\Delta}^-$ in the quotient space $\R_{\Delta}$.  Let
  $\R_{\Delta}^+$ be the $g$-dimensional positive-definite orthogonal
  complement of $\R_{\Delta}^-$.  We call $\R_{\Delta} = \R_{\Delta}^-
  \oplus \R_{\Delta}^+$ the \emph{standard decomposition} of
  $\R_{\Delta}$.
\end{remark}

\begin{defn}[Graph Automorphisms]\label{def:aut-graph}
  Let $\Delta$ be a graph with $r= \lambda_2>0$ and consider the
  $\cO$-module $\cO_\Delta$ in the space $\R_\Delta$ with standard
  decomposition $\R_\Delta = \R_{\Delta}^+ \oplus \R_{\Delta}^-$
  (Remark~\ref{rem:decomposition}).  The symmetry group $\aut(\Delta)$
  of the graph acts on the Hermitian $\cO$-module $\cO^{\Delta}_r$ by
  permuting the basis vectors $\basis_\Delta = \lbrace e_i \colon i
  \in \Delta \rbrace$.  These graph symmetries preserve the adjacency
  matrix $A_\Delta$ and the gram matrix $\gram(\basis_\Delta) =
  \gram(\spanning_\Delta) = rI-A_\Delta$.  As a consequence,
  $\aut(\Delta)$ acts on the quotient $\cO_\Delta$ permuting the
  spanning vectors $\spanning_\Delta$
  (Corollary~\ref{cor:symmetries}).  Since automorphisms of $\Delta$
  fix the principal eigenvector in $\R^\Delta$ and its image $v_- \in
  \R_\Delta$, the group $\aut(\Delta)$ stabilizes the line
  $\R_{\Delta}^- = \R v_-$, so it also acts on the positive definite
  orthogonal complement $\R_{\Delta}^+ \simeq \R^g$.  In other words,
  we have a representation
  \begin{equation*}
    \aut(\Delta) 
    \to \aut(\R_{\Delta}^+) 
    \subseteq \aut(\R_{\Delta}) 
    \simeq O(g, 1). 
  \end{equation*}
\end{defn}

Given $S \subseteq \Delta$, we define
\begin{equation}\label{eq:vector-sums}
e_S = \sum_{i \in S} e_i 
\text{\; and \;}
\alpha_S = \sum_{i \in S} \alpha_i.
\end{equation}
When $\Delta$ is $k$-regular, $\one = e_{\Delta}$ is a principal
eigenvector of $A_{\Delta}$ and $\ip{e_{\Delta}}{e_i} = r - k$ in
$\R^\Delta$ for each $i \in \Delta$, so $e_{\Delta}$ has norm $v (r -
k) < 0$ and its image $\alpha_{\Delta} \in \cO_{\Delta}$ spans the
line $\R_{\Delta}^-$. In some examples $\alpha_{\Delta}$ behaves like
a Weyl vector.

\begin{proposition}[Weyl vector-like]\label{prop:weyl-like}
  Let $\Delta$ be a $k$-regular graph.  The action of $\aut(\Delta)$
  on the Lorentzian space $\R_\Delta$ permutes $\spanning_\Delta =
  \{\alpha_i\}$ and fixes $\alpha_\Delta$. The corresponding fixed
  point $\proj\alpha_\Delta$ in hyperbolic space is equidistant from
  the hyperplanes $\alpha_i^{\bot}$.
\end{proposition}

\begin{proof}
  Since the action permutes the spanning set, it fixes their sum
  $\alpha_\Delta$, and the hyperplanes $\alpha_i^\bot$ are equidistant
  from $\proj \alpha_\Delta$ because $\ip{\alpha_\Delta}{\alpha_i} =
  r-k$ for all $i\in \Delta$.
\end{proof}

\begin{proposition}[Reflections and Automorphisms]\label{prop:cox-aut-delta}
  Let $\Delta$ be a graph with $\lambda_2=1$ or $2$ and let
  $(\Z_\Delta,\spanning_\Delta)$ be its Lorentzian representation.  If
  $\lambda_2=1$ then there is a representation of the semi-direct
  product
  \begin{equation*} 
    \cox(\Delta,2,\infty) \rtimes \aut(\Delta) \to \aut(\Z_\Delta) \subseteq
    \aut(\R_\Delta) \simeq O(g, 1).
  \end{equation*}
  Similarly, if $\lambda_2=2$ then there is a representation of the
  semi-direct product
  \begin{equation*} 
    \cox(\Delta,2) \rtimes \aut(\Delta) \to \aut(\Z_\Delta) \subseteq
    \aut(\R_\Delta) \simeq O(g, 1).
  \end{equation*}
  In particular, these groups also act on the hyperbolic space $\ball
  \R_{\Delta} \simeq \hyp^g_\R$.
\end{proposition}

\begin{proof}
  The Coxeter generator $r_i$ acts on $\Z_\Delta$ as the (order~$2$)
  reflection $R_{\alpha_i}$ in the root $\alpha_i$
  (Proposition~\ref{prop:cox-aut-delta}), and there is an action of
  $\aut(\Delta)$ on $\Z_\Delta$ (Definition~\ref{def:aut-graph}). To
  combine these, the action of the element $(r_i, \sigma)$ on $x \in
  \Z_\Delta$ is defined by the formula $(r_i, \sigma) x = R_
    {\alpha_i} \sigma(x)$.
\end{proof}

\section{Strongly regular graphs}\label{sec:srg}

Strongly regular graphs have uniform structural properties that are
highly compatible with the graph-to-lattice construction.  In this
section we review some general facts about strongly regular graphs.
See \cite{BrouwerVanMaldeghem22} for additional details.

\begin{defn}[Strongly regular graphs]\label{def:srgs}
  A graph $\Delta$ is called a \emph{strongly regular graph} with
  parameters $(v, k, \lambda, \mu)$ if it is a $v$-vertex $k$-regular
  graph of diameter~$2$, where every pair of distinct adjacent
  vertices in $\Delta$ have $\lambda$ common neighbors and every pair
  of distinct non-adjacent vertices have $\mu$ common neighbors.  To
  simplify the statements, we exclude the trivial cases, such as
  complete graphs (Remark~\ref{rem:indistinct}), that have fewer than
  $3$~eigenvalues.  In this more restrictive sense, a strongly regular
  graph can be defined as a $k$-regular graph with exactly
  $3$~eigenvalues.  The spectrum of a strongly regular graph $\Delta$
  with parameters $(v,k,\lambda,\mu)$ is usually written
  $\spec(\Delta) = \{k^1 r^f s^g\}$ with $k > r > 0 > s$.  When the
  multiplicities $f$ and $g$ are not equal, one can show that all of
  these values are integers.  The complement of a strongly regular
  graph is also strongly regular.
\end{defn}

\begin{remark}[Parameters and Graphs]\label{rem:parameters-graphs}
  As Brouwer and Van Maldeghem comment, ``strongly regular graphs live
  on the boundary between the crystalline world and the random world''
  \cite[1.1.18]{BrouwerVanMaldeghem22}.  Even once one has found a set
  of parameters satisfying all of the necessary identities
  (Remark~\ref{rem:identities}), there may be no graph, a unique graph
  or a quite large collection of examples with these parameters.  We
  are particularly interested in the parameter sets where there is a
  unique graph since their Lorentzian representations are often lattices 
  with large automorphism groups.
\end{remark}

The adjacency matrix of a strongly regular graphs generates a small
rank algebra.  Let $J$ be the $v \times v$ ``all-ones'' matrix and
$\one$ be the $v \times 1$ ``all-ones'' vector.
  
\begin{defn}[Bose-Mesner algebras]\label{def:bose-mesner}
  Let $\Delta$ be a connected strongly regular graph with parameters
  $(v,k,\lambda,\mu)$, adjacency matrix $A$, and spectrum $\spec(A) =
  (k^1, r^f, s^g)$ with $k > r > 0 > s$.  The \emph{Bose-Mesner
  algebra} of $\Delta$ is commutative associative algebra generated by
  the identity matrix $I$, the ``all ones matrix'' $J$ and its
  adjacency matrix $A=A_\Delta$.  It is a free $\Z$-module of rank~$3$
  with multiplication rules $J^2 = v J$, $A J = J A = k J$ and
  \begin{equation}
    A^2 - (\lambda-\mu) A + (\mu-k)I = (A - rI)(A-sI) = \mu J 
    \label{eq:adj-matrix-eq2}
  \end{equation}
\end{defn}

The parameters associated to a strongly regular graph are highly
constrained. For instance, they satisfy the basis identity
\begin{equation}\label{eq:basis-identity}
\mu ( v - k -1) = k(k - \lambda - 1).
  \end{equation}
  
\begin{remark}[Identities]\label{rem:identities}
  Let $\Delta$ be a connected non-trivial strongly regular graph with
  parameters $(v,k,\lambda,\mu)$, adjacency matrix $A = A_\Delta$, and
  spectrum $\spec(\Delta)=\{k^1\ r^f\ s^g\}$ where $k>r>0>s$ are the
  three distinct eigenvalues of $\Delta$ having multiplicities $1$,
  $f$, and $g$.  Here are some of the many identities that hold among these
  variables, which we record for later use.
  \begin{enumerate}
  \item[(a)] $r+s = \lambda-\mu$ and $rs = \mu-k$
  \item[(b)] $1+f+g = v$ and $k+rf+sg = 0$
  \item[(c)] $k^2 - (\lambda-\mu)k + (\mu-k)v = (k-r)(k-s) = v\mu$
  \item[(d)] 
  $rv-r+k = k(r+1)(k-r)/\mu$
  \end{enumerate}
  \noindent
  The proofs are elementary. It follows from (b) that 
  when $f \neq g$, the eigenvalues are in $\Z$. And note 
  that identity (d) is a consequence of (a) and (c).
\end{remark}

Strongly regular graphs whose eigenvalues are integers have been
coarsely classified.  The following theorem was originally proved by
Sims (unpublished) and then reproved by Neumaier \cite{Neumaier79}.

\begin{thm}[$s=-m$]\label{thm:classification}
  The strongly regular graphs $\Delta$ with smallest eigenvalue
  $s=-m$, $m\geq 2$ integral, are the following:
  \begin{itemize}
  \item[(a)] Complete multipartite graphs with classes of size $m$,
  \item[(b)] Latin square graphs $LS_m(n)$,
  \item[(c)] Steiner graphs $S_m(n)$,
  \item[(d)] Finitely many other graphs.
  \end{itemize}
\end{thm}

The special case with $s=-2$ was established earlier by Seidel
\cite{Seidel68}.

\begin{thm}[$s = -2$]\label{thm:s=-2}
  The only strongly regular graphs with smallest eigenvalue $s=-2$ are
  the cocktail party graphs (= complement of a set of disjoint edges),
  the triangular graphs $T(n+2) = S_2(n)$, the lattice graphs $L_2(n)
  = LS_2(n)$, the Shrikhande graph, the three Chang graphs and the
  graphs of Petersen, Clebsch, and Schl\"afli.
\end{thm}

The complement of a strongly regular graph with $s=-m$ is a strongly
regular graph with $r=m-1$.  Thus, when each of the graphs in
Theorem~\ref{thm:s=-2} is replaced with its complement, this becomes a
classification of strongly regular graphs with second eigenvalue
$r=1$.

\begin{cor}[$r=1$]\label{cor:r=1}
  The only connected strongly regular graphs with second eigenvalue
  $r=1$ are the triangle graph complements, the lattice graph
  complements, and the complements of the last seven specific graphs 
  listed in Theorem~\ref{thm:s=-2}.
\end{cor}

There are similar lists for strongly regular graphs with $s=-3$ and
their complements with $r=2$.  For definitions of the strongly regular
graphs listed in Theorem~\ref{thm:classification},
Theorem~\ref{thm:s=-2}, and Corollary~\ref{cor:r=1}, see
\cite{BrouwerVanMaldeghem22}.  Here we are primarily interested in the
structure of the classification.

\begin{remark}[Structure of the Classification]\label{rem:structure}
  According to Theorem~\ref{thm:classification}, the strongly regular
  graphs with $s=-m$ (or their complements with $r=m-1$) belong to a
  small number of infinite families with extra structure, or they are
  among a small number of exceptions.  The exceptions can be split
  into two types.  There are a finite number of graphs which have the
  parameters of Steiner graphs or Latin square graphs without
  satisfying the precise definition; these are called
  \emph{pseudo-Steiner graphs} and \emph{pseudo-Latin square graphs}.
  And then there are a finite number of graphs with \emph{exceptional
  parameter sets}.  For $s=-2$ (and $r=1$) there $3$ exceptional
  parameter sets.  For $s=-3$ (and $r=2$) there are $64$ exceptional
  parameter sets.  The example sections include lattices built from an
  infinite family (Section~\ref{sec:kneser}), and a lattice built from
  a strongly regular graph with an exceptional parameter set
  (Section~\ref{sec:mclaughlin}).
\end{remark}

\section{Lattices from strongly regular graphs}\label{sec:srg-lattices}

This section records general facts about the Lorentzian lattices
constructed from strongly regular graphs with integer eigenvalues.
Note that we assume $r \in \Z$ (and $\cO=\Z$) for simplicity of
exposition, and it is common for strongly regular graphs to have
integer eigenvalues.  The basic observation is that for interesting
graphs, the corresponding Lorentzian lattices $ \cO_\Delta$ are also
interesting, particularly when $\Delta$ is highly symmetric and the
multiplicity of $\lambda_2$ is large. For more on this perspective see
the introduction of \cite{NeumaierSeidel83}.

\begin{notation}[Strongly Regular]\label{not:strongly-regular}
  Throughout this section, $\Delta$ is a non-trivial connected
  strongly regular graph with parameters $(v,k,\lambda,\mu)$ and
  spectrum $\{k^1 \; r^f \; s^g \}$ with $k>r>0>s$ and $r,s \in \Z$.
  The form on the singular space $\R^\Delta_r \simeq \R^{g,1,f}$ has
  $3$ eigenvalues, $(r-s) > 0 > (r-k)$ of multiplicity $g$, $f$, and
  $1$, respectively, and the Lorentzian quotient $\R_\Delta =
  \R_{\Delta,r} \simeq \R^{g,1}$ contains the Lorentzian $\Z$-lattice
  $\Z_\Delta$ with its $v$-element spanning set $\spanning_\Delta =
  \{\alpha_i\}$.  Let $\alpha_\Delta= \sum_{i\in \Delta} \alpha_i \in
  \Z_\Delta$ be the image of the all $1$'s vector $\one =
  (1,1,\cdots,1)^\tr \in \Z^\Delta$ in the quotient.  Since $\Delta$
  is $k$-regular, $\one$ is an $(r-k)$-eigenvector in $\Z^\Delta$.
  The $\R$-span and the orthogonal complement of $\alpha_\Delta$
  determine the standard decomposition
  (Remark~\ref{rem:decomposition}).
\end{notation}

Some interesting and highly symmetric lattices arise from strongly
regular graphs where the multiplicity of its second eigenvalue is
large.

\begin{remark}[Large Radical]\label{rem:large-radical}
  When the multiplicity $f$ is large, the dimension of the radical of
  $\Z^\Delta$ is large, and the size of the spanning set
  $\spanning_\Delta$ is much larger than the dimension of the
  Lorentzian space $\R_\Delta$ containing it. For example, the
  strongly regular McLaughlin graph $\Delta$ (studied in
  Section~\ref{sec:mclaughlin}) has $275$ vertices and spectrum
  $\{112^1 \; 2^{252} \; (-28)^{22}\}$.  The $275$-dimensional space
  $\R^\Delta_r$ has a $252$-dimensional radical and the quotient lattice
  $\Z_\Delta \subset \R_\Delta \simeq \R^{22,1}$ has a $275$ element
  spanning set.  The symmetries of this lattice include a semi-direct
  product of an infinite $275$-generator simply-laced Coxeter group
  and the $1\ 796\ 256\ 000$ element automorphism group
  $\aut(\Delta)$, which contains the sporadic McLaughlin finite simple
  group as an index~$2$ subgroup
  (Proposition~\ref{prop:cox-aut-delta}).
\end{remark}

For lattices constructed from strongly regular graphs, there is a
uniform way to construct elements in the radical using special
elements.

\begin{definition}[Special Elements]\label{def:r-special}
  For each vertex $i \in\Delta$ we define a non-negative column vector
  $\tau_{i} \in \Z^\Delta$ which is the $i$-th column of the matrix
  $A-sI$.  
  If we let $N(i)$ denote the set of neighbors of vertex $i \in \Delta$ 
  then
  \[
  \tau_{i} = e_{N(i)} - s e_i
  \]
  where $e_{N(i)}$ is $\sum_{i\in N(i)} e_i$ (Equation~\ref{eq:vector-sums}).
  Note that the linear combination is non-negative because $-s>0$.
\end{definition}

\begin{lemma}[Radical Elements]\label{lem:elts-in-rad}
  When $\Delta$ is a strongly regular graph with parameters
  $(v,k,\lambda,\mu)$, the basis elements $e_i$ and the special
  elements $\tau_j$ and have constant inner product
  $\ip{e_i}{\tau_{j}} = - \mu$ for all $i,j \in \Delta$. Thus $\tau_i
  - \tau_j \in \rad(\Z^\Delta_r)$ for all $i, j \in \Delta$.  In fact,
  the vectors $\tau_i - \tau_j$ span the radical of $\R^\Delta_r$.
\end{lemma}

\begin{proof}
  The $e_i$ are column vectors of the identity matrix $I$ and the
  $\tau_j$ are column vectors of $A-sI$.  Using
  Equation~\ref{eq:adj-matrix-eq2} in Definition~\ref{def:bose-mesner}
  we compute
  \[ 
  \llbracket \ip{e_i}{\tau_j} \rrbracket = I^\tr (rI-A)(A-sI) = -\mu J
  \]
  which shows that $\ip{e_i}{\tau_j}=-\mu$ for all $i,j\in \Delta$.
  One has $\tau_i - \tau_j \in \rad(\Z^\Delta_r)$ because this
  difference is orthogonal to all of the basis vectors.  
  
  To see that the vectors $\tau_i - \tau_j$ span the radical, write $Z
  = \R\lbrace \tau_i - \tau_j \colon i, j \in \Delta \rbrace$.  Note
  that $ \R^{\Delta} \simeq \R^v = \R \mathbf{1} \bot \R^f \bot \R^g$
  where
  \[
  \R \mathbf{1} = \op{ker}(A - k I), \;
  \R^f = \op{ker}( A - r I), \; \R^g = \op{ker}(A - s I)
  \]
  and where $\bot$ means orthogonal with respect to the standard dot
  product on $\R^v$.  Note that $(A - sI)$ kills $\R^g$ and acts on
  $\R \mathbf{1}$ and $\R^f$ as nonzero scalars, so
  \[
  \R \lbrace \tau_i \colon i \in \Delta \rbrace
  =\op{im}(A - s I) = \R \mathbf{1} \bot \R^f.
  \]
  Since $Z$ together with a single $\tau_i$ spans this space
  $\op{im}(A - s I)$, we have $\op{dim}(Z) \geq f$.  But each $\tau_i
  - \tau_j$ has coordinate sum zero, so is orthogonal to $\mathbf{1}$
  with respect to the standard dot product on $\R^v$, so $Z \subseteq
  \R^f$. It follows that $ Z = \R^f = \op{rad}(\R^\Delta_r)$.
\end{proof}

\begin{defn}[Common Image $\tau$]\label{def:tau}
  Since the special elements $\tau_i$ differ by elements of the
  radical, they have a common image $\tau$ in the quotient lattice
  $\Z_\Delta$. By Lemma~\ref{lem:elts-in-rad} the inner product
  $\ip{\alpha_i}{\tau} = - \mu$ for all $i \in \Delta$ and, if we
  write $\alpha_{N(i)}$ for the image of $e_{N(i)}$, we have $v$
  different expressions of $\tau$: for any vertex $i\in \Delta$, we
  have $\tau = \alpha_{N(i)} - s \alpha_i$.  Adding all $v$ of these
  expressions together yields
  \begin{equation}\label{eq:tau-alpha}
      v \tau = (k - s) \alpha_\Delta
  \end{equation} 
  which shows that $\tau$ and $\alpha_\Delta$ are rationally related
  with $\proj \tau = \proj \alpha_\Delta$.  By
  Proposition~\ref{prop:weyl-like} $\tau$ has Weyl vector-like
  properties.  The inner product of both sides with $-\alpha_i$
  produces $v\mu = (k - s) ( k - r)$, identity (c) in
  Remark~\ref{rem:identities}, which is equivalent to the basic
  identity \eqref{eq:basis-identity} among the four parameters of a
  strongly regular graph.  Equation \eqref{eq:tau-alpha} implies that
  $\Z \lbrace \tau, \alpha_{\Delta} \rbrace$ is a $\Z$-module of rank
  $1$. To find a generator, we write $\tfrac{k-s}{v} = \tfrac{p}{q}$
  with $\gcd(p,q)=1$, and choose $a, b \in \Z$ with $ap+bq=1$. Then
  $\Z \lbrace \tau, \alpha_{\Delta} \rbrace = \Z\tau_\Delta$, where
  \[
  \tau_\Delta = a\tau + b\alpha_\Delta = \tfrac{1}{p} \tau = \tfrac{1}{q} \alpha_\Delta.
  \]
\end{defn}

One can project the $\alpha_i$ onto $(\alpha_\Delta)^\bot$ to get the
Euclidean representation of $\Delta$.

\begin{remark}[Euclidean Representation]\label{rem:euclid-rep}
  For a strongly regular graph $\Delta$, the $v$ vectors $(v\alpha_ i
  - \alpha_\Delta)$ span the $g$-dimensional positive definite
  subspace $(\alpha_\Delta)^\bot$, and the unit vectors $\xi_i$ in
  these directions satisfy
  \begin{equation*}
    \ip{\xi_i}{\xi_j} = \begin{cases}
      1 & \text{if \;} i = j,\\
      s/k & \text{if \;} \{ i, j \} \in \edg_\Delta, \\
      -(s + 1)/(v - k - 1) & \text{otherwise}.
    \end{cases}
  \end{equation*}
  This is the \emph{Euclidean representation} of $\Delta$. See
  \cite{Bondarenko14, BPR17, BrouwerVanMaldeghem22}.
\end{remark}

Alternatively, one can project $\alpha_\Delta$ onto the hyperplanes
$\alpha_i^\bot$.  This results in $v$ points in hyperbolic space where
the pairwise distances also only depend on $\Delta$.  To our
knowledge, this variation remains completely unexplored in the
literature.

\begin{defn}[Points on Hyperplanes]\label{def:hyperplane-pts}
  Let $q_i$ be the linear projection of $r \alpha_\Delta$ onto
  $\alpha_i^\bot$ and let $p_i = \proj q_i$.  Concretely, for each
  $i\in \Delta$,
  \[
  q_i =
  \ip{\alpha_i}{\alpha_i} \alpha_\Delta - \ip{\alpha_\Delta}{\alpha_i}
  \alpha_i = r \alpha_\Delta + (k-r)\alpha_i.
  \]
  The vector $q_i$ has negative norm, so the point $p_i \in \ball
  \R_\Delta \simeq \hyp^g_\R$.  Moreover, these points have pairwise
  distances that only depends on whether $\{i,j\}$ is an edge of
  $\Delta$ or not.  Note that $p_i$ is the projection of $\proj
  \alpha_{\Delta}$ on the totally geodesic hyperplane
  $\alpha_i^{\bot}$ in $\ball \R_{\Delta}$.  We write
  \begin{equation*}
    d_e = d(p_i,p_j) \text{\; if \;} \{ i, j \} \in
    \edg_\Delta \text{\; and \;} d_{n e} = d(p_i,p_j) \text{\; if \;}
    \{ i, j \} \notin \edg_\Delta.
  \end{equation*}
\end{defn}

The next lemma computes the inner products in terms of the basic
parameters.

\begin{lemma}\label{lem:de-dne}
  The vectors $q_i \in \alpha_i^\bot$ have inner products
  \begin{equation*}
    \ip{q_i}{q_j} = \begin{cases} 
    -r(k-r)(rv+k-r) & \text{if \;} i=j,\\ 
    -r(k-r)(rv+k-r + (k-r)(r+1)/r ) & \text{if \;} \{i,j\} \in \edg_\Delta, \\ 
    -r(k-r)(rv+k-r + (k-r)) & \text{if \;} \{ i, j \} \notin \edg_\Delta,    
    \end{cases}
  \end{equation*}
  all of which are negative.  As a consequence, one has
  \begin{equation*}
    \op{cosh}( d_e) = 1 + \tfrac{\mu}{k r} \text{\; and \;} \op{cosh}(
    d_{n e}) = 1 + \tfrac{\mu}{k ( r + 1)}.
  \end{equation*}
\end{lemma}

\begin{proof}
  The values of $\ip{q_i}{q_j}$ are simple computations, and they are
  strictly negative since all three factors (after the initial minus
  sign) are positive.  Note that because the vectors $q_i$ all have
  the same length, the hyperbolic distance between 
  $p_i$ and $ p_j$ (Definition~\ref{def:hyperbolic}) can be
  computed as
  \begin{equation*}
    \cosh d(p_i,p_j) = \tfrac{\ip{q_i}{q_j}}{\ip{q_i}{q_i}}.
  \end{equation*}
  Using the computed values of $\ip{q_i}{q_j}$, one has
  \begin{equation*}
    \cosh d_e = 1 + \tfrac{ (k - r)(r + 1)/r}{ r v + k - r} 
    \text{\; and \;}
    \cosh d_{n e} = 1 + \tfrac{ (k - r)}{ r v + k - r} .
  \end{equation*}
  The lemma follows after some simplification using
  identity (d) of Remark~\ref{rem:identities}.
\end{proof}

This hyperbolic representation of $\Delta$ can be used to prove
structural results about strongly regular graphs, similar to earlier
uses the Euclidean representation.

\begin{prop}[Edges]\label{prop:edge}
  Let $\Delta$ be a strongly regular graph, let $\{i,j\}$ be an edge
  of $\Delta$, and let $\Delta_{i j} = \Delta \setminus ( N(i) \cup
  N(j))$ be the induced subgraph on vertices not connected to $i$ or
  $j$.  If $\Delta_{i j}$ contains an edge, then $r \geq 1$.  If
  $\Delta_{i j}$ contains a vertex of valency at least $2$, then $r >
  1$.
\end{prop}

\begin{proof}
  Let $q_i \in \alpha_i^\bot$ be the vectors of
  Definition~\ref{def:hyperplane-pts}.  The vector $q_i + q_j$
  represents the midpoint $m_{ij} = \proj(q_i+q_j)$ of the hyperbolic
  geodesic between $p_i = \proj q_i$ and $p_j = \proj q_j$.  Let $1
  \in \Delta_{ij}$ be a vertex not connected to $i$ or $j$ and let
  $d_0 = d( m_{i j},p_1)$.  Since $p_1$ is equidistant from $p_i$ and
  $p_j$, the triangle with vertices $p_1$, $p_i$, $p_j$ is isosceles
  and the triangle with vertices $p_1$, $p_i$, $m_{ij}$ is a right
  triangle.  By the hyperbolic Pythagorean theorem
  \begin{equation*}
    \cosh d_{ne} = \cosh d (p_i, p_1) 
    = \cosh d(p_i, m_{i j}) \cosh d( p_1, m_{i j} ) 
    = \cosh (d_e/2) \cosh d_0.
  \end{equation*}
  If $\{ 1, 2 \}$ is an edge in $\Delta_{i j}$, then $d(p_t, m_{i j})
  \leq d_0$ for $t = 1, 2$, and so $d(p_1, p_2) \leq 2 d_0$ by the
  triangle inequality .  It follows that
  \begin{equation*}
    \cosh (d_e/2) = \cosh (d(p_1, p_2)/2) \leq \cosh d_0 = \cosh
    d_{ne}/ \cosh(d_e/2).
  \end{equation*}
  By Lemma~\ref{lem:de-dne} and the double-angle formula, this is
  equivalent to
  \begin{equation*}
    \tfrac{1}{2} \left( 2 + \tfrac{\mu}{k r} \right) = \tfrac{1}{2}
    (\cosh(d_e) + 1) = \cosh^2 \tfrac{d_e}{2} \leq \cosh d_{ne} = 1 +
    \tfrac{\mu}{k (r+ 1)}
  \end{equation*}
  which simplifies to $\mu/(2 k r) \leq \mu/(k (r + 1))$ or $(r + 1)
  \leq 2 r$ and so $r \geq 1$.  This proves the first statement.  When
  equality holds and $r = 1$, we have $d(p_1, p_2)/2 = d_0 =
  d(p_1,m_{ij}) = d(m_{ij},p_2)$, which forces $m_{i j}$ to be the
  midpoint of the hyperbolic geodesic between $p_1$ and $p_2$.
  Suppose $1, 2, 3$ are distinct vertices of $\Delta_{ij}$ and $\{ 1,
  2 \}$ and $\{ 2, 3 \}$ are edges in $\Delta_{ij}$.  We already know
  that $r \geq 1$.  If equality holds and $r = 1$, then $m_{i j}$ is
  the midpoint of the geodesic between $p_1$ and $p_2$ and the
  midpoint of the geodesic between $p_2$ and $p_3$. But this forces
  $p_1 = p_3$, contradicting the fact that $d(p_1, p_3)$ is $d_{n e}$
  or $d_e$ and both are positive.
\end{proof}

Proposition~\ref{prop:edge} can also be proved using standard
techniques such as eigenvalue interlacing. The point here is to show
that the hyperbolic geometry of the configuration can be useful to
establish results about graphs.  For more graph theory applications 
of Lorentzian representations, see \cite{BaMc-graphs}.  Geometrically, 
the graph $\Delta_{ij}$ is the subgraph determined by the elements of the
spanning set orthogonal to both $\alpha_i$ and $\alpha_j$, so we call
it the \emph{orthogonal complement} of the edge $\{i,j\}$.

\begin{remark}
  The Kneser graphs $\kneser(n,2)$, whose lattices are studied in the
  next section, have $ r = 1$ and they are consistent with
  Proposition~\ref{prop:edge}.  The orthogonal complement of any edge
  in $\kneser(n,2)$ has exactly four vertices and two disjoint edges.
  In particular, no vertex has valency more than $1$.
\end{remark}

\section{Example: Kneser graphs}\label{sec:kneser}

This section constructs lattices from triangle graph
complements. These strongly regular graphs are special types of
Kneser graphs.

\begin{figure}
  \includegraphics{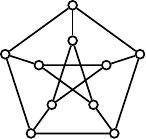}
  \caption{The Petersen graph $\petersen$.}
  \label{fig:petersen}
\end{figure}

\begin{defn}[Kneser Graphs]\label{def:kneser-graphs}
  The \emph{Kneser graph} $\kneser(n,m)$ is the graph $\Delta$ whose
  vertices are the $m$-element subsets of 
  $\lbrace 1, 2, \dotsb, n \rbrace$
  with an edge if and only if the subsets are disjoint.
  When $m=2$, the
  Kneser graph $\Delta_n = \kneser(n,2)$ is a strongly regular graph
  with parameters
  \[
  (v, k, \lambda, \mu) = \left(\tbinom{n}{2}, \tbinom{n-2}{2},
  \tbinom{n-4}{2}, \tbinom{n-3}{2} \right)
  \]
  and for $n\geq 5$ it is a connected graph with spectrum
  $\spec(\Delta_n)= \{k^1 \; r^f \; s^g\}$ where $r = 1$ with
  multiplicity $f=\binom{n}{2}-n$ and $s=3-n$ with multiplicity
  $g=n-1$.  Permutations of $\lbrace 1, \dotsb, n \rbrace$ induce
  automorphisms of the graph $\Delta_n$.  In fact, $\aut(\Delta_n)$ is
  the symmetric group.  The graph $\kneser(n,2)$ is the complement of
  the triangular graph $T(n)$ whose vertices are $2$-element subsets
  of $\{1,\ldots,n\}$ and an edge if and only if the subsets overlap.
  The smallest connected graph in the $\kneser(n,2)$ family is the
  \emph{Petersen graph} $\petersen = \kneser(5, 2)$. See
  Figure~\ref{fig:petersen}.
\end{defn}

\begin{defn}[Kneser Lattices]\label{def:kneser-lattices}
  For $n\geq 5$, let $L_n = \Z_{\Delta_n}$ be the Lorentzian
  $\Z$-lattice of signature $(n-1,1)$ constructed from the graph
  $\Delta_n = \kneser(n,2)$ and let $\spanning_n =
  \spanning_{\Delta_n} = \{ \alpha_{ij} \colon 1\leq i<j \leq n \}$ be
  its $\binom{n}{2}$-element spanning set
  (Proposition~\ref{prop:lor-spaces}).  Note that $\cO = \Z$ since
  $r=\lambda_2=1$ is an integer.  The matrix $\gram(\spanning_n) = I -
  A$ where $A$ is the adjacency matrix of $\Delta_n$.  Concretely,
  \[
  \ip{\alpha_{ij}}{\alpha_{k \ell}} =
  \begin{cases}
  \ \ 1 & \textrm{if $\size{\{i,j\} \cap \{k,\ell\}} = 2$} \\
  \ \ 0 & \textrm{if $\size{\{i,j\} \cap \{k,\ell\}} = 1$}\\
  -1 & \textrm{if $\size{\{i,j\} \cap \{k,\ell\}} = 0$}
  \end{cases}
  \]
  Since the elements $\alpha_{ij} \in \spanning_n$ are short roots of
  norm $1$, reflections in these roots generate a right-angled Coxeter
  group $\cox(\Delta_n,2,\infty) \subseteq \Refl_*(L_n)$
  (Definition~\ref{def:refl}) and the semi-direct product
  $\cox(\Delta_n,2,\infty) \rtimes S_n \subseteq \aut(L_n)$
  (Proposition~\ref{prop:cox-aut-delta}).  It is unclear whether these
  are equalities or proper inclusions.
\end{defn}

\begin{remark}[Nested Lattices]\label{rem:nested}
  Since all of these graphs have the same second eigenvalue
  ($r=\lambda_2=1$), by Theorem~\ref{thm:latt-functor} the
  inclusions
  $\Delta_5 \into \Delta_6 \into \cdots$ lead to 
  inclusions
  $\spanning_5 \into \spanning_6 \into \cdots$ and $L_5 \into L_6
  \into \cdots$ compatible with $\aut(\Delta_n)$-actions on
  $\spanning_n$ and on $L_n$.
\end{remark}

The lattices $L_n$ exist with convenient explicit coordinates inside
$\I_{n,1}$.

\begin{defn}[Odd Coordinates]\label{def:odd-coord}
  Recall that $I_{n,1}$ denotes the unique odd integral self-dual
  $\Z$-lattice (Definition~\ref{def:important-lattices}).  We write
  vectors in $I_{n,1}$ with respect to an ordered orthogonal basis
  $\basis_{n,1} = \{e_0,e_1, \ldots, e_n\}$ where the norm of $e_0$ is
  $-1$, and the norm of $e_i$ is $1$ for all $i>0$.  For $x = \sum_0^n
  x_i e_i \in \I_{n,1}$, the coordinates $x_i \in \Z$ are recorded as
  the column vector $x = (x_0\mid x_1\ x_2\ \dotsb\ x_n)^\tr$ whose
  norm is $-x_0^2 + x_1^2 + \dotsb + x_n^2$. By listing the negative
  norm basis vector $e_0$ first, $I_{n-1, 1}$ naturally sits inside
  $I_{n, 1}$ as the vectors whose last coordinate is $0$.
\end{defn}

\begin{defn}[Sublattices]\label{def:odd-sublattices}
  Let $\spanning'_n$ be the $\binom{n}{2}$ vectors $\alpha'_{ij} = e_0
  + e_i + e_j$ with $i<j$ and $i,j \in \lbrace 1, 2, \dotsb, n
  \rbrace$ inside $\I_{n,1}$ in the ordered basis of
  Definition~\ref{def:odd-coord}, and let $L'_n$ be the $\Z$-span of
  $\spanning'_n$.  For example $\alpha'_{23} \in \spanning'_5$ is the
  vector $(1 \mid 0\ 1\ 1\ 0\ 0)^\tr$ in $\I_{5,1}$.  The vectors
  $\alpha'_{ij}$ have norm $1$ and the same inner products as the
  $\alpha_{ij} \in \spanning_n$ in $L_n$
  (Definition~\ref{def:kneser-lattices}).  In particular,
  $\gram(\spanning_n) = \gram(\spanning'_n)$ and the bijective map
  $\spanning_n \to \spanning'_n$ sending $\alpha_{ij} \to
  \alpha'_{ij}$ extends to a $\Z$-lattice isometry $L_n \to L'_n
  \subset \I_{n,1}$ (Proposition~\ref{prop:gram-spanning}).  We
  identify the vectors $\alpha_{ij}$ with the concrete vectors
  $\alpha'_{ij} = e_0 + e_i + e_j$ and the lattice $L_n$ with this
  explicit sublattice $L'_n \subset \I_{n,1}$ removing the primes.
\end{defn}

The lattice $L_n$ is the orthogonal complement of a specific vector in
$\I_{n,1}$.

\begin{proposition}[Complement]\label{prop:kneser-complement}
  For all $n \geq 5$, the Kneser lattice $L_n = \gamma^\bot \cap
  \I_{n,1}$, where $\gamma$ is the vector $\gamma = (2 \mid 1^n)^\tr =
  2 e_0 + e_1 + \cdots + e_n \in \I_{n,1}$. As a consequence, the
  discriminant of $L_n$ is $n-4$.
\end{proposition}

\begin{proof}
Since $\ip{\gamma}{\alpha_{i j} } = 0$ for all $i, j$, one has
$L_n \subseteq \gamma^{\bot} \cap I_{n, 1}$.
On the other hand, let $x$ be a vector in $\gamma^\bot \cap \I_{n,1}$
  with coordinates $x = (x_0 \mid x_1\ x_2\ \cdots\ x_n)^\tr$ and let
  $x' = x - \sum_{i=2}^n x_i \alpha_{1i}$, a vector of the form $x' =
  (a \mid b\ 0\ \cdots\ 0)^\tr$ with $a,b \in \Z$.  Since
  $\ip{x}{\gamma}=0$, and $\ip{\alpha_{ij}}{\gamma}=0$ for all $i,j$,
  we have $\ip{x'}{\gamma}=0$ which implies $b=2a$.  But the vector
  \[
  x' = (a \mid 2a\ 0\ \cdots\ 0)^\tr = a\cdot (1 \mid
  2\ 0\ \cdots\ 0)^\tr = a\cdot (\alpha_{12} +
  \alpha_{13}-\alpha_{23})
  \]
  is in $L_n$, the $\Z$-span of $\spanning_n$, so $x = x' +
  \sum_{i=2}^n x_i \alpha_{1i}$ is also in $L_n$.  The final assertion
  follows from the standard fact that a primitive sublattice $M =
  \gamma \Z$ in a self-dual lattice $I_{n , 1}$ has the same
  discriminant as its orthogonal complement $M^\bot$
  \cite[Proposition~1.2]{Ebeling13}.
\end{proof}

The next two propositions show that the point $\proj \tau$ is a
Weyl-like vector for these lattices: it is the unique point fixed by
the automorphism group, and the mirrors $\{ \alpha_{ij}^{\bot} \colon
\alpha_{ij} \in \spanning_n\}$ are precisely the short root mirrors
closest to $\proj \tau$.

\begin{proposition}[Unique Fixed Point]\label{prop:kneser-unique-fixed}
  In the Lorentzian representation of the Kneser graph $\Delta_n =
  \kneser(n,2)$ with $n\geq 5$, the action of the symmetric group
  $\aut(\Delta_n)$ on the hyperbolic space $\hyp_\R^{n-1}$ fixes the
  point $\proj \tau$ and only this point.
\end{proposition}
\begin{proof}
  Note that $\aut(\Delta_n)$ acts on $L_n$ by permuting the last $n$ coordinates, so 
  \[
  (n \mid 2^n) =
  \tfrac{2}{n-1} \alpha_\Delta = \tfrac{2}{n-3} \tau
  \]
  in $L_n$ is fixed by $\aut(\Delta_n)$.  If $x = (a \mid
  b_1\ b_2\ \cdots\ b_n)\in \R L_n$ is a vector of negative norm with
  $\proj x$ fixed by $\op{Aut}(\Delta_n)$, then $a \neq 0$, and this
  implies that $x$ must have the form $x = (a \mid b^n)$.  By
  Proposition~\ref{prop:kneser-complement}, $\ip{x}{\gamma}= -2a+nb=0$
  and $a = \tfrac{n}{2}b$.  Thus $x = (a \mid b^n) = \tfrac{b}{2}(n
  \mid 2^n)$ and $\proj x = \proj \tau$.
\end{proof}

\begin{proposition}[Closest Mirrors]\label{l-first-shell-around-tau}
  In the Lorentzian representation of the Kneser graph $\Delta_n =
  \kneser(n,2)$ with $n\geq 5$, the fixed point $\proj \tau \in \ball
  L_n$ is not contained in any short root mirror, and the short
  root mirrors closest to $\proj \tau$ are precisely the
  $\binom{n}{2}$ mirrors $\alpha_{ij}^{\bot}$ with $\alpha_{ij}$ in
  the spanning set.
\end{proposition}

\begin{proof}
  Let $\tau_0 = (n \mid 2^n) = \tfrac{2}{n-3}\tau$. Note that $\tau_0
  = 2 \gamma + (n-4) e_0$ where $\gamma = (2 \mid 1^n)$.  Since
  $\ip{x}{\gamma}=0$ for any $x \in L_n$
  (Proposition~\ref{prop:kneser-complement}),
  \[
  \ip{x}{\tau_0} = (n-4)\ip{x}{e_0} = -(n-4)x_0.
  \]
  Take a short root $x$ with $d( \proj \tau, x^{\bot}) \leq d( \proj
  \tau, \alpha_{i j}^{\bot})$.  Then by using the distance formula in
  Remark~\ref{rem:dist-hyperplane},
  \[
  (n - 4) \abs{x_0}
  = \abs{\ip{\tau_0}{ x}} \leq \abs{\ip{\tau_0}{\alpha_{i j}}} = (n - 4)
  \]
  and $x_0 = 0, \pm 1$.  If $x_0 = 0$, then $\sum_{i = 1}^n x_i = 0$
  (since $\ip{x}{\gamma}=0$).  Since $x_i$'s are not all zero, at
  least two of them are nonzero integers making the norm of $x$ at
  least $2$, so $x$ is not a short root.  This shows that $\proj \tau$
  is not on any short root mirrors. Finally, suppose that $x$ is a
  short root with $\size{x_0}=1$.  Since $x$ and $-x$ determine the
  same mirror, we may assume that $x_0=1$ without loss of generality.
  To be a short root, the sum $\sum_{i = 1}^n x_i^2 = 2$, which
  implies that exactly two of the coordinates $x_i$ and $x_j$ with
  $0<i<j$ are equal to $\pm 1$ and the rest of coordinates are $0$.
  But since $\ip{x}{\gamma} = 0$, we have $x_i + x_j = 2x_0 = 2$, so
  $x_i = x_j =1$ and $x = \alpha_{ij}$.
\end{proof}

The smallest Lorentzian Kneser lattice is of special interest.

\begin{rem}[Petersen Lattice]\label{rem:petersen-lattice}
  When $n=5$, the Kneser lattice $L_5 =\Z_\petersen$ where $\petersen
  = \kneser(5,2)$ is the Petersen graph.  By
  Proposition~\ref{prop:kneser-complement} $\Z_\petersen$ is
  unimodular and thus $\Z_\petersen \simeq \I_{4,1}$
  (Definition~\ref{def:important-lattices}).  One can also construct a
  (complex) Lorentzian $\eisen$-lattice $\eisen_\petersen = \eisen
  \otimes_{\Z} \Z_\petersen$ by tensoring $\Z_\petersen$ with the
  Eisenstein integers and promoting the $\Z$-valued bilinear form to a
  $\eisen$-valued Hermitian form.  This complex lattice and its
  complex reflection group $\Refl(\eisen_\petersen)$ have a close
  connection to the moduli space of cubic surfaces in $\proj^3_\C$.
  See Example~\ref{ex:eisen-petersen} in
  Section~\ref{sec:ball-quoients}.
\end{rem}

\section{Example: McLaughlin graph}\label{sec:mclaughlin}

Many of the strongly regular graphs with $r=2$ and an exceptional
parameter set (Remark~\ref{rem:structure}) have Euclidean
representations inside the Leech lattice $\Lambda$. Standard
references include \cite[Section~5.5]{Wilson09},
\cite[Chapter~10]{BrouwerVanMaldeghem22}, and the original article
\cite{Wilson89}.  In this section we highlight that the Lorentzian
lattices constructed from these graphs can be found inside the
remarkable lattice $\II_{25,1}$.  As an illustration we provide full
details for the Lorentzian lattice of the McLaughlin graph and prove
the analogs of Propositions~\ref{prop:kneser-unique-fixed}
and~\ref{l-first-shell-around-tau} for this lattice.

\begin{remark}[History]\label{rem:history}
  The Lorentzian lattice $\II_{25,1}$ has been studied extensively by
  Conway and Borcherds, and its reflection group $\Refl(\II_{25, 1})$
  is, arguably, the most beautiful example of a Lorentzian reflection
  group.  See \cite{ConwaySloane99} Chapters 24, 26, 27 and
  \cite{Borcherds85-leech} for details.  The lattice $\II_{25, 1}$ has
  been used to study lattices and hyperbolic reflection groups in
  dimensions around $24$ to great success by Borcherds who
  demonstrated that many Lorentzian lattices with interesting
  reflection groups occur inside $\II_{25,1}$
  \cite{Borcherds87-aut-LL, Borcherds00-refl-LL}.  In fact, the way we
  realize McLaughlin lattice inside $\II_{25, 1}$ goes back to
  \cite{Borcherds87-aut-LL}.  For additional details about the Leech
  lattice and the reflection group of the lattice $\II_{25,1}$, see
  \cite{ConwaySloane99, Borcherds85-leech, Borcherds87-aut-LL}.
\end{remark}

\begin{table}
    \centering
    \begin{tabular}{|c|c|l|}
    \hline
    $(v,k,\lambda,\mu)$ & Spectrum & Notes\\
    \hline
    \hline
      $(50,7,0,1)$      & $7^1 2^{28} (-3)^{21}$ & $\Delta_{50}=$ Hoffman-Singleton graph \\ \hline
      $(56,10,0,2)$     & $10^1 2^{35} (-4)^{20}$ & $\Delta_{56}=$ Gewirtz graph \\ \hline
      $(77,16,0,4)$     &  $16^1 2^{55} (-6)^{21}$ & $\Delta_{77}=$ $M_{22}$-graph \\ \hline
      $(81,20,1,6)$     & $20^1 2^{60} (-7)^{20}$  & $\Delta_{81}=$ Brouwer-Haemers graph \\ \hline
      $(100,22,0,6)$    &  $22^1 2^{77} (-8)^{22}$ & $\Delta_{100}=$ Higman-Sims graph \\ \hline
      $(105,32,4,12)$   & $32^1 2^{84} (-10)^{20}$  & $\Delta_{105}=$ inside $\Delta_{275}$ \\ \hline
      $(112,30,2,10)$   &  $30^1 2^{90} (-10)^{21}$ & $\Delta_{112}=$ inside $\Delta_{275}$ \\ \hline
      $(120,42,8,18)$   &  $42^1 2^{99} (-12)^{20}$ & $\Delta_{120}=$ inside $\Delta_{176}$ \\ \hline
      $(162,56,10,24)$  &  $56^1 2^{140} (-16)^{21}$ & $\Delta_{162}=$ inside $\Delta_{275}$ \\ \hline
      $(176,70,18,34)$  &  $70^1 2^{154} (-18)^{22}$ & $\Delta_{176}=$ inside $\Delta_{275}$ \\ \hline
      $(253,112,36,60)$ &  $112^1 2^{230} (-26)^{22}$ & $\Delta_{253}=$ $M_{23}$-graph \\ \hline
      $(275,112,30,56)$ &  $112^1 2^{252} (-28)^{22}$ & $\Delta_{275}=$ McLaughlin graph \\ \hline
      $(276,140,58,84)$ & $140^1 2^{252} (-28)^{23}$ & Not unique (McLaughlin $2$-graph) \\ \hline
    \end{tabular}
    \vspace*{1mm}
    \caption{Some exceptional strongly regular graphs with $r=2$.}
    \label{tab:r=2}
\end{table}

\begin{definition}[Exceptional Graphs]\label{def:except-graph}
  The exceptional parameter sets for strongly regular graphs with
  $r=2$ are known, but a full classification of the corresponding
  graphs is not yet complete (Remark~\ref{rem:structure}).  For
  roughly half of these parameter sets, strongly regular graphs have
  been found and many of these describe a unique graph up to
  isomorphism.  The first twelve lines of Table~\ref{tab:r=2} are
  exceptional parameter sets that describe a unique graph and we write
  $\Delta_v$ for the unique graph in this list with $v$ vertices.  See
  \cite{BrouwerVanMaldeghem22} for descriptions of the graphs and
  references to the literature.
\end{definition}

The McLaughlin graph is listed in the penultimate line of
Table~\ref{tab:r=2}.
  
\begin{defn}[The McLaughlin graph]\label{def:mcl-graph}
  The \emph{McLaughlin graph} $\Delta = \Delta_{275}$ is the unique
  strongly regular graph with parameters $(275,112,30,56)$ and its
  spectrum is $\spec(\Delta) = \{112^1 \; 2^{252} \;(-28)^{22}\}$.  Its
  automorphism group $\aut(\Delta)$ contains the McLaughlin sporadic
  finite simple group $\McL$ as an index~$2$ subgroup.  The size of
  $\McL$ is $2^7.3^6.5^3.7.11$.  See \cite[p.~292]{ConwaySloane99}.
\end{defn}

All graphs in Table~\ref{tab:r=2} have Euclidean representations in
the Leech lattice.

\begin{defn}[Leech Lattice]\label{def:leech}
  The \emph{Leech lattice} $\Lambda$ is the unique even self-dual
  $\Z$-lattice in dimension $24$ with no vectors of norm $2$,
  i.e. $\Lambda(2)=\emptyset$.  It has \emph{standard coordinates} in
  $\tfrac{1}{\sqrt{8}} \Z^{24}$, but with the scale factor
  $\tfrac{1}{\sqrt{8}}$ is usually left implicit and suppressed.  The
  $196\ 560$ vectors in $\Lambda(4)$ have three `shapes'.  There are
  $2^{12}\cdot 24 = 98\ 304$ norm~$4$ vectors $(\mp 3,\pm 1^{23})$
  with the $\mp 3$ in any position, and the upper signs in a code word
  of the extended binary Golay code.  There are also $2^2\cdot
  \binom{24}{2} = 1\ 104$ norm~$4$ vectors of shape $((\pm 4)^2
  0^{22})$.  Finally, there are $2^7 \cdot 759 = 97\ 152$ norm~$4$
  vectors of shape $((\pm 2)^8 0^{16})$ with an even number of minus
  signs and where the locations of the non-zero coordinates form an
  octad of the $S(5,8,24)$ Steiner system used to construct the
  extended binary Golay code.
\end{defn}

For a combinatorial construction of the McLaughlin graph using the
$S(4,7,23)$ Steiner system, see \cite[\S
  10.61]{BrouwerVanMaldeghem22}. Here is a construction inside the
Leech lattice.

\begin{defn}[McLaughlin inside Leech]\label{def:mcl-in-leech}
  Let $u_0, u_1, u_2 \in \Lambda$ be Leech vectors that form a
  triangle with side lengths $2$, $2$, and $\sqrt{6}$, and let
  $\Delta_{275}$ be the set of $275$ Leech vectors distance $2$ from
  $u_i$ for $i=0,1,2$.  For all $x \neq y \in \Delta_{275}$, the
  vector $x-y$ has norm $4$ or norm $6$ (see Lemma~\ref{lem:2-norms}).
  The \emph{McLaughlin graph} has vertex set $\Delta_{275}$ and there
  is an edge $\{x,y\}$ if and only if $x-y$ has norm $6$.
\end{defn}

For explicit computations, we choose $u_0$, $u_1$, and $u_2$ as in
\cite{Wilson89}.

\begin{remark}[Explicit coordinates]\label{rem:coordinates}
  Fix $u_0 = (0^{24})$, $u_1 = (4^2\ 0^{22})$, and $u_2 =
  (3\ (-1)^{23})$ as our standard triangle, where the exponents
  indicate repeated consecutive coordinate values. Thus, $u_1=4 e_1 +
  4 e_2$.  The vectors $u_1-u_0$ and $u_2-u_0$ have norm $4$, and one
  can check that $u_2-u_1$ has norm $6$.  Let
  $\Delta=\Delta_{275}$ be
  \begin{equation*}
    \Delta = \{ x \in \Lambda(4) \colon (x - u_1)^2 = (x -
    u_2)^2= 4 \}.
  \end{equation*}
  The $275$ vectors in $\Delta$ split into three groups.  There are
  $22$ vectors of the form $(4\ 0\ (-4)\ 0^{21})$.  There are $176$
  vectors of the form $(3\ 1\ (-1)^{15}\ 1^7)$.  There are $77$
  vectors of the form $(2\ 2\ (-2)^6\ 0^{16})$.  In all three cases,
  the first two entries are fixed, while the remaining $22$ entries
  are permuted to form a Leech vector.  This split corresponds to the
  $22+176+77$ combinatorial construction of the McLaughlin graph
  \cite[\S 10.61]{BrouwerVanMaldeghem22}.
\end{remark}

\begin{remark}[Sphere containing $\Delta$]\label{rem:sphere}
  Let $U$ be the $2$-dimensional subspace containing the triangle with
  vertices $u_i$ and let $c$ be the circumcenter of this triangle, the
  unique point in $U$ equidistant from the $u_i$.  One can verify that
  the point $c=\tfrac15(u_0+2u_1+2u_2)$ and the vectors $u_i-c$ have
  norm $\tfrac85$.  And since $u_0$ is the origin (and $U$ is a linear
  subspace) the points equidistant from the $u_i$ lie in the affine
  subspace $U^\bot + c$.  In particular, the $275$ vectors $x-c$ with
  $x \in \Delta$ have norm $\tfrac{12}{5}$, so the $x$'s lie in a
  $21$-dimensional sphere centered at $c$ in $U^\bot+c$.
\end{remark}

One can verify by computer, or case-by-case, that the Leech vectors in
$\Delta$ are a $2$-distance set in this sphere, but one can also use
the geometry of the Leech lattice.

\begin{lemma}[Norms]\label{lem:2-norms}
  If $x, y\in \Delta$ and $x \neq y$, then the norm of $x-y$ is $4$ or
  $6$.
\end{lemma}

\begin{proof}
  Since the distance from $x$ to $y$ is at most twice the distance to
  $c$, the norm of $x-y$ is at most $4\cdot \tfrac{12}{5} =
  \tfrac{48}{5} < 10$.  In particular, as Leech vectors, the norm of
  $x-y$ is $4$, $6$, or $8$.  If $x-y$ has norm $8$, then the midpoint
  $m=\tfrac12(x+y)$ is a deep hole of the Leech lattice of type
  $A_1^{24}$ \cite[Theorem~8 in Chapter~23]{ConwaySloane99}, meaning
  there are $48$ Leech vectors $z$ at distance $\sqrt{2}$ from $m$ and
  the norm $2$ vectors $z-m$ are pairwise opposite or orthogonal.
  Since the triangle with vertices $u_i$, $x$, and $y$ is an isosceles
  right triangle with a right angle at $u_i$, all three vertices are
  equidistant from $m$, the midpoint of the hypotenuse.  Thus $u_1$
  and $u_2$ are Leech vectors at distance $\sqrt{2}$ from $m$, but the
  vectors $u_1-m$ and $u_2-m$ are neither opposite nor orthogonal
  since their inner product is $-1$, contradiction.  Thus $x-y$ cannot
  have norm $8$.
\end{proof}

\begin{figure}
  \centering
  \includegraphics[width=0.8\linewidth]{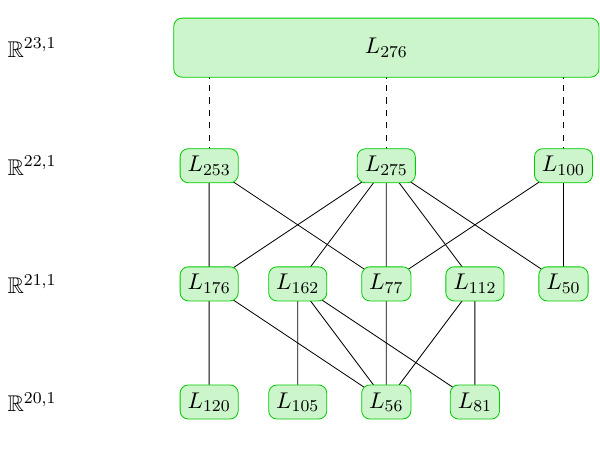}
  \caption{Exceptional lattices with $r=2$ under inclusion.}
  \label{fig:r=2}
\end{figure}

We now shift from Euclidean representations to Lorentzian lattices.

\begin{definition}[Exceptional Lattices]\label{def:except-lattices}
  We write $L_v$ for the Lorentzian lattice of the graph $\Delta_v$
  listed in Table~\ref{tab:r=2}.  By Theorem~\ref{thm:latt-functor}
  inclusion relations among these graphs result in inclusion relations
  between their Lorentzian lattices. The lines in Figure~\ref{fig:r=2}
  are lattice inclusions that result from subgraph inclusions that we
  found in the literature.  There may be other graph inclusions that
  we have missed. All twelve of these examples produce
  lattices in $\R^{n,1}$ where $n$ is between $20$ and $22$.  
\end{definition}

\begin{remark}[McLaughlin $2$-graph]\label{rem:mcl-2-graph}
  The final line of Table~\ref{tab:r=2} is the parameter set
  $(276,140,58,84)$ with spectrum $\{140^1 \; 2^{252} \;
  (-28)^{23}\}$.  Nozaki has shown that there are at least ten
  thousand non-isomorphic strongly regular graphs with these
  parameters \cite[\S 10.61]{BrouwerVanMaldeghem22}, but all of these
  were found by process known as \emph{Seidel-switching}.  The finite
  collection of graphs switching equivalent to a given graph is called
  a \emph{$2$-graph}. Let $L_{276}$ be the Lorentzian lattice
  constructed from any one of these graphs.  One can show, in this
  example, that the Seidel-switching process can be carried out
  geometrically inside the lattice $L_{276}$, so that the ten thousand
  non-isomorphic graphs construct the same lattice $L_{276}$ despite
  having non-equivalent spanning sets.  The statement made above about
  the lattice remaining unchanged is not true in general.  There exist
  distinct strongly regular graphs that are Seidel-switching
  equivalent, but which produce non-isomorphic lattices.  We hope to
  explore issues related to Seidel-switching in the context of
  Lorentzian representations of graphs (e.g. in $L_{276})$, in a
  future article.  See also Remark~\ref{rem:mcl-2-graph-2}.  In
  Figure~\ref{fig:r=2}, the unique lattice $L_{276}$ in $\R^{23,1}$ is
  listed at the top.  The inclusion lines from $M_{23}$ lattice
  $L_{253}$, the McLaughlin lattice $L_{275}$ and the Higman-Sims
  lattice $L_{100}$ are dashed because their spanning sets include
  into different spanning sets.
\end{remark}

It turns out that all of the lattices in Figure~\ref{fig:r=2} are
sublattices of $\II_{25, 1}$. We show this explicitly for the
$22$-dimensional McLaughlin lattice $L_{275}$ and explain the related
hyperbolic geometry inside the Weyl chamber of $\II_{25, 1}$.  The
following description of the simple roots and fundamental domain of
the reflection group of $\II_{25,1}$ is due to Conway and Borcherds.
See \cite{Borcherds85-leech}.

\begin{defn}[Leech roots and the Leech cusp]\label{def:leech-roots}
  The lattice $\II_{25,1}$ can be viewed as the orthogonal direct sum
  $\II_{25,1} = \Lambda \bot H$, of the Leech lattice $\Lambda$ in
  $\R^{24}$ and a hyperbolic cell $H = \Z \rho_1 + \Z \rho$ where
  $\rho_1$ and $\rho$ are null vectors with $\ip{\rho_1}{ \rho} = -1$.
  We write the vectors of $\II_{25,1}$ in the form
  \[
  z = x + m \rho_1 + n \rho = (x ; m, n) \text{\; where \;} 
  x \in \Lambda 
  \text{\; and \;}
  m, n \in \Z.
  \]
  Note that $z^2 = x^2 - 2 m n$.  In these coordinates $\rho = (0; 0,
  1)$ and $\proj \rho$ is in the boundary of $\hyp^{25}_\R$.  The
  vector $\rho$ (and the point $\proj \rho$) is called a \emph{Leech
  cusp}.  The norm $2$ vectors of $\II_{25, 1}$ are the \textit{roots}
  of $\II_{25, 1}$ and the hyperplanes orthogonal to them (in $\R
  \II_{25, 1}$ or in $\hyp^{25}_\R$) are the \textit{mirrors} of
  $\II_{25, 1}$.  The roots $\alpha \in \II_{25,1}$ satisfying
  $\ip{\alpha}{ \rho} = -1$ are in bijection with the vectors of the
  Leech lattice.  For each $x \in \Lambda$, define a \emph{Leech root}
  \[
  \alpha_{x, \Lambda} := (x; 1, \tfrac{1}{2}x^2 - 1).
  \]
  These are precisely the roots in $\II_{25,1}$ with
  $\ip{\alpha}{\rho}=-1$.  In particular, each norm~$4$ vectors $w \in
  \Lambda(4)$ determines a Leech root $\alpha_{w, \Lambda} = (w; 1,
  1)$.  It is useful fact that for all $x, y \in \Lambda$,
  \begin{equation}
  \ip{\alpha_{x, \Lambda}}{\alpha_{y, \Lambda}}
  = 2 - \tfrac{1}{2}(x - y)^2.
  \label{eq:Leech-root-ip}
  \end{equation}
  The Leech mirrors, i.e. the mirrors of the Leech roots, are in bijection
  with the walls of the domain
  \[
  D = \{ \proj z \in \hyp^{25}_\R \colon
  \ip{z}{ \alpha_{x, \Lambda}} \leq 0 
  \text{\; for all \;} x \in \Lambda \}.
  \]
  This is a fundamental domain for the action of the reflection group
  of $\II_{25, 1}$ ``around" the cusp $\rho$.  Even though $\rho$ is
  in the closure of this fundamental domain $D$, there are no mirrors
  of $\II_{25,1}$ passing through $\rho$.  This is a consequence of
  the fact that the Leech lattice has no norm $2$ vectors.
  Geometrically, the Leech mirrors are precisely the mirrors closest
  to $\rho$ according to horocyclic distance, and they are all tangent
  to the unit horoball centered at $\rho$.
\end{defn}  

\begin{remark}[Reflection Group of $\II_{25,1}$]\label{rem:refl-even-25-1}
  The Leech roots $\lbrace \alpha_{x, \Lambda} \colon x \in \Lambda
  \rbrace$ index reflections $\{R_x \colon x\in \Lambda \}$ that
  generate $\Refl(\II_{25,1})$ and $\Refl(\II_{25,1})$ can be viewed
  as an infinite rank Coxeter group with these reflections as its
  ``simple system''\footnote{Technically the Leech roots merely form a
  linearly dependent \emph{root basis}, but this is enough to
  establish that the representation of the Coxeter group is faithful.}
  and only three types of relations: (1) the generators $R_x$ have
  order $2$, (2) the product $R_x R_y$ has order $2$ when the $x-y$
  has norm $4$, and (3) the product $R_xR_y$ has order $3$ when $x-y$
  has norm $6$.  See \cite[Chapter 27]{ConwaySloane99}.  Conway's
  aphorism for this situation is: ``The Leech lattice is the Dynkin
  diagram of $\II_{25,1}$''.
\end{remark}

\begin{definition}[McLaughlin Roots]
  The $275$ norm~$4$ vectors in $\Delta$
  (Remark~\ref{rem:coordinates}) corresponding to $275$ Leech roots of
  the form
  \begin{equation*}
    \spanning' = \{ \alpha_{x, \Lambda} 
    = (x; 1,1) \colon x \in \Delta \}.  
  \end{equation*}
  We call the elements in $\spanning'$ \emph{McLaughlin roots}, their
  mirrors \emph{McLaughlin mirrors}, and we write $L'$ for the
  submodule of $\II_{25,1}$ generated by $\spanning'$.
\end{definition}

\begin{proposition}[Sublattice]\label{prop:sublattice} 
  The spanning sets in $\Z_\Delta$ and in $L'$ have the same 
  Gram matrix: $\gram(\spanning_\Delta) = \gram(\spanning')$.
  So the map that sends $\alpha_x \mapsto \alpha_{x, \Lambda}$ 
  defines an isomorphism $\Z_{\Delta} \to L'$ taking $S_{\Delta}$ to $S'$. 
\end{proposition}

\begin{proof} 
  Equation \eqref{eq:Leech-root-ip} and
  Definition~\ref{def:mcl-in-leech} implies that
  $\gram(\spanning_\Delta) = \gram(\spanning')$. It follows that there
  is a natural map $\Z^\Delta \to L'$ given by $e_x \mapsto \alpha_{x,
    \Lambda}$ for $x \in \Delta$. As in the proof of
  Theorem~\ref{thm:latt-functor}, $L'$ is a submodule of a Lorentzian
  space and the image of $\one \in \Z^\Delta$ is a vector of negative
  norm, so $L'$ is non-degenerate (Remark~\ref{rem:negative-norm}) and
  the map $\Z^{\Delta} \to L'$ induces an isometry $\Z_\Delta \to L'$
  such that $\alpha_x \mapsto \alpha_{x,
    \Lambda}$(Proposition~\ref{prop:maps-onto}).
\end{proof}

Based on Proposition~\ref{prop:sublattice}, we identify $\Z_\Delta$
with the sublattice $L'$ of $\II_{25,1}$.  In particular we identify
$\spanning_\Delta$ with $\spanning'$ and $\alpha_{x, \Lambda}$ with
$\alpha_x$.  The next result shows that $\Z_\Delta$ is the orthogonal
complement of a positive definite sublattice $M$ inside $\II_{25,1}$.
Let $M$ be the $\Z$-span of the Leech roots $\alpha_j := \alpha_{u_j}$
for $j = 0, 1, 2$, where $u_j$'s are as in
Definition~\ref{def:mcl-in-leech}.  In the coordinates of
Remark~\ref{rem:coordinates}
  \[
    \alpha_0 = (0; 1, -1), \;\;
    \alpha_1 = (u_1; 1, 1), \; 
    \text{ and \;}
    \alpha_2 = (u_2; 1, 1).
  \]
The sublattice $M$ is a root lattice of type $A_1 \sqcup A_2$ and
it has discriminant~$6$. 

\begin{proposition}[Complement]\label{prop:mcl-complement} 
  A Leech root $\alpha_{x}$ is a McLaughlin root if and only if
  $\alpha_{x} \in M^\bot$.  In fact, $\Z_\Delta = M^\bot$ and the
  $\alpha_x \in \spanning_\Delta$ are the Leech roots in $\Z_\Delta$.
\end{proposition}

\begin{proof}
  Let $\alpha_x = (x; 1, *)$ be a Leech root and recall that
  $\ip{\alpha_{x}}{ \alpha_j} = 2 - \tfrac{1}{2} ( x - u_j)^2$.  We
  have $\alpha_x \in M^{\bot}$ if and only if $\ip{\alpha_{x}}{
    \alpha_j} = 0$ for $j = 0, 1, 2$ if and only if $(x - u_j)^2 = 4$
  for $j = 0, 1, 2$ if and only if $x \in \Delta$, i.e. $\alpha_x$ is
  a McLaughlin root.  Since the McLaughlin roots $\spanning_\Delta$
  are in $M^{\bot}$ and they span $\Z_{\Delta}$, we have $\Z_{\Delta}
  \subseteq M^{\bot}$. Since $\II_{25,1}$ is self-dual and the
  sublattice $M$ has discriminant $6$, its orthogonal complement
  $M^{\bot}$ has discriminant $-6$.  By definition $\Z_{\Delta}$ is
  spanned by the $275$ McLaughlin roots.  We ran the LLL-reduction
  algorithm on these $275$ vectors to obtain a basis of $\Z_{\Delta}$.
  We then verified that the Gram matrix of this basis has determinant
  $-6$, which proves that $\Z_{\Delta} = M^{\bot}$.
\end{proof}

\begin{definition}[Weyl-like Vector $\tau_\Delta$]
  Let $\alpha_\Delta$, $\tau$, and $\tau_\Delta$ be defined as in
  Definition~\ref{def:tau}.  In the coordinates of
  Remark~\ref{rem:coordinates}, the vector $\tau_\Delta = (14, 6,
  (-2)^{22} ; 5, 5)$ in $\Z_\Delta$ has norm~$-10$.  In addition,
  $\alpha_\Delta = 55 \cdot \tau_\Delta$ and $\tau = 28\cdot
  \tau_\Delta$.
\end{definition}

\begin{prop}[Unique Fixed Point]\label{prop:mcl-tau}
  The vector $\tau_\Delta$ has the following properties:
  \begin{enumerate}
  \item[(a)] $\ip{\tau_\Delta}{ \alpha_{x}} = -2$ for all $\alpha_x
    \in \spanning_\Delta$.
  \item[(b)] $\proj \tau_\Delta$ is equidistant from the $275$
    McLaughlin mirrors.
  \item[(c)] $\proj\tau_\Delta$ is the unique point in $\ball \R_\Delta
    \simeq \hyp^{22}_\R$ fixed by $\McL$.
  \item[(d)] $\tau_\Delta = 2\rho + (\alpha_0 + 2\alpha_1 +
    2\alpha_2)$.
  \item[(e)] $\proj\tau_\Delta$ is the projection of the Leech cusp
    $\proj\rho$ onto $\ball \R_\Delta$.
  \item[(f)] $\proj\tau_\Delta$ belongs to the facet of the boundary
    of Weyl chamber $D \cap \hyp^{22}_\R$.
\end{enumerate}
\end{prop}

\begin{proof}
  For (a), $\ip{\tau}{\alpha_x}=-\mu = -56$ by
  Lemma~\ref{lem:elts-in-rad}, so $\ip{\tau_\Delta}{\alpha_x} =
  -56/28=-2$.  Part (b) is immediate from (a).  For (c), $\McL \subset
  \aut(\Delta)$ fixes $\proj\tau_\Delta = \proj \alpha_\Delta$ by
  Proposition~\ref{prop:weyl-like}, and we get an induced non-trivial
  Euclidean representation of $\McL$ on $\R^{22}$.  If there were
  additional fixed points, that would induced a non-trivial
  representation of $\McL$ on a proper subspace of $\R^{22}$,
  contradicting the fact that the smallest non-trivial representation
  of the McLaughlin simple group is $22$-dimensional.  Part (d) can be
  verified via the explicit coordinates. One can also argue as
  follows.  Let $u = 2 \rho + (\alpha_0 + 2 \alpha_1 + 2 \alpha_2)$.
  Since $\ip{u}{\alpha_j} = 0$ for $j = 0, 1, 2$, we have $u \in
  M^{\bot} = \Z_{\Delta}$. Furthermore, for each $x \in \Delta$, one
  has $\ip{\alpha_x}{\alpha_j} = 0$, so $\ip{u}{\alpha_x} = \ip{2
    \rho}{\alpha_x} = - 2$. By part (a), $u - \tau_{\Delta}$ is in
  $\Z_{\Delta}$ and is orthogonal to each $\lbrace \alpha_x \colon x
  \in \Delta \rbrace$, so $u - \tau_{\Delta} = 0$. This proves (d).
  Next, since $(\alpha_0 + 2\alpha_1 + 2\alpha_2) \in M$ and
  $\tau_\Delta \in M^\bot$, the projection of their difference $2\rho$
  to $\R M^\bot$ is $\tau_\Delta$.  This proves (e).  Finally, note
  that if $\alpha_{x}$ is any Leech root distinct from $\alpha_0,
  \alpha_1, \alpha_2$, then
  \begin{equation}
    \ip{\tau_\Delta}{\alpha_{x}} = \ip{2 \rho}{\alpha_{x}} + 2
    \ip{\alpha_1}{ \alpha_{x}} + 2 \ip{\alpha_2}{ \alpha_{x}} +
    \ip{\alpha_0}{ \alpha_{x}} \leq -2
    \label{eq-rlambdatau}
  \end{equation}
  since $\ip{\rho}{\alpha_{x}} = -1$ and $ \ip{\alpha_j}{ \alpha_{x}}
  \leq 0$. This proves (f).
\end{proof}

The statement about Leech roots can be extended to all root of
$\II_{25,1}$.

\begin{proposition}[Closest Mirrors]\label{prop:no-mirrors-McL}
  There are no mirrors of $\Z_{\Delta}$ passing through the point
  $\tau_\Delta$.  The $275$ McLaughlin mirrors $\{ \alpha_x^{\bot}
  \colon x \in \Delta \}$ are precisely the mirrors of $\Z_{\Delta}$
  that are closest to $\tau_\Delta$ and they are equidistant from
  $\tau_\Delta$.
\end{proposition}

\begin{proof}
  Let $\alpha \in \Z_{\Delta}(2)$ be a root of $\Z_\Delta$. By
  Proposition~\ref{prop:mcl-complement}, $\ip{\alpha}{\alpha_j} = 0$
  for $j = 0, 1, 2$, so Proposition~\ref{prop:mcl-tau}(d) implies that
  \[
  \ip{\tau_\Delta}{\alpha} = \ip{2\rho}{\alpha}.
  \]
  Since $\alpha$ is also a root of $\II_{25,1}$ and all roots of
  $\II_{25, 1}$ have nonzero dot product with $\rho$
  (Definition~\ref{def:leech-roots}), we have
  $\ip{\tau_\Delta}{\alpha} \neq 0$, which shows that no mirrors
  contain $\tau_\Delta$.  The minimum possible value of
  $\abs{\ip{\tau_\Delta}{\alpha}}$ is $2$ and it is attained exactly
  when $\ip{\rho}{\alpha} = \pm 1$.  Thus $\alpha^{\bot}$ is closest
  to $\tau$ if and only if, up to a sign, $\alpha \in \Z_{\Delta}$ is
  a Leech root and this is equivalent to $\alpha$ being a McLaughlin
  root by Proposition~\ref{prop:mcl-complement}.
\end{proof}

\begin{remark}[McLaughlin $2$-graph, revisited]\label{rem:mcl-2-graph-2}
  The lattice $L_{276}$ (Remark~\ref{rem:mcl-2-graph}) is the
  sublattice of $\II_{25,1}$ generated by finite set of the Leech
  roots orthogonal to Leech roots $\alpha_x$ and $\alpha_y$ where
  $x-y$ has norm $6$.  There are exactly $552$ Leech roots in
  $L_{276}$ which come in pairs.  If $\Delta$ is one of the $10,000$
  strongly regular graphs with parameters $(276,140,58,84)$, then its
  Lorentzian lattice is isometric to $L_{276}$ and its spanning set is
  a choice of $276$ of these $552$ Leech roots, one from each pair.
\end{remark}

The construction of $L_{275}$ inside $\II_{25,1}$ can be placed in a
more general context.

\begin{remark}[General Procedure]
  Let $A$ be a set of Leech roots that form a simple system of a
  finite Coxeter group, let $M \subseteq \II_{25,1}$ be the positive
  definite rank~$k$ sublattice they span and let $M^{\bot}$ be the
  orthogonal sublattice in $\II_{25,1}$ of signature $(25 - k, 1)$.
  Let $\Delta = \{ x \in \Lambda \colon \alpha_x \in M^{\bot} \}$.
  Let $\rho_A$ and $\rho_{\Delta}$ be the linear projections of $\rho$
  on $M$ and $M^{\bot}$ respectively.  Then
  \[
  \rho = \rho_A + \rho_{\Delta}
  \]
  The vector $\rho$ is called the Weyl vector of the fundamental
  domain $D$\footnote{Indeed, there is a important Kac-Moody algebra,
  called the fake Monster Lie algebra, for which the root lattice is
  $\II_{25, 1}$, a set of simple roots is the set of Leech roots and
  the Weyl vector is $\rho$.}.  Since $\rho$ has inner product $-1$
  with all the Leech roots, it follows that $\rho_A$ is a Weyl vector
  for the simple system $A$. Thus $\rho_{\Delta}$ has inner product
  $-1$ with the Leech roots $\{ \alpha_x \colon x \in \Delta \}$ and
  should be called the Weyl vector of $(\Z_{\Delta}, \{ \alpha_x \colon
  x \in \Delta \})$.  When $A=\{\alpha_0,\alpha_1,\alpha_2\}$ of type
  $A_2 \bot A_1$, the graph $\Delta \subseteq \Lambda$ 
    is the McLaughlin
  graph, $(\Z_{\Delta}, \{ \alpha_x \colon x \in \Delta \})$ is the
  Lorentzian representation of $\Delta$ realized inside $\II_{25, 1}$
  and $\rho_{\Delta} = \tfrac{1}{2} \tau_\Delta$.  So $\tfrac12
  \tau_\Delta$ is the analog of Weyl vector.  A similar description
  holds for Lorentzian representations of some of the other
  exceptional strongly regular graphs listed in Table~\ref{tab:r=2}.
  When $A = A_3$, then $\Delta$ is the Higman-Sims and the entire 
  discussion in this section can be repeated in this context.
  When $A =
  A_2$, then $\Delta$ is the collection of $552$ roots that form the
  McLaughlin $2$-graph.
\end{remark}

\section{Complex Hermitian lattices from bipartite graphs}\label{sec:bipartite}

In this section the spanning set of a real Lorentzian lattice
constructed from a bipartite graph is modified to create a spanning
set for a complex Lorentzian lattice.  Roughly speaking, all of the
real spanning vectors are rescaled to have integer norms and the
spanning vectors indexed by one class of vertices in the bipartite
graph are rotated by a root of unity.  This gives a complex analog of the 
constructions in \S\ref{sec:graph-lattice}.
  
\begin{defn}[Bipartite graphs]\label{def-directed-bipartite-graphs}
  A \emph{bipartite graph} is a graph $\Delta$ with a non-trivial
  vertex partition $\Delta = \Delta^- \sqcup \Delta^+$ such that every
  edge of $\Delta$ contains one vertex in $\Delta^-$ and one vertex in
  $\Delta^+$.  We call $\Delta^-$ and $\Delta^+$ the \emph{parts} of
  $\Delta$ and set $v_\pm = \abs{\Delta^\pm}$.  Since many of our
  examples come from incidence graphs of point-line configurations, we
  call the vertices in $\Delta^-$ \emph{points} and the vertices in
  $\Delta^+$ \emph{lines}.  Equivalently, a bipartite graph is a graph
  with a surjective \emph{sign function}
  \[
  \sgn: \Delta \onto \{-1,+1\}
  \]
  where the induced subgraphs on $\Delta^{\pm} := \sgn^{-1}(\pm 1)$
  have no edges.  When $\Delta$ is connected, the partition of
  $\Delta$ into parts and the corresponding function $\sgn$ are unique
  up to switching the $\pm$ signs.  We write $\vecDelta =
  (\Delta,\sgn)$ for a \emph{directed} bipartite graph where the
  labeling of the parts is fixed.  A \emph{directed edge} $(i,j)$ is
  an edge $\{i,j\}$ of $\Delta$ with $i\in \Delta^-$ and $j \in
  \Delta^+$ and we write $\edg_\vecDelta$ for the set of directed
  edges.  The group $\aut(\vecDelta)$ contains the automorphisms of
  $\Delta$ that send $\Delta^+$ to $\Delta^+$ and $\Delta^-$ to
  $\Delta^-$.  If we fix a linear order on $\Delta$ so that the points
  in $\Delta^-$ come before the lines in $\Delta^+$, the adjacency
  matrix of $\Delta$ has the form
  \begin{equation*}
    A = A_\Delta = \begin{bmatrix} 0 & N \\ N^\tr & 0 \end{bmatrix}
  \end{equation*}
  where $N = \llbracket N_{m \ell} \rrbracket$ is the $v_-
  \times v_+$ matrix with $N_{m \ell} =1$ if $(m, \ell)
  \in \edg_\vecDelta$ is a directed edge of $\Delta$ and $N_{ m \ell}
  = 0$ otherwise.  As usual, if $A$ is $k$-regular then the all 
  $1$'s column vector $\one$ is a $k$-eigenvector.
  Let $\bigraph$ be the category of connected directed bipartite
  graphs $\vecDelta = (\Delta, \sgn)$ whose arrows are sign-preserving
  inclusions $\iota\colon \vecDelta' \into \vecDelta$ of connected
  induced subgraphs. For each $r>0$ let $\bigraph(r)$ be the full
  subcategory where we restrict to graphs with $\lambda_2=r$.
\end{defn} 
  
In many examples the parts have the same size, but this is not
required.  The general construction uses constants $r$, $\zeta$, $c$
and a ring $\cO$: the constant $r$ comes from $\vecDelta$, the root of
unity $\zeta$ needs to be chosen, $c$ is defined from $r$, and $\cO$
is defined using $r$ and $\zeta$.

\begin{defn}[$r$, $\zeta$, $c$, and $\cO$]\label{def:constants}
  Let $r>0$ be a real algebraic integer such as $\sqrt{2}$ or
  $\sqrt{3}$, let $\zeta$ be a root of unity, and let $\zeta'=
  \overline{\zeta}=\zeta^{-1}$ be its conjugate and
  inverse.\footnote{We use $\zeta'$ for visual clarity: the bar on
  $\overline{\zeta}$ can be hard to see and the exponent of
  $\zeta^{-1}$ needs spacing that obscures symmetries in the
  formulas.}  The real lattice construction uses the ring of integers
  $\cO(r)$ of the number field $\Q(r)$.  The complex lattice
  construction uses the ring of integers $\cO = \cO(r \zeta, r
  \zeta')$ of the number field $\Q(r\zeta,r\zeta')$.\footnote{When
  $\zeta=\pm 1$, $\Q(r\zeta,r\zeta') = \Q(r)$ and $\cO(r \zeta, r
  \zeta')$ reduces to $\cO(r)$.  In particular, the real construction
  is a special case of the complex construction. To produce complex
  examples, we choose $\zeta \neq \pm 1$.} For simplicity, we assume
  that $\cO$ is a PID unless we state otherwise.  The number $\zeta$
  is the \emph{twisting constant}. There is also a \emph{stretching
  constant} $c \in \cO(r)$ selected so that $n=rc$ is the smallest
  positive integer in $r\cO(r)$.  In other words, $c = r^{-1} \min \{
  r \cO(r) \cap \N \}$.\footnote{The intersection is non-empty.  Let
  $p(x) \in \Z[x]$ be the minimal monic polynomial with $p(r)=0$ and
  write $p(x) = x q(x) + p(0)$. Note $p(0) \neq 0$ by the minimality
  of $p$.  For $q(r) \in \cO(r)$ we have $r q(r) = -p(0)$, so the
  abelian group $r\cO(r) \cap \Z$ contains all integer multiples of
  $p(0)$.} Because $r$ and $n$ are positive real, so is $c$, and
  $\sqrt{c}$ is a well-defined.  Moreover, since $rc=n \in \N$, we
  have $\Q(r)=\Q(c)$, and since $(r\zeta)(c\zeta') = (r\zeta')(c\zeta)
  = n$ we have $\Q(r \zeta, r\zeta') = \Q(c \zeta, c \zeta')$.  Note
  that all four elements $r\zeta$, $r\zeta'$, $c\zeta$, and $c\zeta'$
  are in $\cO$: they are in the number field and they are products of
  the algebraic integers $r$, $c$, $\zeta$, and $\zeta'$. For the
  remainder of the section let $\vecDelta = (\Delta, \sgn)$ be a
  directed connected bipartite graph with second largest eigenvalue
  $r=\lambda_2>0$ and let $\zeta \neq \pm 1$ be a root of unity so
  that $\cO = \cO(r \zeta, r \zeta')$ is a PID.  In addition, $g = v -
  f - 1$ where $v=\abs{\vecDelta}$ and $f$ is the multiplicity of $r$.
\end{defn}

The construction uses a new basis in the vector space $\C^\Delta$.
  
\begin{definition}[New Basis]\label{def:new-basis}    
  Recall that $\R^\Delta$ and $\C^{\Delta}$ has basis $\{e_i \colon i
  \in \Delta\}$.  We define a new basis $\{\vec{e}_i \colon i \in
  \vecDelta\}$ of $\C^\Delta$ by
  \begin{equation}
  \label{eq:vec_e_i}
     \vec{e}_i = 
     \begin{cases}
         \sqrt{c}\cdot e_i & \textrm{if\;} i \in \Delta^- \\
         \sqrt{c}\cdot\zeta\cdot e_i & \textrm{if\;} i \in \Delta^+ \\
     \end{cases}
  \end{equation}
  In other words, the vectors are scaled by $\sqrt{c}$ and those
  indexed by vertices in $\Delta^+$ are rotated by $\zeta$.  Note that
  $\C \vec{e}_i = \C e_i$ for all $i\in \Delta$, so $\{\vec{e}_i\}$ is
  also a basis of $\C^\Delta$.
\end{definition}

 The general construction given below has close parallels with
 \S\ref{sec:graph-lattice}.

\begin{definition}[Hermitian Spaces]\label{def:c-spaces}
  Recall that the symmetric bilinear form on $\R^\Delta_r$ is defined
  by the Gram matrix $\gram(\{e_i\}) = rI-A$
  (Definition~\ref{def:shift-module}). This matrix can also be used to
  define a Hermitian form $\ip{\;}{\;}$ on the complex vector space
  $\C^\Delta$. We write $\C^\vecDelta_r$ for the space $\C^\Delta$
  with this degenerate Hermitian form.
  \footnote{There is a subtle difference between $\C^\vecDelta_r$ and
  $\C^\Delta_r$ defined in \S\ref{sec:graph-lattice}. They are the
  same $\C$-vector space with the same basis $\lbrace e_i \colon i \in
  \Delta \rbrace$. The space $\C^{\Delta}_r$ has a symmetric bilinear
  form and the space $\C^{\vecDelta}_r$ has a hermitian form, both
  defined by the Gram matrix $(r I - A)$ for the basis $\lbrace e_i
  \colon i \in \Delta \rbrace$.}
  The signature of $\C^\vecDelta_r$ is the signature of $\R^\Delta_r$
  which is $(g,1,f)$, and the non-degenerate quotient $\C_\vecDelta
  =\C^\vecDelta_r / \rad(\C^\vecDelta_r)$ has signature $(g,1)$
  (Remark~\ref{rem:rank-signature}).  We write $\alpha_i$ for the
  image of $e_i$ in $\R_\Delta$, $\vec{\alpha}_i$ for the image of
  $\vec{e}_i$ in $\C_\vecDelta$, and $\spanning_{\vecDelta} = \lbrace
  \vec{\alpha}_i \colon i \in \Delta \rbrace$ for the spanning set.
  For the Hermitian form $\ip{\;}{\;}$ on $\C_\vecDelta$ we have the
  complex version of Equation~\ref{eq:rI-A} in
  Definition~\ref{def:shift-module}:
  \begin{equation*}
    \ip{\vec{e}_i}
       {\vec{e}_j} = 
    \ip{\vec{\alpha}_i}
       {\vec{\alpha}_j} = 
       \begin{cases}          c\cdot r = n & \text{\; if } i  = j, \\
         -c\cdot\zeta & \text{\; if } 
         i \in \Delta^-, \; j \in \Delta^+, \text{\; and \;} \lbrace i, j \rbrace \in \edg_{\Delta},
         \\
         -c\cdot\zeta' &\text{\; if } 
         j \in \Delta^-, \; i \in \Delta^+, \text{\; and \;} \lbrace i, j \rbrace \in \edg_{\Delta},
         \\
         0  & \text{\; otherwise}. 
       \end{cases}
  \end{equation*}
  Thus $\gram(\{\vec{e}_i\}) =\gram(\{\vec{\alpha}_i\})$ is
  \[
    B_\zeta := c (rI-A_\zeta)
  \]
  where $A$ is the adjacency matrix, $A_{\zeta} = P_\zeta^\htr A
  P_\zeta$, and $P_\zeta$ is a diagonal matrix with diagonal entries
  $1^{v_-} \zeta^{v_+}$.  As block matrices,
  \[
  B_\zeta :=
  c\cdot 
  \begin{bmatrix}
    rI & -\zeta N\\
    -\zeta' N^\tr & rI
  \end{bmatrix}
  \textrm{,\ \ }
  A_{\zeta} :=
  \begin{bmatrix} 
    0 &  \zeta  N \\ 
    \zeta' N^{\msf{T}} & 0  
  \end{bmatrix}
  \textrm{, and\ \ }
  P_\zeta := 
  \begin{bmatrix} 
    I &  0 \\ 
    0 & \zeta   I  
  \end{bmatrix}.
  \]
  For later use we also note that
  \begin{equation}\label{eq:A-square}
  A_\zeta^\htr = A_\zeta 
  \textrm{\ \ and \ \ }
  A_\zeta^2 = A^2 = 
  \begin{bmatrix} 
  NN^\tr & 0 \\
  0 & N^\tr N
  \end{bmatrix}.
  \end{equation}
\end{definition}

Examples~\ref{ex:eisenstein} and~\ref{ex:gaussian} illustrate how 
stretching and twisting can produce Eisenstein and Gaussian matrices.  
The first calculation uses $r=c=\sqrt{3}$, $\zeta=i$ and $\cO = \eisen$, 
and the second calculation uses $r=c=\sqrt{2}$, $\zeta=\zeta_8$ and $\cO=\gauss$.

\begin{example}[$\eisen$-matrices]\label{ex:eisenstein}
  Consider $\C^2$ with a Hermitian form defined by the real matrix
  $rI-A$ where $A$ is the adjacency matrix of a single edge and
  $r=\sqrt{3}$.  Suppose we create a new basis by rescaling both basis
  vectors by $\sqrt[4]{3}$ and rotating the second basis vector by
  $i$.  The rescaling multiples the form by $\sqrt{3}$ and the
  Gram matrix for the new basis has entries in $\eisen = \Z[\omega]$:
  \[
  \sqrt{3} \cdot
  \begin{bmatrix} 
    1 &  0 \\ 
    0 & i  
  \end{bmatrix}^\htr 
  \begin{bmatrix} 
    \sqrt{3} &  -1 \\ 
    -1 & \sqrt{3}  
  \end{bmatrix}
  \begin{bmatrix} 
    1 &  0 \\
    0 & i  
  \end{bmatrix}
  =
  \sqrt{3}\cdot 
  \begin{bmatrix} 
    \sqrt{3} &  -i \\ 
    i & \sqrt{3}  
  \end{bmatrix} 
  =
  \begin{bmatrix} 
    3 &  -\theta \\ 
    \theta & 3  
  \end{bmatrix}.
  \]
  As usual $\theta = \omega - \omega^2 = i\sqrt{3}$.  The
  $\Z[\sqrt{3}]$-span of the original basis is a
  $\Z[\sqrt{3}]$-lattice in $\R^2 \subset \C^2$.  The $\eisen$-span of
  the new basis is a $2$-dimensional $\eisen$-lattice in $\C^2$.
\end{example}

\begin{example}[$\gauss$-matrices]\label{ex:gaussian}
  Consider $\C^2$ with a Hermitian form defined by the
  matrix $rI-A$ where $A$ is the adjacency matrix of a single edge and
  $r=\sqrt{2}$.  Suppose we create a new basis by rescaling both basis
  vectors by $\sqrt[4]{2}$ and rotating the second basis vector by
  $\zeta_8$, the primitive eighth root of unity.  The 
  rescaling multiples the form by $\sqrt{2}$ and the Gram matrix for
  the new basis has entries in $\gauss = \Z[i]$:
  \[
  \sqrt{2} \cdot
  \begin{bmatrix} 
    1 &  0 \\ 
    0 & \zeta_8  
  \end{bmatrix}^\htr 
  \begin{bmatrix} 
    \sqrt{2} &  -1 \\ 
    -1 & \sqrt{2}  
  \end{bmatrix}
  \begin{bmatrix} 
    1 &  0 \\
    0 & \zeta_8  
  \end{bmatrix}
  = \sqrt{2}\cdot 
  \begin{bmatrix} 
    \sqrt{2} &  \zeta'_8 \\ 
    \zeta_8 & \sqrt{2}  
  \end{bmatrix} 
  =
  \begin{bmatrix} 
    2 &  -1-i \\ 
    -1+i & 2  
  \end{bmatrix}.
  \]
  The $\Z[\sqrt{2}]$-span of the original basis is a
  $\Z[\sqrt{2}]$-lattice in $\R^2 \subset \C^2$.  The $\gauss$-span of
  the new basis is a $2$-dimensional $\gauss$-lattice in $\C^2$.
\end{example}

We are now ready to define the $\cO$-lattice of $\vecDelta$.

\begin{definition}[$\cO$-lattices]\label{def:cplx-construction} 
  Let $\cO = \cO(r \zeta, r \zeta')$ be as in Definition~\ref{def:constants},
  let  $\cO^{\vecDelta}_r = \cO \{\vec{e}_i\}$ be the free $\cO$-module with 
  the Hermitian form defined by the matrix $B_\zeta= c(r I - A_{\zeta})$ and 
  let $\cO_{\vecDelta} = \cO^{\vecDelta}_r/ \rad(\cO^{\vecDelta}_r)$.
  \footnote{The $\cO$-modules $\cO^{\vecDelta}_r$ and
  $\cO_{\vecDelta}$ depend on the choice of $\zeta$, but we leave
  $\zeta$ implicit for notational simplicity.  Since there is a
  canonical choice for $\zeta$ in each of our examples, this slight
  abuse of notation should not be a problem.}
  Since the entries of $B_\zeta$ are in $\cO$, the form on
  $\cO_\vecDelta$ takes values in $\cO$.  As in
  Remark~\ref{remark-real-change-of-rings}, using Appendix
  Lemma~\ref{lemma-change-of-rings} we obtain natural inclusions
  \[
  \cO^{\vecDelta}_r \subseteq \C^{\vecDelta}_r
  \text{\; and \;}
  \cO_{\vecDelta} \subseteq \C_{\vecDelta}.
  \]
  and $\cO_{\vecDelta}$ contains the spanning set
  $\spanning_{\vecDelta} = \lbrace \vec{\alpha}_i \colon i \in \Delta
  \rbrace$.  We say that $(\cO_{\vecDelta}, \spanning_{\vecDelta})$ is
  \textit{complex Lorentzian representation of $\vecDelta$} defined
  over the ring $\cO$. In the more general notation of
  Definition~\ref{definition-Hermitian-R-modules} one has
  \[
  \cO^{\vecDelta}_r = \cO^v_B 
  \text{\; and \;}
  \cO_{\vecDelta} = \cO_{v,B}
  \]
  where $B = B_\zeta$ and where the involution $\sigma_{\cO}$ is
  complex conjugation.  
\end{definition}

\begin{lemma}\label{lemma-O_Delta-in-C_Delta}
  The Hermitian $\cO$-module $\cO_{\vecDelta}$ is the $\cO$-span of a
  $\C$-basis of the vector space $\C_{\vecDelta}$, and thus it is an
  $\cO$-lattice in $\C_{\vecDelta}$.
\end{lemma}

\begin{proof}
  This follows from Lemma \ref{lemma-base-change} in the Appendix, or
  alternatively, by mimicking the proof of
  Lemma~\ref{lemma-O_Delta-in-R_Delta}.
\end{proof}

The following result is the careful statement of
Theorem~\ref{main:clatt-functor}.

\begin{theorem}[Complex Lorentzian lattice functor]\label{thm:clatt-functor}
  Let $r>0$ be a real algebraic integer and let $\zeta \neq \pm 1$ be
  a root of unity so that $\cO=\cO(r \zeta, r \zeta')$ is a PID.  Then
  there is a functor $\bigraph(r) \to \cO$-$\lattice$ sending
  $\vecDelta \mapsto \cO_\vecDelta$.  In other words, (a) if
  $\vecDelta$ is a connected directed bipartite graph with
  $\lambda_2=r$, then $\cO_\vecDelta$ is a Hermitian Lorentzian
  $\cO$-lattice, and (b) if $\vecDelta' \subset \vecDelta$ is a
  connected induced subgraph with the same second eigenvalue
  $\lambda_2=r$, then the inclusion map $\vecDelta' \into \vecDelta$
  induces an $\cO$-linear embedding map $\cO_{\vecDelta'} \into
  \cO_\vecDelta$ between Hermitian Lorentzian $\cO$-lattices sending
  $\spanning_{\vecDelta'} \into \spanning_{\vecDelta}$.
\end{theorem}
  
\begin{proof}
  The $\cO$-module $\cO_\vecDelta$ is an $\cO$-lattice in
  $\C_{\vecDelta}$ by Lemma~\ref{lemma-O_Delta-in-C_Delta} and as
  noted in \ref{def:c-spaces} $\C_{\vecDelta}$ has signature $(g,1)$.
  This proves part (a).  For part (b), the proof of essentially
  identical to the proof of Theorem~\ref{thm:latt-functor}), with the
  obvious changes: changing $\R$ to $\C$, $e_i$ to $\vec{e}_i$,
  $\alpha_i$ to $\vec{\alpha}_i$, and $A$ to $A_{\zeta}$.
\end{proof}

Corollaries~\ref{maincor:eisen} and~\ref{maincor:gauss} of
\S\ref{sec:intro} are immediate by taking $(r, \zeta) = (\sqrt{2},
\zeta_8)$ and $(\sqrt{3}, \zeta_4)$.

\begin{cor}\label{cor:gauss}
  There is a functor $\bigraph(\sqrt{2}) \to \gauss$-$\lattice$.
\end{cor}

\begin{cor}\label{cor:eisen}  
  There is a functor $\bigraph(\sqrt{3}) \to \eisen$-$\lattice$.
\end{cor}

\begin{table}
    \centering
    \begin{tabular}{|c|c|c|c|l|}
    \hline Point-Line & $\size{\vecDelta}$ & Spectrum & Sig. & Notes\\ \hline
    \hline
    $7_3$ & $14$ & $3^1 (\sqrt{2})^6 (-\sqrt{2})^6 (-3)^1$ & $(7,1)$ &
    $\vecDelta_{7}=$ Heawood graph \\ \hline
    $(4_3,6_2)$ & $10$ & $(\sqrt{6})^1 (\sqrt{2})^3 0^2 (-\sqrt{2})^3
    (-\sqrt{6})^1$ & $(6,1)$ & $\vecDelta_6= $ Subdivided $K_4$ \\ \hline
    $4_2$ & $8$ & $2^1 (\sqrt{2})^2 0^2 (-\sqrt{2})^2 (-2)^1$ & $(5,1)$ &
    $\vecDelta_5 =C_8$, an $8$-cycle \\ \hline
    \end{tabular}
    \vspace*{1mm}
    \caption{Three point-line incidence graphs with $\lambda_2 = \sqrt{2}$.}
    \label{tab:sqrt2}
\end{table}

\begin{table}
  \centering
  \begin{tabular}{|c|c|c|c|l|}
    \hline
    Point-Line & $\size{\vecDelta}$ & Spectrum & Sig. & Notes\\
    \hline
    \hline
    $13_4$  & $26$ & $4^1 (\sqrt{3})^{12} (-\sqrt{3})^{12} (-4)^1$ & $(13,1)$ &  $\vecDelta_{13}=$ bimonster graph \\ \hline
    $(9_4,12_3)$ & $21$ & $(\sqrt{12})^1 (\sqrt{3})^8 0^3 (-\sqrt{3})^8 (-\sqrt{12})^1$ & $(12,1)$ & $\vecDelta_{12} = \inc(\textrm{Hesse})$ \\ \hline
    $9_3$ & $18$ & $3^1 (\sqrt{3})^{6} 0^4 (-\sqrt{3})^{6} (-3)^1$  & $(11,1)$ & $\vecDelta_{11} = \inc(\textrm{Pappus})$ \\ \hline
    $8_3$ & $16$ & $3^1 (\sqrt{3})^{4} 1^3 (-1)^3 (-\sqrt{3})^{4} (-3)^1$  & $(11,1)$ & $\vecDelta_{11}' = \inc(\textrm{M\"obKan})$ \\ \hline
    $(9_2,6_3)$ & $15$ & $(\sqrt{6})^1 (\sqrt{3})^4 0^5 (-\sqrt{3})^4 (-\sqrt{6})^1$ & $(10,1)$ & $\vecDelta_{10} = \inc(\textrm{dual $K_{3,3}$})$ \\ \hline        
    $6_2$ & $12$ & $2^1 (\sqrt{3})^2 1^2 0^2 (-1)^2 (-\sqrt{3})^2 (-2)^1$ & $(9,1)$ & $\vecDelta_{9}= C_{12}$, a $12$-cycle \\ \hline
    \end{tabular}
    \vspace*{1mm}
    \caption{Six point-line incidence graphs with $\lambda_2 = \sqrt{3}$.}
    \label{tab:sqrt3}
\end{table}

Tables~\ref{tab:sqrt2} and~\ref{tab:sqrt3} list some bipartite graphs
with $\lambda_2=\sqrt{2}$ and $\lambda_2=\sqrt{3}$, respectively.  We
discuss these point-line graphs in \S\ref{sec:proj-planes} and some of
their lattices in \S\ref{sec:ball-quoients}.

\section{Standard Decomposition}\label{sec:std-decomposition}

Since the spectrum of a bipartite graph is symmetric, connected
bipartite graphs have two distinct $1$-dimensional eigenspaces,
leading to a different version of a standard decomposition
(Remark~\ref{rem:decomposition}).  Throughout this section $\vecDelta
= (\Delta, \sgn)$ is a connected directed bipartite graph with
$r=\lambda_2>0$, $\zeta \neq \pm 1$ is a root of unity, and
$\C_\vecDelta$ is the complex Lorentzian space defined in
Definition~\ref{def:c-spaces}.

\begin{defn}[Standard Decomposition]\label{def:c-decomposition}
  Similar to the real construction given in
  Remark~\ref{rem:decomposition}, the principal eigenvector of
  $A_\zeta$ spans a $1$-dimensional negative definite subspace in
  $\C^\vecDelta_r$, but now with eigenvalue $c(\lambda_2-\lambda_1)<0$
  for $B_\zeta$ because of rescaling.  This $\lambda_1$-eigenspace
  projects to the $1$-dimensional negative definite subspace in the
  quotient space $\C_{\vecDelta} =
  \C^\vecDelta_r/\rad(\C^\vecDelta_r)$.  In addition, there is a
  $1$-dimensional $(-\lambda_1)$-eigenspace.  If $x$ is the principal
  $\lambda_1$-eigenvector of $A_\zeta$ and we change the sign in the
  $\Delta^+$ coordinates of $x$, the resulting vector $y$ is a
  $(-\lambda_1)$-eigenvector of $A_\zeta$.  This eigenspace projects
  to the $1$-dimensional positive definite subspace in the quotient
  space $\C_{\vecDelta}$ with eigenvalue $c(\lambda_1 - \lambda_2)>0$
  for $B_{\zeta}$.  These two vectors span a signature $(1,1)$
  subspace $\C_{\vecDelta}^{1,1} \subset \C_{\vecDelta}$.  Let
  $\C_{\vecDelta}^+$ be its $(g-1)$-dimensional positive-definite
  orthogonal complement.  We call $\C_{\vecDelta} =
  \C_{\vecDelta}^{1,1} \oplus \C_{\vecDelta}^+$ the \emph{standard
  decomposition} of $\C_{\vecDelta}$.
\end{defn}

When $\vecDelta$ is $k$-regular, the standard decomposition has a
concrete description.  For each subset $S \subset \Delta$, let
\[
\vec{e}_S = \sum_{i\in S} \vec{e}_i
\text{\; and \;} 
\vec{\alpha}_S = \sum_{i\in S} \vec{\alpha}_i.
\]
Also, let $\one^- = \vec{e}_{\Delta^-}$ and $\one^+
=\vec{e}_{\Delta^+}$ with images $\alpha_- = \vec{\alpha}_{\Delta^-}$
and $\alpha_+ =\vec{\alpha}_{\Delta^+}$, dropping the vector symbols
since there is no danger of confusion with something defined in
$\R_\Delta$.
 
\begin{remark}[Standard Decomposition for $k$-Regular Graphs]\label{rem:c-reg-decomp}
  If $\vecDelta$ is $k$-regular, then $\one = \one^- + \one^+$ is a
  $k$-eigenvector and $\one^- - \one^+$ is a $(-k)$-eigenvector for
  the adjacency matrix $A$. Moreover, if $x$ is a
  $\lambda_i$-eigenvector of $A$, then $P^\htr_\zeta x$ is a
  $\lambda_i$-eigenvector of $A_\zeta = P^\htr_\zeta A P_\zeta$.
  Thus, $\one^- + (\zeta') \one^+$ is a $k$-eigenvector of $A_\zeta$
  with image $\alpha_\vecDelta :=\alpha_- + (\zeta')\alpha_+$ in
  $\C_\vecDelta$.  And $\one^- - (\zeta') \one^+$ is a
  $(-k)$-eigenvector of $A_\zeta$ with image $\alpha_- - (\zeta')
  \alpha_+$ in $\C_\vecDelta$.  In particular, $\C^{1,1}_\vecDelta$ is
  the $\C$-span of $\{\alpha_-,\alpha_+\}$.
  Finally, by Equation \eqref{eq:vec_e_i} $\one^- + (\zeta') \one^+$
  in the $\{\vec{e}_i\}$ basis and $\sqrt{c}\cdot (\one^- +
  \one^+)=\sqrt{c}\cdot \one$ in the $\{e_i\}$ basis are the same
  vector in $\C^\vecDelta_r$. As a consequence $\alpha_\vecDelta$ and
  $\sqrt{c}\cdot \alpha_\Delta$ are the same vector in $\C_\vecDelta$,
  so that $\alpha_\vecDelta$, $\alpha_\Delta$ and $\tau$
  (Definition~\ref{def:tau}) represent the same projective point.
\end{remark}

The remarks in Definition~\ref{def:aut-graph} hold for complex
Lorentzian representations using the complex standard decomposition
(Definition~\ref{def:c-decomposition}).
  
\begin{remark}[Graph Automorphisms]\label{rem:c-graph-auts}
  From the construction of $\cO^\vecDelta$ and $\cO_{\vecDelta}$ 
  (Definition~\ref{def:cplx-construction}) it is clear that
  $\aut(\vecDelta)$ acts on the singular $\cO$-module
  $\cO^{\vecDelta}$ by permuting the basis vectors $\{ \vec{e}_i
  \colon i \in \Delta\}$ and on the lattice $\cO_{\vecDelta}$ by
  permuting the spanning vectors $\{ \vec{\alpha}_i \colon i \in
  \Delta \}$ preserving the parts.  As a consequence, every automorphism 
  in $\aut(\vecDelta)$ fixes $\alpha_-$ and $\alpha_+$, 
  it fixes $\C^{1,1}_\vecDelta$ pointwise, and there is an 
  induced action on the positive definite space $\C_{\vecDelta}^+ \simeq \C^{g-1}$.  
  In short, we have a representation
  \begin{equation*}
    \aut(\vecDelta) \to \aut(\C_{\vecDelta}^+) \subseteq
    \aut(\C_{\vecDelta}) \simeq U(g,1). 
  \end{equation*}
\end{remark}

Automorphisms in $\aut(\Delta) \setminus \aut(\vecDelta)$ which swap
$\Delta^+$ and $\Delta^-$ (if they exist) permute the basis vectors in
$\cO^{\vecDelta}$ and spanning vectors in $\cO_\vecDelta$ in a way
that does not, in general, extend to an isometry of $\cO_{\vecDelta}$.
This is due to the complex nature of the off-diagonal entries of
$B_\zeta$.  There is, however, a projectivized automorphism of the
spanning set in the ball.

\begin{remark}[When does $\aut(\Delta)$ act on the ball?]\label{rem:aut-delta-action}
  Suppose $\vecDelta$ is a directed bipartite graph with a point-line
  reversing automorphism $\sigma$ of the underlying (undirected) graph
  $\Delta$.  For example, all of the exceptional bipartite graphs in
  \S\ref{sec:ball-quoients} have this property. There is an
  automorphism of $\C^\vecDelta$ that sends $\vec{e}_i$ to
  $\vec{e}_{\sigma(i)}$ for $i\in \Delta^-$ and $\vec{e}_i$ to
  $\zeta^2\cdot \vec{e}_{\sigma(i)}$ for $i\in \Delta^+$.  This
  descends to an automorphism of $\C_\vecDelta$ that sends
  $\vec{\alpha}_i$ to $\vec{\alpha}_{\sigma(i)}$ for $i\in \Delta^-$
  and $\vec{\alpha}_i$ to $\zeta^2\cdot \vec{\alpha}_{\sigma(i)}$ for
  $i\in \Delta^+$.  In the ball $\ball \C_\vecDelta \simeq \hyp^g_\C$,
  this becomes an automorphism of the projectivized spanning set whose
  order is that of $\sigma \in \aut(\Delta)$.  Note that any such
  $\sigma$ together with $\aut(\vecDelta)$ generate $\aut(\Delta)$.
  Thus $\aut(\Delta)$ acts on $\ball\C_\Delta$.
\end{remark}

Finally, we remark that when both the real and the complex
constructions can be applied, the lattices they produce are closely
connected.

\begin{remark}[Complexification and Real Form]
  Let $(\vecDelta, \sgn)$ be a connected directed bipartite graph with
  second largest eigenvalue $r > 0$ and let $\zeta$ be a suitable a
  root of unity.  If the rings of integers $\cO = \cO(r)$ and $\cO' =
  \cO(r \zeta, r \zeta')$ in the number fields $\Q(r)$ and
  $\Q(r\zeta,r\zeta')$ are both PIDs, then we can construct two
  distinct lattice representations of the graph: let $(\cO_\Delta,\{
  \alpha_i \})$ be the real Lorentzian representation in $\R_\Delta$
  and let $(\cO'_\vecDelta,\{ \vec{\alpha}_i\})$ be the complex
  Lorentzian representation in $\C_\vecDelta$.  Here $\cO_{\Delta}$ is
  an $\cO$-lattice and $\cO'_{\vecDelta}$ is an $\cO'$-lattice.  We
  think of $\cO'_{\vecDelta}$ as the ``complexification" of
  $\cO_{\Delta}$ and we think of $\cO_{\Delta}$ as the ``real form" of
  $\cO'_{\vecDelta}$.  The conversions go as follows.
  
  Starting with $\cO_\Delta$, tensoring with $\R$ produces
  $\R_\Delta$, and tensoring over $\C$ produces $\C_\vecDelta$
  (uniquely extending the symmetric bilinear form on $\R_\Delta$ to a
  Hermitian form that leaves the norms on $\R_\Delta$ unchanged).
  Once we stretch and twist the spanning set $\{\alpha_i\}$ to get the
  new spanning set $\{\vec{\alpha}_i\}$, the $\cO'$-span of
  $\{\vec{\alpha}_i\}$ is $\cO'_\vecDelta$.  This is essentially the
  construction given in this section.

  In the other direction, given $\cO'_\vecDelta$ with spanning set
  $\{\vec{\alpha}_i\}$, tensoring with $\C$ produces $\C_\vecDelta$
  with the same spanning set.  In this complex vector space, we can
  shrink and untwist $\{\vec{\alpha}_i\}$ to find a spanning set
  $\{\alpha_i\}$ with a real Gram matrix.  The $\R$-span of
  $\{\alpha_i\}$ is $\R_\Delta$ and the $\cO$-span of $\{\alpha_i\}$
  is $\cO_\Delta$.
\end{remark}

\section{Complex lattices from symmetric 2-designs}\label{sec:sym-2-designs}

We now present complex analogs of the results in \S\ref{sec:srg} and
\S\ref{sec:srg-lattices}.  Symmetric $2$-designs can be viewed as a
bipartite analog of strongly regular graphs, with a uniform structure
compatible with the bipartite-graph-to-Hermitian-lattice construction.
Concretely we show that the complex hyperbolic representations of
symmetric $2$-designs have properties similar to those of the real
hyperbolic representations of strongly regular graphs using nearly
identical proofs.

\begin{definition}[Symmetric $2$-designs]\label{def:sym-2-designs}
  A \emph{symmetric $2$-design} is a $k$-regular bipartite graph
  $\Delta$ where any two vertices distance~$2$ apart have exactly
  $\mu$ common neighbors.\footnote{A bipartite graph with these
  properties is, technically speaking, the incidence graph of a
  symmetric $2$-design.  Each can be constructed from the other and we
  choose to focus on the graph.}  Since every vertex is $k$-regular
  and every edge has a vertex of each type, $\abs{\Delta^+} =
  \abs{\Delta^-} = g$ for some $g$.  We say $\Delta$ is a symmetric
  $2$-design with \emph{parameters} $(g,k,\mu)$.  In the literature,
  the standard parameters are denoted by $(v,k,\lambda)$.  We deviate
  from this notation to highlight the similarities with the strong
  regular graphs discussed in \S\ref{sec:srg-lattices}.  The
  parameters $(g, k, \mu)$ of a symmetric $2$-design must satisfy the
  basic identity
  \begin{equation}
    \mu(g - 1) = k(k - 1).
    \label{eq-basic-identity-among-parameters-of-a-2-design}
  \end{equation}
  The Bruck-Chowla-Ryser theorem provides additional constraints on
  the parameters of a symmetric $2$-design \cite[Theorem
    8.3.1]{BrouwerVanMaldeghem22}.  A \emph{directed} symmetric
  $2$-design $\vecDelta = (\Delta,\sgn)$ is a symmetric $2$-design
  $\Delta$ with a sign function $\sgn$
  (Definition~\ref{def-directed-bipartite-graphs}).
\end{definition}

\begin{definition}[Identities]\label{def:c-identities}
  Let $A$ be the adjacency matrix of a symmetric $2$-design $\Delta$
  and let $N$ be the $g$-by-$g$ incidence matrix that records which
  elements of $\Delta^-$ are connected to which elements of $\Delta^+$
  (Definition~\ref{def-directed-bipartite-graphs}).  From the
  definition of a symmetric $2$-design, it is clear that
  \begin{equation}\label{eq:NNt}
    NN^\tr = N^\tr N = (k-\mu)I + \mu J,      
  \end{equation}
  a square matrix with all diagonal entries equal to $k$ and all
  off-diagonal entries equal to $\mu$.  The spectrum of a symmetric
  $2$-design with parameters $(g,k,\mu)$ is $\spec(\Delta) = \{k^1 \;
  r^f \; (-r)^f \; (-k)^1\}$ with $k>r>0$ and $r = \sqrt{k-\mu}$.
  There are a number of identities similar to those in
  Remark~\ref{rem:identities}.  Since $v = 2g = 2f+2$ we have $g = v -
  f - 1$, and if we write $s=-r$, we have
  \begin{equation}
    r + s =   0, \;\;\;\; r s = \mu - k.
    \label{eq-design-r-and-s}
  \end{equation}
  As a consequence, the second largest eigenvalue $r$ is either an
  rational integer or $\Q(r)$ is a real quadratic extension of $\Q$.
  One can also compute that
  \begin{equation}\label{eq:c-adj-matrix}
      (A_\zeta - rI)(A_\zeta - sI) =     A_\zeta^2 + (\mu-k)I = 
      A^2 + (\mu-k)I = \mu 
      \begin{bmatrix}
          J & 0 \\
          0 & J
      \end{bmatrix}
  \end{equation}
  using Equations~\ref{eq:A-square} and~\ref{eq:NNt}. Compare this
  equation to Equation~\ref{eq:adj-matrix-eq2} in
  Definition~\ref{def:bose-mesner}.
\end{definition}

For any directed symmetric $2$-design $\vecDelta$, there are special
elements $\vec{\tau}_i$ with uniform properties.  This is the complex
version of Definition~\ref{def:r-special}.

\begin{defn}[Special elements]\label{def:c-special}
  Let $\vecDelta = (\Delta, \sgn)$ be a $(g, k, \mu)$ directed
  symmetric $2$-design with spectrum $(k^1\; r^f \; s^f \; (-k)^1)$
  with $s=-r$, and let $(\cO_{\vecDelta}, \vec{\alpha}_i)$ be its
  complex Lorentzian representation defined over $\cO = \cO(r \zeta, r
  \zeta')$ for some fixed root of unity $\zeta$.  For each vertex $i
  \in\vecDelta$ we define a vector $\vec{\tau}_{i} \in \C^\Delta$
  whose coordinates with respect to the ordered basis $\{\vec{e}_i\}$
  is the $i$-th column of the matrix
  \[
  A_\zeta-sI = 
  \begin{bmatrix}
    -sI & \zeta N \\
    \zeta' N^\tr & -sI
  \end{bmatrix}.
  \]
  Concretely, if $N(i)$ denotes the set of neighbors of $i \in \vecDelta$, 
  \begin{equation}\label{eq:vec-tau}
    \vec{\tau}_{i} = 
    \begin{cases}
      - s \vec{e}_i + \zeta'\vec{e}_{N(i)} & \textrm{ if $i \in \Delta^-$} \\
      \zeta\vec{e}_{N(i)} - s \vec{e}_i  & \textrm{ if $i \in \Delta^+$} 
    \end{cases}
  \end{equation}
\end{defn}

The special elements $\vec{\tau}_i$ satisfy a complex version of
Lemma~\ref{lem:elts-in-rad}.

\begin{lemma}[Radical Elements]\label{lem:c-elts-in-rad}
  Let $\vecDelta$ be a directed symmetric $2$-design with parameters
  $(g,k,\mu)$.  The elements $\vec{e}_i$ and $\vec{\tau}_j$ have inner
  product
  \[
  \ip{\vec{e}_i}{\vec{\tau}_j} = 
  \begin{cases}
      -c\mu & \textrm{if $\sgn(i)=\sgn(j)$},\\
      0 &  \textrm{if $\sgn(i)\neq\sgn(j)$}.\\
  \end{cases} 
  \]
  Thus $\vec{\tau}_j - \vec{\tau}_{j'} \in \rad(\C^{\vecDelta}_r)$
  for all $j$ and $j'$ with $\sgn(j)=\sgn(j')$.  In fact, these vectors 
  span the radical of $\C^\vecDelta_r$.
\end{lemma}

\begin{proof}
  The elements $\vec{e}_i$ are column vectors of the identity matrix $I$, 
  and the elements $\vec{\tau}_j$ are column vectors of $A_\zeta-sI$.  
  Using Equation~\ref{eq:c-adj-matrix} in Definition~\ref{def:sym-2-designs}
  we have
  \[ 
  \llbracket \ip{\vec{e}_i}{\vec{\tau}_j} \rrbracket = I^\tr\cdot c (rI-A_\zeta)\cdot (A_\zeta-sI) = -c \mu 
  \begin{bmatrix}
      J & 0 \\
      0 & J 
  \end{bmatrix}
  \]
  which shows that for all $i,j\in \Delta$ we have
  $\ip{\vec{e}_i}{\vec{\tau}_j}=-c\mu$ when $i$ and $j$ have the same
  sign and $0$ when they have the opposite signs.
  
  When $\sgn(j) = \sgn(j')$, one has $\vec{\tau}_j - \vec{\tau}_{j'}
  \in \rad(\C^{\vecDelta}_r)$ because $\vec{\tau}_j - \vec{\tau}_{j'}$
  is orthogonal to all of the basis vectors $\vec{e}_i$.

  The proof that these vectors span the radical is similar to the
  proof in Lemma~\ref{lem:elts-in-rad}. We write $Z = \C\lbrace
  \vec{\tau}_i - \vec{\tau}_j \colon i, j \in \Delta^+ \text{\; or \;}
  i, j \in \Delta^- \rbrace$.  Let $u^{\pm} = \mathbf{1}^- \pm \zeta'
  \mathbf{1}^+$. So $\C \lbrace u^+, u^- \rbrace = \C \lbrace
  \mathbf{1}^+, \mathbf{1}^- \rbrace$.  Note that $\C^v = \C u^+ \bot
  \C^g_+ \bot \C^g_- \bot \C u^-$ where
  \[
  \C u^{\pm} = \op{ker}(A_{\zeta} \mp k I)
  \text{\;and \; } \C^g_{\pm} = \op{ker}(A_{\zeta} \mp r I)
  \]
  and where $\bot$ means orthogonal with respect to the standard
  Hermitian dot product on $\C^v$.  Note that $(A_{\zeta} + r I)$
  kills $\C^g_-$ and acts on the other three subspaces as nonzero
  scalars, so
  \[
  \C \lbrace \vec{\tau}_i \colon i \in \Delta \rbrace
  =\op{im}(A_{\zeta} + r I) = 
  \C u^+ \bot \C^g_+  \bot \C u^-.
  \]
  Since $Z$ together with a single $\vec{\tau}_i$ from $\Delta^{\pm}$
  each span this space $\op{im}(A_{\zeta} + r I)$, we have
  $\op{dim}(Z) \geq g$.  But if $i, j \in \Delta^+$ or $i, j \in
  \Delta^-$, then $\vec{\tau}_i - \vec{\tau}_j$ is orthogonal to
  $\mathbf{1}^+$ and $\mathbf{1}^-$ with respect to the standard
  Hermitian dot product on $\C^v$, so $Z \subseteq \C^g_+$. It follows
  that $ Z = \C^g_+ = \op{rad}(\C^\Delta_r)$.
\end{proof}

As in Definition~\ref{def:tau}, the special vectors $\vec{\tau}_i$
coalesce in the quotient, but here there are two images, one for each
part of $\vecDelta$.  The subscripts on the images, which may appear
to have been switched, have been chosen to make sure that $\proj
\tau_{\varepsilon}$ is the point of intersection of the mirrors
$\lbrace \vec{\alpha}_i^{\bot} \colon i \in \Delta^{\varepsilon}
\rbrace$ in complex hyperbolic space.

\begin{figure}
  \includegraphics{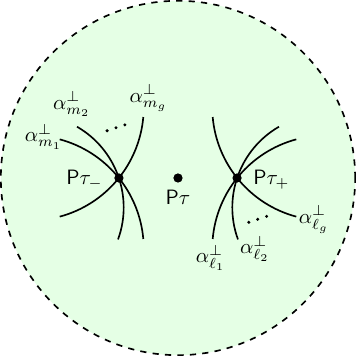}
  \caption{A schematic drawing of the mirrors corresponding to points
    $\Delta^- = \{ m_1, \dotsb, m_g \}$ passing through $\proj\tau_-$
    and the mirror corresponding to the lines $\Delta^+ = \{ {\ell_1},
    \dotsb, \ell_g \}$ passing through $\proj\tau_+$, and the midpoint
    $\proj\tau$ between them.}\label{fig:mirrors}
\end{figure}

\begin{defn}[Common Images $\tau_{\pm}$]\label{def:taupm}
  Let $\tau_+ \in \C_{\vecDelta}$ be the common image of the vectors
  $\{ \vec{\tau}_i \colon i \in \Delta^- \} \subset \C^\vecDelta$ and
  let $\tau_- \in \C_{\vecDelta}$ be the common image of the vectors
  $\{ \vec{\tau}_j \colon j \in \Delta^+ \} \subset \C^\vecDelta$.
  The vector symbol has been dropped since there is no danger of
  confusion with something defined in $\R_\Delta$.  By summing the $g$
  equations for $\tau_+$ and the $g$ equations for $\tau_-$, given in
  Equation~\ref{eq:vec-tau} in Definition~\ref{def:c-special} we see
  that
  \begin{equation}\label{eq:alpha-tau}
    \begin{array}{lcrcr}
        g \cdot \tau_+ &=& (-s)\cdot \alpha_- &+& (k\;\zeta')\cdot \alpha_+ \\
        g \cdot \tau_- &=&  (k\;\zeta) \cdot \alpha_- &+& (-s)\cdot \alpha_+
    \end{array}
  \end{equation}
  Both $\tau_+$ and $\tau_-$ are in $\C^{1,1}_\vecDelta$, the subspace
  spanned by $\{\alpha_-,\alpha_+\}$ and pointwise fixed by
  $\aut(\vecDelta)$ (Remark~\ref{rem:c-reg-decomp}). The equations in
  Lemma~\ref{lem:c-elts-in-rad} show that for $\varepsilon \in \{ \pm
  1 \}$,
  \begin{equation}\label{eq:ip-alpha-tau}
  \ip{\vec{\alpha}_i}{\tau_{\varepsilon}}
  = \begin{cases}
  0 & \text{if \;} \varepsilon = \sgn(i), \\
  -c\; \mu
  & \text{if \;} \varepsilon \neq \sgn(i).
  \end{cases}
  \end{equation}
  In other words, $\tau_- \in \bigcap \{\vec{\alpha}_i^\bot \colon i
  \in \Delta^- \}$ and $\tau_+ \in \bigcap \{\vec{\alpha}_i^\bot
  \colon i \in \Delta^+ \}$.  The subscript on $\tau_\varepsilon$
  indicates it is in the intersection of the $\Delta^\varepsilon$
  mirrors.  Note that $\proj\tau_{\varepsilon}$ is the \emph{unique}
  projective point in the intersection of the $\Delta^\varepsilon$
  mirrors because the $g$ spanning vectors defining these mirrors are
  pairwise orthogonal.  See Figure~\ref{fig:mirrors} for a schematic
  illustration.
\end{defn}

\begin{defn}[Midpoint $\tau_\vecDelta$]\label{def:midpoint}
  If we define $\tau_\vecDelta := \tau_+ + (\zeta')\tau_-$ on analogy
  with $\alpha_\vecDelta := \alpha_- + (\zeta')\alpha_+$
  (Remark~\ref{rem:c-reg-decomp}), keeping in mind the subscript
  switch, then from Equation~\ref{eq:alpha-tau} we get the symmetric
  $2$-design analog of Equation~\ref{eq:tau-alpha} in
  Definition~\ref{def:tau}:
  \begin{equation}
      g\cdot \tau_\vecDelta = (k-s)\cdot \alpha_\vecDelta.
      \label{eq:tau-and-alpha-vecDelta}
  \end{equation}
  Taking inner product of both sides of
  Equation~\eqref{eq:tau-and-alpha-vecDelta} with $\vec{\alpha}_i$
  (and dividing by $-c$) for some $i \in \Delta^-$ yields $g \mu = (k
  - r)(k - s)$ which is equivalent to the basic identity
  \eqref{eq-basic-identity-among-parameters-of-a-2-design} among the
  parameters of a symmetric $2$-design.  This is exactly analogous to
  the real case.  See Definition~\ref{def:tau}.
    
  Projectively, $\proj\tau_\vecDelta = \proj\alpha_\vecDelta =
  \proj\alpha_\Delta = \proj\tau_\Delta = \proj\tau$.  From
  Equations~\ref{eq:alpha-tau} and~\ref{eq:ip-alpha-tau} one computes
  that the common norm of $\tau_-$ and $\tau_+$ is the negative
  integer $(-s)(-c\;\mu) = -n\;\mu \in \Z_{<0}$.  Similarly, one
  computes that $\ip{\tau_+}{\zeta'\tau_-} = - c\; \mu\; k \in
  \R_{<0}$ is a negative real number.  As a consequence,
  $\tau_\vecDelta = \tau_+ + (\zeta')\tau_-$ represents the midpoint
  $\proj\tau_\vecDelta = \proj\tau$ of the real geodesic joining
  $\proj\tau_+$ and $\proj\tau_-$.
  
  If $\Delta$ has an automorphism $\sigma$ that interchanges
  $\Delta^+$ and $\Delta^-$, then, as discussed in
  Remark~\ref{rem:aut-delta-action}, the group $\op{Aut}(\Delta) =
  \langle \op{Aut}(\vecDelta), \sigma \rangle$ acts on the ball $\ball
  \C_{\vecDelta} \simeq \ball^g_{\C}$.  The action of
  $\op{Aut}(\vecDelta)$ pointwise fixes the $\ball^1_{\C}$ spanned by
  $\proj \tau_-$ and $\proj \tau_+$.  On the other hand, $\sigma$
  interchanges $\proj \tau_+$ and $\proj \tau_-$.  So the midpoint
  $\proj \tau$ is fixed by $\sigma$, hence by $\op{Aut}(\Delta)$.  In
  many examples $\proj \tau$ is the unique point in the ball
  $\ball^g_\C$ fixed by $\op{Aut}(\Delta)$.
 \end{defn}

The space $\C^{1,1}_\vecDelta$ is an analog of the Coxeter axis of a
bipartite Coxeter element.

\begin{remark}[$\C^{1,1}_\vecDelta$ and the Coxeter axis]
  Let $\vecDelta$ be a directed bipartite graph and let
  $W=\cox(\vecDelta,2)$ be the corresponding simply-laced Coxeter
  group.  It has a natural linear action on the vector space
  $\R_{\Delta,2}$ (Definition~\ref{def:shift-module}).  We write $w_-$
  for the product of the reflections $\{R_{\alpha_i} \colon i \in
  \Delta^-\}$ and $w_+$ for the product of the reflections
  $\{R_{\alpha_i} \colon i \in \Delta^+\}$.  Since both sets of
  reflections pairwise commute, the products are well-defined.  The
  product $w_- \cdot w_+$ is called the \emph{bipartite Coxeter
  element} of $W$. When the graph $\vecDelta$ is neither a Dynkin
  diagram nor an extended Dynkin diagram, the group $W$ is neither
  spherical nor affine and the standard faithful action of $W$ on a
  real vector space preserves a proper convex cone called the
  \emph{Tits cone}, and the bipartite Coxeter element $w=w_- \cdot
  w_+$ acts on the Tits cone with north-south dynamics.  For a
  description of the dynamical properties of bipartite Coxeter
  elements, see \cite{McMullen2002}.  The action of $w$ has a
  repelling fixed null line $\R z_-$, an attracting fixed null line
  $\R z_+$, and the subspace $\R^{1,1}_\vecDelta = \R\{z_-,z_+\}$ has
  signature $(1,1)$.  This is the \emph{Coxeter plane}, or the
  \emph{Coxeter axis} when working projectively.  The involutions
  $w_-$ and $w_+$ act on the Coxeter axis and they fix projective
  points that one could call $\proj\tau_-$ and $\proj\tau_+$.
\end{remark}

\section{Examples: Finite Projective Planes}\label{sec:proj-planes}

The incidence graph of the points and lines of a finite projective
plane is a classic and motivating example of a symmetric $2$-design.
This section contains a general discussion of lattices one can build
from these graphs, and from point-line configurations inside them that
are $2$-designs.  The details of certain specific examples with
exceptional properties are postponed until \S\ref{sec:ball-quoients}.

\begin{definition}[Finite Projective Planes]\label{def:p2fq}
  Let $\vecDelta = \inc(\proj^2_\F)$ be the incidence graph of the
  finite projective plane over finite field $\F = \F_q$ of size $q$,
  where $q \geq 2$ is a prime power.  To keep the notation simple
  while maintaining distinct names, we write $\vecDelta = \vecDelta_g$
  for these graphs where the integer $g=q^2+q+1$.  The graph
  $\vecDelta_g$ is a $(g,k,\mu)$-symmetric $2$-design with $k=q+1$,
  $\mu=1$, $r = \sqrt{k-\mu} = \sqrt{q}$ and spectrum
  \[
  \spec(\vecDelta_g) = \{(q+1)^1 \;
  (\sqrt{q})^{q^2+q} \;
  (-\sqrt{q})^{q^2+q}\;
  (-q-1)^1
  \}.
  \]
  It has points $\points=\Delta^- = \{m_1,\ldots,m_g\}$ and lines
  $\lines = \Delta^+ = \{\ell_1,\ldots,\ell_g\}$.  Note that the
  signature of $\C_\vecDelta$ is $\C^{g,1}$, so the use of $g$ as an
  index is consistent with Notation~\ref{not:dim-index} used in the
  final section.
\end{definition}

We give the incidence graphs of small projective planes special names. 

\begin{remark}[Heawood and the bimonster]\label{rem:heawood-bimonster}
  When $q=2$, the $14$-vertex graph $\vecDelta_7$ is generally known
  as the \emph{Heawood graph}.  When $q=3$, the $26$-vertex graph
  $\vecDelta_{13}$ does not have a standard name.  We call it the
  \emph{bimonster graph} since it plays a key role in Allcock's
  `monstrous proposal' and the Conway-Ivanov-Norton-Simons
  presentation of the bimonster group, the monster finite simple group
  wreath $\Z/2$; see \cite{ConwaySimons2001, Allcock09-monstrous,
    AB18-generators, AB25-deflation}.
\end{remark}

The $6$-cycle is a symmetric $2$-design closely connected to projective planes.
 
\begin{remark}[The hexagon and $\F_1$]\label{rem:f1}
  The graph $\vecDelta = \inc(\proj^2_{\F})$ is an example of a
  spherical building of type $A_2$ and the many hexagons it contains
  are its apartments.  When Jacques Tits first introduced buildings in
  1957, he suggested that these apartments should be viewed as the
  corresponding building over the ``field with one element''
  (\cite{Tits57}). In other words, the hexagon $C_6$, a
  $(3,2,1)$-symmetric $2$-design, should be viewed as
  `$\inc(\proj^2_{\F_1})$'.
\end{remark}

Bi-regular subsets of $\vecDelta =\inc(\proj^2_\F)$ are point-line
configurations. These are $2$-designs that need not be symmetric.

\begin{definition}[Point-Line Configurations]
  Let $\vecDelta' \subset \vecDelta$ be a connected induced bipartite
  subgraph of a finite projective plane incidence graph with parts
  $\points' \subset \points$ and $\lines' \subset \lines$.  If
  $\vecDelta'$ has $a=\size{\points'}$ points, $b=\size{\lines'}$
  lines, and it is bi-regular, with every point vertex in $\points'$
  connected to $c$ line vertices and every line vertex in $\lines'$
  connected to $d$ point vertices, then $\vecDelta'$ is called an
  \emph{$(a_c,b_d)$-configuration}. When the configuration is
  \emph{self-dual} with $a=b$ and $c=d$, $\vecDelta'$ is a symmetric
  $2$-design, and the notation is simplified to $a_c$.
\end{definition}

\begin{example}[Small Examples]
  All of the bipartite graphs listed in Tables~\ref{tab:sqrt2}
  and~\ref{tab:sqrt3} are point-line configurations inside the Heawood
  graph (Table~\ref{tab:sqrt2}) and the bimonster graph
  (Table~\ref{tab:sqrt3}).  The Heawood graph is a $7_3$-configuration
  and the bimonster graph is a $13_4$-configuration.  The others are
  famous configurations associated with the names Hesse $(9_4,12_6)$,
  Pappus $9_3$, and M\"obius-Kantor $8_3$.  The rest are elementary
  geometric configurations.
\end{example}

\begin{example}[Lattices from $\vecDelta_g$]\label{ex:proj-plane-lattices}
  Let $\vecDelta_{g}$ be the incidence graph of a finite projective
  plane as defined in Definition~\ref{def:p2fq}.  When $q=2$ one can
  apply Corollary~\ref{cor:gauss} to produce a Gaussian lattice and
  when $q=3$ one can apply Corollary~\ref{cor:eisen} to produce an
  Eisenstein lattice.  More generally, When $q$ is a prime, $q \equiv
  3 \bmod 4$, and we choose $\zeta = i$, the ring of integers of
  $\Q(r\zeta,r\zeta')$ is $\cO = \Z[(1+\sqrt{-q})/2]$.  Let $\cO_g =
  \cO_{\vecDelta_g}$ be the complex Lorentzian representation of the
  graph $\vecDelta_{g}$ defined over $\cO$.  These $\cO$-modules were
  studied by the first author in \cite{Basak14-proj-planes}.  For
  values of $q$ where $\cO$ is not a PID, our results as stated do not
  literally apply.  Nevertheless, it was shown in
  \cite{Basak14-proj-planes} that the $\cO$-modules $\cO_g$ are
  $\cO$-lattices, in fact free $\cO$-modules, not merely projective.
  This gives a family of $\cO$-lattices of signature $(g, 1)$ that are
  $\sqrt{-q}$-modular, i.e. $\sqrt{-q}\cdot \cO_{g}^\vee = \cO_{g}$,
  where $\cO_g^\vee$ is the dual lattice
  (Definition~\ref{def:self-dual}).  Write $\F = \F_q$.  By
  Remark~\ref{rem:c-graph-auts}, the group $\aut(\vecDelta) =
  \op{PGL}_3(\F)$ acts on the lattice $\cO_g$ transitively permuting
  the $g$ orthogonal mirrors $\{ \vec{\alpha}_l^{\bot} \colon l \in
  \Delta^+ \}$ indexed by the lines of $\proj^2_\F$, transitively
  permuting the $g$ orthogonal mirrors $\{ \vec{\alpha}_m^{\bot}
  \colon m \in \Delta^-\}$ indexed by the points of $\proj^2_\F$, and
  fixing the subspace $\C^{1,1}_{\vecDelta_g}$ spanned by $\tau_-$ and
  $\tau_+$ as schematically shown in Figure~\ref{fig:mirrors}.
\end{example}

\begin{remark}[Leech-like Lattices]
  If $q$ is a $3 \bmod 4$ prime and there is some element in
  $\Q(\sqrt{-q})$ of norm $g$,\footnote{Equivalently, the ternary
  quadratic form $x^2 + q y^2 - g z^2$ represents $0$ over $\Z$.
  Number theory conjectures, such as the Bateman-Horn conjecture,
  imply there are infinitely many such primes.} then one can find a
  primitive norm zero vector $\rho$ in $\cO_{g}$ that is fixed by
  $\op{PGL}_3(\F)$.  This is true, for example, when $q = 3$, $47$,
  $59$, $71, \dotsb$.
  It follows that $\cO_{g}$ contains a copy of the positive definite
  lattice $\rho^{\bot}/\rho$ and, by taking its rescaled real form,
  one finds a positive definite even self-dual $\Z$-lattice $M_{q}$ of
  rank $2 g - 2 = 2q(q + 1)$ that has $\op{PGL}_3(\F)$ in its symmetry
  group. When $q=3$ the $\Z$-lattice $M_3$ is the Leech
  lattice~$\Lambda$. See \cite{Basak14-proj-planes} \S3.6 and \S3.7
  for details.
\end{remark}

\section{Examples: Complex lattices and ball quotients}\label{sec:ball-quoients}

This final section discusses several pairs of complex Lorentzian
lattices (and one pair of quaternionic Lorentzian lattices) that have
already appeared in the literature.  Their similarities prompted this
study and they can now be constructed in a more systematic fashion.
For a summary of these similarities, see
Remark~\ref{rem:similarities-summary}. Let $\clatt$ be one of the
lattices under consideration.  We begin by establishing common
notation for $\clatt$, its projectivized reflection group $\Gamma$
with mirrors $\eu{H}$ and the corresponding ball quotient $X$.

\begin{table}
  \centering
  \begin{tabular}{|c|c|c|c|c|}
    \hline
    Example &Graph & $\lambda_2$ & Lattice & Moduli Space\\
    \hline \hline
    \ref{ex:eisen-petersen} & Petersen & $1$ & $\eisen_4$  & cubic surfaces in $\proj_\C^3$\\
    \ref{ex:6-gon} & $6$-cycle & $1$ & $\eisen_3$ & $6$ points in $\proj_\C^1$ \\ \hline
    \ref{ex:gauss-heawood} & Heawood & $\sqrt{2}$ & $\gauss_{7}$ & \\
    \ref{ex:8-gon} & $8$-cycle & $\sqrt{2}$ & $\gauss_{5}$ & $8$ points in $\proj_\C^1$ \\ \hline
    \ref{ex:monstrous} & bimonster & $\sqrt{3}$ & $\eisen_{13}$ &  \\
    \ref{ex:12-gon} & $12$-cycle & $\sqrt{3}$ & $\eisen_{9}$ & $12$ points in $\proj_\C^1$ \\ \hline
  \end{tabular}
  \vspace*{3mm}
  \caption{Three pairs of exceptional lattices and sublattices.}
  \label{tab:lattice-pairs}
\end{table}

\begin{remark}[Lattice Pairs]\label{rem:example-lattices}
  Three pairs of exceptional complex Lorentzian lattices are listed in
  Table~\ref{tab:lattice-pairs}.  The larger lattice is constructed
  from a graph and the smaller lattice is a sub-lattice constructed
  from the largest induced cycle contained in the graph.  Note that
  the graph and the cyclic subgraph always share a common second
  largest eigenvalue $\lambda_2$.
  In each case, there is a ring $\cO$ (which is the Eisenstein
  integers $\eisen$ or the Gaussian integers $\gauss$) so that the
  complex Lorentzian lattice $\clatt=\cO_\vecDelta$ is the complex
  Lorentzian representation of a directed bipartite graph $\vecDelta$
  (Definition~\ref{def:cplx-construction}), or $\clatt =\cO_\Delta$ is
  a `complexified' version of a real Lorentzian representation of a
  graph $\Delta$ (Definition~\ref{def:shift-module}) with the bilinear
  form naturally promoted to a Hermitian form.
\end{remark}

As noted in Table~\ref{tab:lattice-pairs}, some of the ball quotients
$X$ have moduli space interpretations, namely, the space of cubic
surfaces in $\proj^3_\C$ and the spaces $\mc{M}_{0,n}$ of $n$ points
on $\proj^1_\C$ for $n = 6, 8, 12$.  One of these ball quotients is
conjecturally related to the Monster simple group.

\begin{definition}[Roots and Reflections]\label{def:example-refl}
 The Hermitian $\cO$-lattice $\clatt$ under consideration comes with a
 distinguished spanning set $\{\vec{\alpha}_i \colon i \in \Delta\}$.
 Since $\lambda_2$ is $1$, $\sqrt{2}$, or $\sqrt{3}$ and $c = r =
 \lambda_2$ in the bipartite cases, the vectors $\vec{\alpha}_i$ have
 integer norm $r^2=\ip{r}{r}=1$, $2$, or $3$.
 We call the vectors in $\clatt(r^2)$ the \textit{roots of $\clatt$}
 because $\op{Aut}(\clatt)$ contains a $\zeta_d$-complex reflection
 for each vector in $\clatt(r^2)$.  The largest integer $d>0$ with
 this property only depends on $r^2$; when $r^2=1,2,3$, one has
 $d=6,4,3$, respectively.  Let $\Refl(\clatt) \subseteq U(g,1)$ be the
 group generated by the order~$d$ complex reflections in the roots,
 and let $\Gamma = \proj \Refl (\clatt) \subseteq \proj U(g, 1)$ be
 the projectivized reflection group that acts faithfully on the
 complex hyperbolic space $\ball \clatt \simeq \hyp^g_\C$, where
 $\ball \clatt$ is a shorthand for $\ball (\C \otimes_{\cO} \clatt)$.
\end{definition}

\begin{definition}[Mirrors and the Ball Quotient]\label{def:example-ball}
  Let $\eu{H}$ be the union of the mirrors of $\Gamma$, i.e. the
  hyperplanes in $\ball \clatt$ determined by the orthogonal
  complements of roots of $\clatt$, and let $\ball \clatt^{\circ} =
  \ball \clatt - \eu{H}$.  We write
  \[
  X = \Gamma \backslash \ball \clatt \text{\; \and \;} X^{\circ} 
  = \Gamma \backslash \ball \clatt^{\circ} = X - (\Gamma \backslash \eu{H})
  \]
  for the corresponding ball quotients.
\end{definition}

\begin{notation}[Indexed by Dimension]\label{not:dim-index}
  To keep the notation simple, we write $\vecDelta = \vecDelta_g$ or
  $\Delta = \Delta_g$ for the graphs, as in Definition~\ref{def:p2fq},
  and
  \[
  (\clatt, \Gamma, \eu{H}, X) = (\cO_g, \Gamma_g, \eu{H}_g, X_g)
  \]
  where the integer $g$ is the complex dimension of the ball
  $\ball\clatt \simeq \hyp^g_\C$.  The ball quotient $X_g$ has complex
  dimension $g$ and the lattice $\cO_g$ has signature $(g,1)$. If
  $\vecDelta$ is a symmetric $2$-design then this $g$ is the same as
  the $g$ used earlier.
\end{notation}

The first pair of examples is a complexification of the real
construction.

\begin{example}[Eisenstein lattice from the Petersen graph]\label{ex:eisen-petersen}
  Let $\petersen = \Delta = \kneser(5,2)$ be the Petersen graph.  The
  complex lattice we want to consider here is
  \[
  \eisen_4  
  = \eisen \otimes_{\Z} \Z_{\petersen} \simeq \eisen^{4, 1}
  \]
  with an $\eisen$-valued Hermitian form obtained by naturally
  extending the $\Z$-valued symmetric bilinear form of
  $\Z_{\eu{P}}$.\footnote{
  The $\eisen$-modules $\eisen_{4}$ considered here and the
  $\eisen$-module $\eisen_{\eu{P}}$ in the notation of
  \S\ref{sec:graph-lattice} are the same $\eisen$-module, but with
  different forms.  One is Hermitian and the other is symmetric
  bilinear.  See the footnote to Definition~\ref{def:c-spaces}.}
  Note that $\Z_{\eu{P}}$ is the smallest among the lattices discussed
  in \S\ref{sec:kneser}.  Let $\Gamma_4$ be the group generated by the
  order six complex reflections in the `short roots' $\eisen_4(1)$ of
  $\eisen_4$ and let $\eu{H}_4$ be the mirrors in $\ball \eisen_4
  \simeq \hyp^4_\C$ orthogonal to them.  The ball quotient $X_4 =
  \Gamma_4 \backslash \ball \eisen_4$ is the space of cubic surfaces
  in $\proj^3_\C$. One can find an explicit correspondence between a
  neighborhood of $\tau$ in $\mathbb{B} \eisen_4$ (with the
  configuration of mirrors $\{ \alpha_i^{\bot} \colon i \in \petersen
  \}$ arranged around $\tau$) and a neighborhood of the Clebsch
  surface in the Sylvester family of cubic surfaces such that the
  point $\tau$ corresponds to the Clebsch surface and the ten mirrors
  closest to $\tau$ correspond to ten simplest deformations of the
  Clebsch surface to a surface with a single nodal singularity.  This
  leads to a presentation of the orbifold fundamental group of $X_4$
  as a quotient of Artin group of $\petersen$ \cite{ABL} and an
  explicit description of the monodromy action of the ten generators
  of $\pi_1^{\op{orb}}(X_4, \proj \tau)$ on the $27$ lines of the
  Clebsch surface in terms of $\petersen$~\cite{B-monodromy}.
\end{example}

The largest induced cycle in the Petersen graph $\petersen$ is a $6$-cycle.

\begin{example}[Eisenstein lattice from a $6$-cycle]\label{ex:6-gon}
  The Petersen graph contains an induced cycle $\Delta'=C_6$ with the
  same second eigenvalue $\lambda_2=1$, so by
  Theorem~\ref{thm:latt-functor}(b) the $\Z$-lattice $\Z_{\Delta'}$
  embeds as a sublattice of $\Z_{\eu{P}}$.  Similarly, the
  $\eisen$-lattice
  \[
  \eisen_3 = \eisen \otimes_{\Z} \Z_{\Delta'}
  \]
  with its $\eisen$-valued Hermitian form obtained by extending the
  $\Z$-valued symmetric bilinear form on $\Z_{\Delta'}$, embeds as a
  sublattice of $\eisen_4$ (Example~\ref{ex:eisen-petersen}).  To
  describe $\eisen_3$ as a concrete sublattice of $\eisen_4$, fix a
  vertex $i \in \petersen$.  Then $\eisen_3$ is the orthogonal
  complement of the short root $\vec{\alpha}_i$ in $\eisen_4$ and is
  spanned by the six roots $\vec{\alpha}_j$, as $j$ varies over the
  hexagon in $\petersen$ that is disconnected from the vertex $i$.
  Let $\Gamma_3$ be the group generated by the order six complex
  reflections in the short roots of $\eisen_3$ and let $X_3 \simeq
  \Gamma_3 \backslash \ball \eisen_3$ be the corresponding ball
  quotient.  Then $X_3$ is the space of six points in $\proj^1_{\C}$.
  This example is used in \cite{ABL} in the same way the moduli space
  of twelve points in $\proj^1_{\C}$ (Example~\ref{ex:12-gon}) is used
  in the monstrous proposal \cite{AB25-deflation}.
\end{example}

The next pair of examples are from graphs with $\lambda_2=\sqrt{2}$.

\begin{remark}[Bipartite graphs with $\lambda_2=\sqrt{2}$]
  The three point-line configurations in Table~\ref{tab:sqrt2} have
  incidence graphs with $\lambda_2=\sqrt{2}$, and by
  Corollary~\ref{cor:eisen} these produce Lorentzian
  $\gauss$-lattices.  There are nested inclusions among these
  configurations, graphs, and lattices.
\end{remark}

The second exceptional pair in Table~\ref{tab:sqrt2} comes from the
largest and smallest of these graphs.  The largest lattice is a
special case of Example~\ref{ex:proj-plane-lattices}.

\begin{example}[Gaussian lattice from the Heawood graph]\label{ex:gauss-heawood}
  Let $\vecDelta_7 = \inc(\proj^2_\F)$ be the Heawood graph with
  $\F=\F_2$ and $\lambda_2 = \sqrt{2}$. By Corollary~\ref{cor:gauss}
  there is a Lorentzian Gaussian representation
  $(\gauss_{7},\{\vec{\alpha}_i\})$ of $\vecDelta_7$ in $\C^{7,1}$.
  One can show that
  \[
  \gauss_7 \simeq 3 D_4^{\gauss} \oplus \II_{1,1}^{\gauss, 1 + i}
  \]
  where $D_4^{\gauss}$ is a ``complex form" of the $D_4$ lattice
  defined over $\gauss$ (see \cite{Basak20-barnes-wall} \S 2.2 for a
  definition of $D_4^{\gauss}$).  In other words, $\gauss_7$ is
  isomorphic to three copies of the Gaussian $D_4$ lattice plus a
  hyperbolic cell $\II_{1,1}^{\gauss, 1 + i}$ (Definition
  \ref{def:hyp-cell}). This example and some related ones were studied
  in the dissertation of Jeremiah Goertz (Iowa State PhD thesis, 2016,
  unpublished) under the supervision of the first author.  Some of
  these results appear in \cite{Basak20-barnes-wall}.
\end{example}

The largest induced cycle in the Heawood graph is an $8$-cycle.

\begin{example}[Gaussian lattice from an $8$-cycle]\label{ex:8-gon}
  The Heawood graph contains an induced cycle $\vecDelta'=C_8$ with
  the same second eigenvalue $\lambda_2=\sqrt{2}$, so by
  Corollary~\ref{cor:gauss}, the $\gauss$-lattice $\gauss_5 =
  \gauss_{\vecDelta'}$ embeds as a sublattice of $\gauss_7$
  (Example~\ref{ex:gauss-heawood}).  The smaller lattice decomposes as
  \[
  \gauss_5 \simeq 2 D_4^{\gauss} \oplus \II_{1,1}^{\gauss,1+i}
  \]
  where $D_4^\gauss$ and $\II_{1,1}^{\gauss,1+i}$ are as in
  Example~\ref{ex:gauss-heawood}.  The corresponding reflection group
  $\Gamma_5$ is in the Deligne-Mostow list, and the corresponding
  $5$-ball quotient $X_5$ turns out to be the moduli space of
  $8$-points on $\proj^1_\C$ (Eduard Looijenga; private
  communication).
\end{example}

Before continuing on to our final pair of examples, we pause to note that 
the Gaussian lattice pair just described can be converted to a pair of 
quaternionic lattices over the Hurwitz ring $\hurwitz$ of integral quaternions.

\begin{example}[Hurwitzian lattice from the Heawood graph]\label{ex:hurwitz-heawood}
  Let $\vecDelta_7$ be the Heawood graph as in
  Example~\ref{ex:gauss-heawood}.  Just like we built the
  $\eisen$-lattice $\eisen_4 = \eisen \otimes_\Z \Z_\petersen$ from
  the $\Z$-lattice $\Z_\petersen$, using
  Definition~\ref{def:cplx-construction}, one can build a Hermitian
  $\hurwitz$-lattice $\hurwitz_7 = \hurwitz \otimes_{\gauss} \gauss_7$
  from the $\gauss$-lattice $\gauss_7$, where the $\gauss$-valued
  Hermitian form on $\gauss_7$ is naturally extended to a
  $\hurwitz$-valued Hermitian form on $\hurwitz_7$ with signature $(7,
  1)$.  One can show that
  \[
  \hurwitz_7 \simeq \Lambda^{\hurwitz} \oplus
  \II_{1,1}^{\hurwitz, 1 + i} \simeq 3 E_8^{\hurwitz} \oplus
  \II_{1,1}^{\hurwitz, 1 + i}
  \]
  where $\Lambda^{\hurwitz}$ (resp. $E_8^{\hurwitz}$) is the
  quaternionic Leech lattice (resp. quaternionic $E_8$-lattice)
  defined over $\hurwitz$.  The story of this quaternionic lattice and 
  its reflection group acting on quaternionic hyperbolic space is described in
  \cite{Basak07-qLeech}.
\end{example}

As expected, the $8$-cycle in the Heawood graph corresponds to a sublattice.

\begin{example}[Hurwitzian lattice from an $8$-cycle]\label{ex:hurwitz-8-gon}
  The quaternionic Lorentzian representation of the Heawood graph $\hurwitz_7$
  (Example~\ref{ex:hurwitz-heawood}) contains a
  $\hurwitz$-sublattice of $\hurwitz_{5} = \hurwitz \otimes_{\gauss}
  \gauss_5$.  As in Example~\ref{ex:8-gon}, there is a decomposition
  \[
  \hurwitz_5 \simeq 2 E_8^{\hurwitz} \oplus \II_{1,1}^{\hurwitz, 1 + i}.
  \]
\end{example}

The final pair of examples are from graphs with $\lambda_2=\sqrt{3}$.

\begin{figure}
  \centering
  \includegraphics[width=1\linewidth]{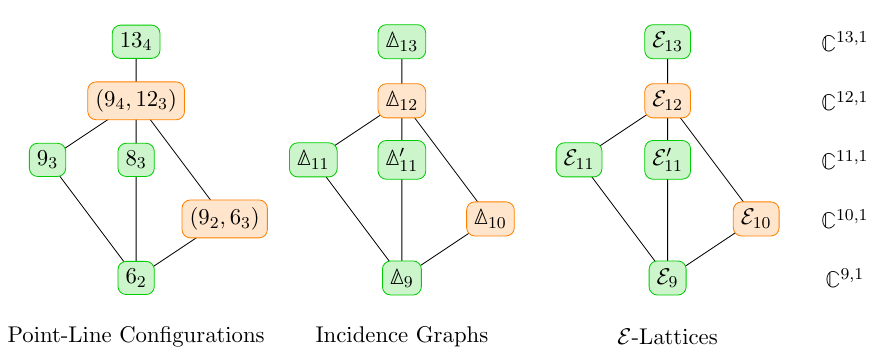}
  \caption{Six point-line configurations (left) that have incidence
    graphs with $\lambda_2=\sqrt{3}$ (middle).  From these we
    construct Lorentzian $\eisen$-lattices in $\C^{n,1}$ (right).}
  \label{fig:eisen-lattices}
\end{figure}

\begin{remark}[Bipartite graphs with $\lambda_2=\sqrt{3}$]
  The six point-line configurations in Table~\ref{tab:sqrt3} have
  incidence graphs with $\lambda_2=\sqrt{3}$, and by
  Corollary~\ref{cor:eisen} these produce Lorentzian
  $\eisen$-lattices.  The inclusion relations among these
  configurations, graphs, and lattices are shown in
  Figure~\ref{fig:eisen-lattices}.  The plan to use the complex
  dimension of the hyperbolic ball as a subscript for these graphs and
  lattice (Notation~\ref{not:dim-index}) fails to distinguish the two
  examples in $\C^{11,1}$, so one receives a prime.  It will be
  interesting to to investigate these sublattices of $\eisen_{13}$ and
  their relationship to the monstrous proposal.
\end{remark}

The third exceptional pair in Table~\ref{tab:sqrt2} comes from the
largest and smallest lattices in Figure~\ref{fig:eisen-lattices}.  The
largest lattice is a special case of
Example~\ref{ex:proj-plane-lattices} and is arguably the most
interesting example of all.

\begin{example}[Eisenstein lattice from the bimonster graph]\label{ex:monstrous}
  Let $\vecDelta_{13} = \inc(\proj^2_\F)$ be the bimonster graph with
  $\F=\F_3$ and $\lambda_2 = \sqrt{3}$. By Corollary~\ref{cor:eisen}
  there is a Lorentzian Eisenstein representation
  $(\eisen_{13},\{\vec{\alpha}_i\})$ of $\vecDelta_{13}$ in
  $\C^{13,1}$.  The signature $(13, 1)$ Eisenstein lattice
  $\eisen_{13}$ and its (projectivized) complex reflection group
  $\Gamma_{13} \subseteq PU(13, 1)$ were first studied by Daniel
  Allcock in \cite{Allcock00-leech}.  He proved that the group
  $\Gamma_{13}$ is arithmetic. Thirteen is the largest dimension $n$
  in which an arithmetic complex reflection group in $U(n,1)$ is known
  and this is the unique known example.  One has
  \[
  \eisen_{13} \simeq \Lambda^{\eisen} \oplus \II_{1,1}^{\eisen,
    \theta} \simeq 3 E_8^{\eisen} \oplus \II_{1,1}^{\eisen, \theta}
  \]
  where $\Lambda^{\eisen}$ (resp. $E_8^{\eisen}$) is the complex Leech
  lattice (resp. complex $E_8$ lattice) defined over $\eisen$ (see
  Definition~\ref{def:important-lattices} or
  \cite{Basak07-bimonster-I} \S 2) and $\II_{1,1}^{\eisen, \theta}$ is
  a hyperbolic cell (Definition~\ref{def:hyp-cell}).
\end{example}

\begin{remark}[Monstrous proposal]
  The (projectivized) complex hyperbolic reflection group $\Gamma_{13}
  = \proj\Refl(\eisen_{13})$ plays central role in Daniel Allcock's
  monstrous proposal conjecture \cite{Allcock09-monstrous} and in the
  work towards proving it \cite{Basak07-bimonster-I, Allcock09-Y555,
    Basak16-bimonster-II, AB16-braid-like, AB18-generators,
    AB25-deflation}.  In these papers, it is proved that that the
  orbifold fundamental group $\pi_1^{\op{orb}}( X_{13}^{\circ},
  \proj\tau)$ of the complex thirteen dimensional ball quotient
  $X_{13}^{\circ}$ is ``braidlike".  It is generated by the
  ``meridians" $\lbrace g_i \colon i \in \vecDelta_{13} \rbrace$
  around the twenty six mirrors $\{\vec{\alpha}_i^{\bot} \colon i \in
  \vecDelta_{13} \}$ closest to $\tau$ and these meridians satisfy all
  the Artin relations of the graph $\vecDelta_{13}$.  So
  $\pi_1^{\op{orb}}(X_{13}^{\circ}, \proj\tau)$ is a quotient of the
  Artin group of $\vecDelta_{13}$.  Roughly speaking $g_i$ is the
  simplest loop in $X_{13}^{\circ}$ based at $\proj\tau$ that go
  around the mirror $\vec{\alpha}_i^{\bot}$ once.  The monstrous
  proposal conjecture predicts that the quotient of the fundamental
  group by the normal subgroup generated by $\lbrace g_i^2 \colon i
  \in \vecDelta_{13} \rbrace$ is the bimonster $(\mbb{M} \times
  \mbb{M})\rtimes \Z/2$ where $\mbb{M}$ is the Monster sporadic simple
  group.  Based on the results of \cite{AB16-braid-like,
    AB18-generators, AB25-deflation}, it is known that this quotient
  is either $((\mbb{M} \times \mbb{M})\rtimes \Z/2)$ or $\Z/2$.
\end{remark}    

The largest induced cycle in the bimonster graph is a $12$-cycle.

\begin{example}[Eisenstein lattice from a $12$-cycle]\label{ex:12-gon}
  The bimonster graph $\vecDelta_{13}$ contains an induced cycle
  $\vecDelta'= \vecDelta_9 = C_{12}$ with the same second eigenvalue
  $\lambda_2=\sqrt{3}$, so by Corollary~\ref{cor:eisen}, the
  $\eisen$-lattice $\eisen_9 = \eisen_{\vecDelta'}$ embeds as a
  sublattice of $\eisen_{13}$ (Example~\ref{ex:monstrous}).  The
  smaller lattice decomposes as
  \[
  \eisen_9 \simeq 2 E_8^{\eisen} \oplus \II_{1,1}^{\eisen, \theta}
  \]
  where $E_8^{\eisen}$ and $\II_{1,1}^{\eisen, \theta}$ are as in
  Example~\ref{ex:monstrous}. The (projectivized) complex reflection
  group $\Gamma_9 = \proj \Refl(\eisen_9) \subseteq \proj U(9, 1)$ is
  the largest example arithmetic lattice in $\proj U(n,1)$ constructed
  by Deligne and Mostow \cite{DeligneMostow86} using monodromy of
  hypergeometric functions and by Thurston \cite{Thurston98} while
  studying the space of triangulations of the $2$-sphere.  The
  $9$-ball quotient $X_9$ turns out to be the the moduli space of
  $12$-points on $\proj^1_\C$.
\end{example}

The moduli space interpretation of $X_9$ plays a key role in the
recent work towards proving the monstrous proposal
\cite{AB25-deflation}.

\begin{remark}[$12$ points on $\proj_\C^1$]
  The automorphism group of the graph $\vecDelta_9$ acts on $\eisen_9$
  and on $\ball\eisen_9 \simeq \hyp^9_\C$ fixing a unique point
  $\proj\tau$ corresponding to the principal eigenvector of its
  adjacency matrix.  In the moduli space picture, the fixed point
  $\proj\tau$ represents a symmetric configuration of twelve points on
  the real $2$-sphere sitting around the equator and the mirrors
  orthogonal to the vectors $\lbrace \vec{\alpha}_i \colon i \in
  C_{12} \rbrace$ corresponds to the simplest possible degenerations
  of this twelve point configuration when two of the neighboring
  points come together.  For more on this, see \cite{AB25-deflation}.
\end{remark}

We conclude with some remarks about the similarities among these examples.

\begin{remark}[Similarities]\label{rem:similarities-summary}
  Let $\vecDelta$ be one of graphs used to construct the eight
  exceptional examples, the six complex examples in
  Table~\ref{tab:lattice-pairs} plus the two quaternionic
  examples. Let $(\clatt,\{\vec{\alpha}_i\})$ be corresponding complex
  or quaternionic Lorentzian representation of these graphs described
  above.  In all of these examples, the spanning vectors
  $\vec{\alpha}_i$'s have norm $r^2 \in \N$ where $r$ is the second
  largest eigenvalues of $\Delta$.  Let $\Gamma = \proj \Refl(\clatt)$
  be the projectivized complex (or quaternionic) reflection group of
  $\clatt$ and let $d$ be the largest integer such that $\Gamma$
  contains $\zeta_d$-reflections in vectors in $\clatt(r^2)$.  In the
  monstrous example \ref{ex:monstrous} and its sub-example described
  in \ref{ex:12-gon}, one has $d = 3$.  In the Gaussian and the
  quaternionic examples described in \ref{ex:gauss-heawood},
  \ref{ex:8-gon}, \ref{ex:hurwitz-heawood}, \ref{ex:hurwitz-8-gon},
  one has $d = 4$.  In Example~\ref{ex:eisen-petersen} related to
  cubic surfaces and its sub-example described in \ref{ex:6-gon} one
  has $d = 6$.  In each of the directed bipartite graphs $\vecDelta$
  discussed in the exceptional examples above, there is an involution
  $\sigma$ of the undirected graph $\Delta$ that interchanges
  $\Delta^+$ and $\Delta^-$.  By Remark~\ref{rem:aut-delta-action} the
  full group $\aut(\Delta)$ acts on $\ball \clatt$ (the last paragraph
  of Definition~\ref{def:midpoint}).  In all eight examples, the
  following facts are known to hold (with one exception):
  \emph{
  \begin{enumerate}
  \item The point $\proj\tau$ is the unique fixed point of the action
    of $\aut(\Delta)$ on $\ball \clatt$.
  \item The mirrors $\{ \vec{\alpha}_i^{\bot} \colon i \in \Delta \}$
    are the closest to $\proj\tau$ (and equidistant from $\proj\tau$)
    among $\{ \vec{\alpha}^{\bot} \colon \alpha \in \clatt(r^2) \}$,
    the mirrors of $\Gamma$.
  \item The complex or quaternionic reflection group $\Gamma$ has
    finite index in the automorphism group of $\clatt$, i.e. $\Gamma$
    is arithmetic.
  \item There is a group homomorphism $\cox(\Delta,d)\to \Gamma$ that
    sends the Coxeter generators to the order~$d$ complex or
    quaternionic reflections in the roots $\{\vec{\alpha}_i \colon i
    \in \Delta \}$.  In particular, these reflections braid (commute)
    when the corresponding vertices in $\Delta$ are (are not)
    connected by an edge.
  \item With the possible exception of Example~\ref{ex:gauss-heawood},
    the reflections in the roots $\{\vec{\alpha}_i \colon i \in \Delta
    \}$ generate $\Gamma$.  So the complex or quaternionic reflection
    group $\Gamma$ is a quotient of $\cox(\Delta, d)$.
  \end{enumerate}
  }
  \noindent

  It is likely that the order~$4$ complex reflections in the roots
  $\lbrace \vec{\alpha}_i \colon i \in \vecDelta_7 \rbrace$ in
  Example~\ref{ex:gauss-heawood} generate $\Gamma_7$, but no one has
  verified this. The proofs of these facts mostly proceed on a
  case-by-case basis.  Proofs are given in the references for the
  individual examples, although there are a few cases that remain
  unpublished.
\end{remark}

The properties listed in Remark~\ref{rem:similarities-summary} lead to
the following analogy: the point $\proj\tau$ behaves like the Weyl
vector and the spanning vectors $\vec{\alpha}_i$ behave like the
simple roots of $\Gamma$.  Even though this analogy is rough and
imprecise, it has been extremely useful when trying to guess results
and proofs in the work on monstrous proposal and the work on
fundamental group of moduli space of cubic surfaces.

\appendix

\section{A little linear algebra}\label{app:algebra}

In this self-contained appendix we describe a general construction
that starts with a commutative ring $R$ with an isomorphism $\sigma_R:
R \to R$ denoted $x \mapsto \bar{x}$ satisfying $\sigma_R^2 = \op{id}$
and a Hermitian matrix $B \in R^{v \times v}$ and builds an $R$-module
$R_{v,B}$ with a ``nondegenerate" $R$-valued Hermitian form. We prove
some basic linear algebra results about $R_{v,B}$.  The Lorentzian
representation of a graph $\Delta$ studied in this article is a
special case of this construction.  We have two kinds of
applications. In the ``real case" studied in the first few sections,
$R$ is often $\Z$ or $\R$ and $B = (r I - A)$ where $A$ is the
adjacency matrix of a connected graph $\Delta$ and $r$ is the second
largest eigenvalue of $A$.  The ``complex case" studied in the last
few sections is slightly more involved; 
see Definitions 
\ref{def:c-spaces}
and 
\ref{def:cplx-construction}.  The generality actually makes
the definitions and arguments simpler so we have chosen to present
them in this manner.

\begin{lemma}\label{lemma-module-rank-equals-matrix-rank}
  Let $R$ be a PID. Let $B \in R^{b \times a}$.  Let $N \subseteq R^b$
  be the $R$-module spanned by the columns of $B$.  Then the
  $R$-module $N$ is free and its rank as an $R$-module is equal to the
  rank of the matrix $B$, i.e., $\op{rank}_R(N) = \op{rank}(B)$.
\end{lemma}

\begin{proof}
  Since $N \leq R^b$, it follows that $N$ is free  
    and $\op{rank}_R(N) \leq b$.  Let $u_1, \dotsb, u_n$
  be a $R$-basis of $N$.  Let $F$ be the fraction field of $R$ and let
  $V$ be the $F$-span of the columns of $B$.  Then $u_1, \dotsb, u_n$
  span $V$. On the other hand note that $u_1, \dotsb, u_n$ are also
  linearly independent over $F$ since any nontrivial dependence
  relation over $F$ would yield a nontrivial dependence relation over
  $R$ by clearing denominators. So $ \op{rank}(B)= \op{dim}_F(V) = n$. 
\end{proof}

\begin{definition}[Hermitian $R$-modules]\label{definition-Hermitian-R-modules}
  An \emph{involution} of a commutative ring $R$ means a ring
  homomorphism $\sigma_R: R \to R$, denoted by $x \mapsto \bar{x}$
  such that $\sigma_R^2 = \op{id}_R$.  Let $R$ be a commutative ring
  with an involution denoted by $x \mapsto \bar{x}$. A matrix $B \in
  R^{v \times v}$ is called Hermitian if $B^{\msf{H}} :=
  \bar{B}^{\msf{T}} = B$.  Note that we allow the possibility that
  $\sigma_R = \op{id}_R$.  As usual, if $\sigma_R = \op{id}_R$,
  Hermitian matrix/form means symmetric matrix/bilinear form.  Let $B
  \in R^{v \times v}$ be a Hermitian matrix for some $v \in \N$.  Let
  $R^v_B$ denote the free $R$-module $R^v$ together with the
  $R$-valued Hermitian form
  \[
  \ip{u}{w} 
  = u^{\msf{H}} B w 
  \text{\; for \;} u, w \in R^v.
  \]
  Let 
  \[
  \op{rad} (R^v_B)
  := \lbrace w \in R^v \colon \ip{R^v_B}{w} = 0 \rbrace.
  \]
  Note that
  \[
  \op{rad}(R^v_B) 
  = \op{ker} \lbrace  B : R^v \to R^v \rbrace.
  \]
  Define 
  \[
  R_{v,B}  = R^v_B/ \op{rad}( R^{v}_B)
  \]
  together with the nondegenerate (i.e. zero radical) Hermitian form
  induced from $R^v_B$, denoted again by $\ip{u}{w}$.
\end{definition}

\begin{lemma}\label{lemma-Lorentzian-rep-is-image-of-Gram}
  Let $B(R^v)$ be the image of the map $B: R^v \to R^v$, i.e.  the
  $R$-span of the columns of $B$.  Then one has an $R$-module
  isomorphism
  \[
  R_{v,B} \simeq B(R^{v}).
  \] 
  If $R$ is a PID, then $R_{v,B}$ is a free $R$-module and
  $\op{rank}_{R}( R_{v,B}) = \op{rank}(B)$.
\end{lemma}

\begin{proof}
  The linear map $B: R^v \to R^v$ induces a short exact sequence of
  $R$-modules
  \begin{equation}
    0 \to \op{rad}  (R^{v}_B)  \to R^v 
    \xrightarrow{B} B(R^v) \to 0.
    \label{eq-ses-LR}
  \end{equation}
  So $R_{v,B}\simeq B(R^v)$ as $R$-modules
  and the lemma follows from 
  Lemma~\ref{lemma-module-rank-equals-matrix-rank}.
\end{proof}

\begin{lemma}[change of rings]\label{lemma-change-of-rings}
  Let $\varphi: R \to S$ be a homomorphism of commutative rings with
  involutions\footnote{All commutative rings have $1$'s, homomorphisms
  are unit preserving, compatible with involutions.}.  Let $B \in R^{v
    \times v}$ such that $B^{\msf{H}} = B$ and $C = \varphi(B)$.  Any
  $S$-module is a $R$-module via the homomorphism $\varphi$.  Let
  $\varphi^v: R^v \to S^v$ be the map induced by $\varphi$.  There is
  a unique $R$-module map $\varphi_v: R_{v,B} \to S_{v,C}$ such that
  \[
  u + \op{rad}( R^v_B)
  \xmapsto{\varphi_v} 
  \varphi^v(u) + \op{rad}( S^v_C) 
  \text{\; for all \;} u \in R^v
  \]
  and one has the following commutative diagram of $R$-modules:
  \[
  \begin{tikzcd}[column sep=large]
    0 \arrow[r] &
    \op{rad}(R^v_B)
    \arrow[r, hook]
    \arrow[d, "\varphi^v"] &
    R^v_B
    \arrow[r, "\pi "]
    \arrow[d, "\varphi^v"] &
    R_{v,B}
    \arrow[r]
    \arrow[d, "\varphi_v"] &
    0
    \\
    0 \arrow[r] &
    \op{rad}(S^v_C)
    \arrow[r, hook] &
    S^v_C
    \arrow[r, "\pi "] &
    S_{v,C}
    \arrow[r] &
    0 .
  \end{tikzcd}
  \]
  Here $\pi$ denotes the quotient maps and the rows are exact.
  Furthermore the maps $\varphi^v: R^v_B \to S^v_C$ and $\varphi_v:
  R_{v,B} \to S_{v,C}$ are compatible with the forms. In particular,
  \[
  \varphi\left( \ip{u}{w} \right)
  = \ip{ \varphi_{v}(u)}{ \varphi_{v}(w)}
  \text{\; for all \;} u, w \in R_{v,B}.
  \]
  As a special case, suppose $\iota: R \hookrightarrow S$ is an
  inclusion of commutative rings with involutions.  Then $R^{v}_B \cap
  \op{rad}( S^{v}_B) = \op{rad}( R^{v}_B)$.  One has a unique
  injective $R$-module map $\iota_{v} : R_{v,B} \hookrightarrow
  S_{v,B}$ such that
  \[
  u + \op{rad}( R^{v}_B) 
  \xmapsto{\iota_{v}} 
  u + \op{rad}( S^{v}_B) 
  \text{\; for all \;} u \in R^{v},
  \]
  and one has the following commutative diagram of $R$-modules:
  \[
  \begin{tikzcd}[column sep=large]
    0 \arrow[r] &
    \op{rad}(R^v_B)
    \arrow[r, hook]
    \arrow[d, hook] &
    R^v_B
    \arrow[r, "\pi "]
    \arrow[d, hook] &
    R_{v,B}
    \arrow[r]
    \arrow[d, hook, "\iota_v"] &
    0
    \\
    0 \arrow[r] &
    \op{rad}(S^v_B)
    \arrow[r, hook] &
    S^v_B
    \arrow[r, "\pi "] &
    S_{v,B}
    \arrow[r] &
    0 
  \end{tikzcd}
  \]
  where the rows are exact and the vertical maps are all injective.
  Furthermore, the map $\iota_{v}: R_{v,B} \to S_{v,B}$ is form
  preserving.
\end{lemma}

\begin{proof}[sketch of proof]
  Let $u \in \op{rad}( R^{v}_B)$. Then $B u = 0$, so applying
  $\varphi$ we obtain $C \varphi^v(u) = 0$, so $ \varphi^{v}(u) \in
  \op{rad}( S^{v}_C)$.  This shows that $\varphi^{v}$ carries
  $\op{rad}( R^{v}_B)$ into $\op{rad}( S^{v}_C)$.  Now consider the
  diagram without the vertical map $\varphi_{v}$.  It has exact rows
  by definition and the left square commutes.  Since the left square
  commutes, the map $\pi \circ \varphi^{v} : R^{v}_B \to S_{v,C}$ has
  $\op{rad}(R^{v}_B)$ in its kernel, so it induces a unique map
  $\varphi_{v} : R_{v,B} \to S_{v,C}$ such that the right square
  commutes.  The compatibility of $\varphi_{v}$ with the forms is
  immediate.

  In the special case, when $\iota : R \to S$ is an inclusion, note
  that $\op{rad}( S^{v}_B)$ and $\op{rad}( R^{v}_B)$ are both kernels
  of the same linear map $B$ on $S^{v}$ and $R^{v}$ respectively, so
  the equality $R^{v}_B \cap \op{rad}( S^{v}_B) = \op{rad}( R^{v}_B)$
  holds.  This implies that the map $\iota_{v}$ is injective.  The
  rest of the special case is immediate from the general case.
\end{proof}

\begin{lemma}[Base change]\label{lemma-base-change}
  Let $R\subseteq S$ be an inclusion of commutative rings with
  involutions.  Assume $S$ is flat over $R$. Let $B \in R^{v \times
    v}$ be a Hermitian matrix.  Then the natural isomorphism $S
  \otimes_R R^v_B \to S^v_B$ induces an isomorphism of short exact
  sequences of $S$-modules
  \[
  \begin{tikzcd}[]
    0 \arrow[r] &
    S\otimes_R \op{rad}(R^v_B)
    \arrow[r]
    \arrow[d, "\simeq"'] &
    S\otimes_R R^v_B
    \arrow[r]
    \arrow[d, "\simeq"'] &
    S\otimes_R R_{v,B}
    \arrow[r]
    \arrow[d, "\simeq"'] &
    0
    \\
    0 \arrow[r] &
    \op{rad}(S^v_B)
    \arrow[r, hook] &
    S^v_B
    \arrow[r] &
    S_{v,B}
    \arrow[r] &
    0 .
  \end{tikzcd}
  \]
  In particular, the natural map
  \[
  S\otimes_R R_{v,B}\longrightarrow S_{v,B},
  \;\; \text{given by \;}
  s\otimes(u+\op{rad}(R^v_B))\longmapsto s u+\op{rad}(S^v_B),
  \]
  is an isomorphism.
\end{lemma}

\begin{proof}
  Due to flatness, $(S \otimes_R \cdot)$ preserves exactness. So the
  exact sequence
  \[
  0\longrightarrow \op{rad}(R^v_B) \longrightarrow R^v_B
  \xrightarrow{B} R^v_B
  \]
  yields an exact sequence
  \[
  0 \longrightarrow S\otimes_R \op{rad}(R^v_B)
  \longrightarrow S\otimes_R R^v_B
  \xrightarrow{1\otimes B} S\otimes_R R^v_B.
  \]
  Under the canonical identification
  $
  \mu: S \otimes_R R^v_B \simeq S^v_B,
  $
  given by $(s \otimes u) \mapsto s u$, 
  the map $1\otimes B$ on $S\otimes_R R^v_B$ corresponds to 
  multiplication by $B$ on $S^v_B$. 
  In other words, we have a commutative square
  \[
  \begin{tikzcd}[column sep=large]
    S\otimes_R R^v_B
    \arrow[r, "1\otimes B"]
    \arrow[d, "\mu"', "\simeq"] &
    S\otimes_R R^v_B
    \arrow[d, "\mu", "\simeq"']
    \\
    S^v_B
    \arrow[r, "B"'] &
    S^v_B .
  \end{tikzcd}
  \]
  Indeed, if $s\in S$ and $u\in R^v$, then $ \mu\bigl((1\otimes
  B)(s\otimes u)\bigr) =\mu(s\otimes B u) =s B u =B(s u).  $ So $\mu$
  isomorphically maps the kernel of $(1 \otimes B)$ to the kernel of
  $B$ and we obtain the following commutative diagram with exact rows
  and where the vertical maps are isomorphisms:
  \[
  \begin{tikzcd}[column sep=large]
    0 \arrow[r] &
    S\otimes_R \op{rad}(R^v_B)
    \arrow[r]
    \arrow[d, "\mu"', "\simeq"] &
    S\otimes_R R^v_B
    \arrow[r, "1\otimes B"]
    \arrow[d, "\mu"', "\simeq"] &
    S\otimes_R R^v_B
    \arrow[d, "\mu", "\simeq"']
    \\
    0 \arrow[r] &
    \op{rad}( S^v_B)
    \arrow[r, hook] &
    S^v_B
    \arrow[r, "B"'] &
    S^v_B .
  \end{tikzcd}
  \]
  Now, applying $(S \otimes_R \cdot)$ to the short exact sequence
  \[
  0 \longrightarrow \op{rad}(R^v_B) \longrightarrow R^v_B
  \longrightarrow R_{v,B} \longrightarrow 0
  \]
  we obtain the following commutative diagram with exact rows:
  \[
  \begin{tikzcd}[]
    0 \arrow[r] &
    S\otimes_R \op{rad}(R^v_B)
    \arrow[r]
    \arrow[d, "\simeq"'] &
    S\otimes_R R^v_B
    \arrow[r]
    \arrow[d, "\simeq"'] &
    S\otimes_R R_{v,B}
    \arrow[r]
    &
    0
    \\
    0 \arrow[r] &
    \op{rad}(S^v_B)
    \arrow[r, hook] &
    S^v_B
    \arrow[r] &
    S_{v,B}
    \arrow[r] &
    0 .
  \end{tikzcd}
  \]
  It follows that there is a unique induced isomorphism between the
  cokernels, that is, from $S \otimes_R R_{v,B}$ to $S_{v,B}$, making
  the required commutative diagram.
\end{proof}


\bibliographystyle{amsalpha}
\bibliography{lor-lat-from-graphs}

@article {McMullen2002,
    AUTHOR = {McMullen, Curtis T.},
     TITLE = {Coxeter groups, {S}alem numbers and the {H}ilbert metric},
   JOURNAL = {Publ. Math. Inst. Hautes \'Etudes Sci.},
  FJOURNAL = {Publications Math\'ematiques. Institut de Hautes \'Etudes
              Scientifiques},
    NUMBER = {95},
      YEAR = {2002},
     PAGES = {151--183},
      ISSN = {0073-8301,1618-1913},
   MRCLASS = {20F55 (11R06 37C85 37D40)},
  MRNUMBER = {1953192},
MRREVIEWER = {Michael\ J.\ Mossinghoff},
       DOI = {10.1007/s102400200001},
       URL = {https://doi.org/10.1007/s102400200001},
}

@incollection {Tits57,
    AUTHOR = {Tits, J.},
     TITLE = {Sur les analogues alg\'ebriques des groupes semi-simples
              complexes},
 BOOKTITLE = {Colloque d'alg\`ebre sup\'erieure, tenu \`a{} {B}ruxelles du
              19 au 22 d\'ecembre 1956},
    SERIES = {Centre Belge de Recherches Math\'ematiques},
     PAGES = {261--289},
 PUBLISHER = {\'Etablissements Ceuterick, Louvain},
      YEAR = {1957},
   MRCLASS = {50.00 (22.00)},
  MRNUMBER = {108765},
MRREVIEWER = {H.\ Freudenthal},
}

@book {Roman08,
    AUTHOR = {Roman, Steven},
     TITLE = {Advanced linear algebra},
    SERIES = {Graduate Texts in Mathematics},
    VOLUME = {135},
   EDITION = {Third},
 PUBLISHER = {Springer, New York},
      YEAR = {2008},
     PAGES = {xviii+522},
      ISBN = {978-0-387-72828-5},
   MRCLASS = {15-01},
  MRNUMBER = {2344656},
}

@article {Wilson89,
    AUTHOR = {Wilson, Robert A.},
     TITLE = {Vector stabilizers and subgroups of {L}eech lattice groups},
   JOURNAL = {J. Algebra},
  FJOURNAL = {Journal of Algebra},
    VOLUME = {127},
      YEAR = {1989},
    NUMBER = {2},
     PAGES = {387--408},
      ISSN = {0021-8693,1090-266X},
   MRCLASS = {20D08 (20C20 20E28)},
  MRNUMBER = {1028461},
       DOI = {10.1016/0021-8693(89)90260-3},
       URL = {https://doi.org/10.1016/0021-8693(89)90260-3},
}

@book {Wilson09,
    AUTHOR = {Wilson, Robert A.},
     TITLE = {The finite simple groups},
    SERIES = {Graduate Texts in Mathematics},
    VOLUME = {251},
 PUBLISHER = {Springer-Verlag London, Ltd., London},
      YEAR = {2009},
     PAGES = {xvi+298},
      ISBN = {978-1-84800-987-5},
   MRCLASS = {20D05},
  MRNUMBER = {2562037},
MRREVIEWER = {Gernot\ Stroth},
       DOI = {10.1007/978-1-84800-988-2},
       URL = {https://doi.org/10.1007/978-1-84800-988-2},
}

@Unpublished{BaMc-graphs,
  author = 	 {Basak, Tathagata and McCammond, Jon},
  title = 	 {Equipartition among regular subgraphs},
  note = 	 {Preprint 2026},
}

@article {Seidel68,
    AUTHOR = {Seidel, J. J.},
     TITLE = {Strongly regular graphs with {$(-1,\,1,\,0)$} adjacency matrix
              having eigenvalue {$3$}},
   JOURNAL = {Linear Algebra Appl.},
  FJOURNAL = {Linear Algebra and its Applications},
    VOLUME = {1},
      YEAR = {1968},
     PAGES = {281--298},
      ISSN = {0024-3795,1873-1856},
   MRCLASS = {05.40},
  MRNUMBER = {234861},
MRREVIEWER = {B.\ R.\ Heap},
       DOI = {10.1016/0024-3795(68)90008-6},
       URL = {https://doi.org/10.1016/0024-3795(68)90008-6},
}

@article {Neumaier79,
    AUTHOR = {Neumaier, A.},
     TITLE = {Strongly regular graphs with smallest eigenvalue {$-m$}},
   JOURNAL = {Arch. Math. (Basel)},
  FJOURNAL = {Archiv der Mathematik},
    VOLUME = {33},
      YEAR = {1979/80},
    NUMBER = {4},
     PAGES = {392--400},
      ISSN = {0003-889X,1420-8938},
   MRCLASS = {05C99 (05C50)},
  MRNUMBER = {564298},
MRREVIEWER = {E.\ A.\ Nordhaus},
       DOI = {10.1007/BF01222774},
       URL = {https://doi.org/10.1007/BF01222774},
}

@article {Krammer09,
    AUTHOR = {Krammer, Daan},
     TITLE = {The conjugacy problem for {C}oxeter groups},
   JOURNAL = {Groups Geom. Dyn.},
  FJOURNAL = {Groups, Geometry, and Dynamics},
    VOLUME = {3},
      YEAR = {2009},
    NUMBER = {1},
     PAGES = {71--171},
      ISSN = {1661-7207,1661-7215},
   MRCLASS = {20F55 (20F10)},
  MRNUMBER = {2466021},
MRREVIEWER = {Robert\ H.\ Gilman},
       DOI = {10.4171/GGD/52},
       URL = {https://doi.org/10.4171/GGD/52},
}

@book {BjornerBrenti05,
    AUTHOR = {Bj\"orner, Anders and Brenti, Francesco},
     TITLE = {Combinatorics of {C}oxeter groups},
    SERIES = {Graduate Texts in Mathematics},
    VOLUME = {231},
 PUBLISHER = {Springer, New York},
      YEAR = {2005},
     PAGES = {xiv+363},
      ISBN = {978-3540-442387; 3-540-44238-3},
   MRCLASS = {05-01 (05E15 20F55)},
  MRNUMBER = {2133266},
MRREVIEWER = {Jian-yi\ Shi},
}

@book {Ebeling13,
    AUTHOR = {Ebeling, Wolfgang},
     TITLE = {Lattices and codes},
    SERIES = {Advanced Lectures in Mathematics},
   EDITION = {Third},
      NOTE = {A course partially based on lectures by Friedrich Hirzebruch},
 PUBLISHER = {Springer Spektrum, Wiesbaden},
      YEAR = {2013},
     PAGES = {xvi+167},
      ISBN = {978-3-658-00359-3; 978-3-658-00360-9},
   MRCLASS = {11H31 (11T71 51F15)},
  MRNUMBER = {2977354},
MRREVIEWER = {Caleb\ McKinley\ Shor},
       DOI = {10.1007/978-3-658-00360-9},
       URL = {https://doi.org/10.1007/978-3-658-00360-9},
}

@article {Goldman25,
    AUTHOR = {Goldman, Katherine M.},
     TITLE = {{${\rm CAT}(0)$} and cubulated {S}hephard groups},
   JOURNAL = {J. Lond. Math. Soc. (2)},
  FJOURNAL = {Journal of the London Mathematical Society. Second Series},
    VOLUME = {111},
      YEAR = {2025},
    NUMBER = {1},
     PAGES = {Paper No. e70050, 33},
      ISSN = {0024-6107,1469-7750},
   MRCLASS = {20F65 (57M60)},
  MRNUMBER = {4840638},
MRREVIEWER = {Lawk\ Mineh},
       DOI = {10.1112/jlms.70050},
       URL = {https://doi.org/10.1112/jlms.70050},
}

@book {AdkinsWeintraub92,
    AUTHOR = {Adkins, William A. and Weintraub, Steven H.},
     TITLE = {Algebra},
    SERIES = {Graduate Texts in Mathematics},
    VOLUME = {136},
      NOTE = {An approach via module theory},
 PUBLISHER = {Springer-Verlag, New York},
      YEAR = {1992},
     PAGES = {x+526},
      ISBN = {0-387-97839-9},
   MRCLASS = {00A05 (13-01 15-01 16-01 20-01)},
  MRNUMBER = {1181420},
MRREVIEWER = {Gerard\ M.\ Enright},
       DOI = {10.1007/978-1-4612-0923-2},
       URL = {https://doi.org/10.1007/978-1-4612-0923-2},
}

@article {BonadrenkoRadchenko13,
    AUTHOR = {Bondarenko, Andriy V. and Radchenko, Danylo V.},
     TITLE = {On a family of strongly regular graphs with {$\lambda=1$}},
   JOURNAL = {J. Combin. Theory Ser. B},
  FJOURNAL = {Journal of Combinatorial Theory. Series B},
    VOLUME = {103},
      YEAR = {2013},
    NUMBER = {4},
     PAGES = {521--531},
      ISSN = {0095-8956,1096-0902},
   MRCLASS = {05E30},
  MRNUMBER = {3071380},
MRREVIEWER = {Kaishun\ Wang},
       DOI = {10.1016/j.jctb.2013.05.005},
       URL = {https://doi.org/10.1016/j.jctb.2013.05.005},
}

@article {Bondarenko14,
    AUTHOR = {Bondarenko, Andriy},
     TITLE = {On {B}orsuk's conjecture for two-distance sets},
   JOURNAL = {Discrete Comput. Geom.},
  FJOURNAL = {Discrete \& Computational Geometry. An International Journal
              of Mathematics and Computer Science},
    VOLUME = {51},
      YEAR = {2014},
    NUMBER = {3},
     PAGES = {509--515},
      ISSN = {0179-5376,1432-0444},
   MRCLASS = {52C17 (05E30)},
  MRNUMBER = {3201240},
MRREVIEWER = {Marek\ Lassak},
       DOI = {10.1007/s00454-014-9579-4},
       URL = {https://doi.org/10.1007/s00454-014-9579-4},
}

@book {BrouwerVanMaldeghem22,
    AUTHOR = {Brouwer, Andries E. and Van Maldeghem, H.},
     TITLE = {Strongly regular graphs},
    SERIES = {Encyclopedia of Mathematics and its Applications},
    VOLUME = {182},
 PUBLISHER = {Cambridge University Press, Cambridge},
      YEAR = {2022},
     PAGES = {xvii+462},
      ISBN = {978-1-316-51203-6},
   MRCLASS = {05-02 (05E30 51A50 51E24 94Bxx)},
  MRNUMBER = {4350112},
MRREVIEWER = {Ulrich\ Tamm},
       DOI = {10.1017/9781009057226},
       URL = {https://doi.org/10.1017/9781009057226},
}

@article {BPR17,
    AUTHOR = {Bondarenko, A. V. and Prymak, A. and Radchenko, D.},
     TITLE = {Non-existence of {$(76,30,8,14)$} strongly regular graph},
   JOURNAL = {Linear Algebra Appl.},
  FJOURNAL = {Linear Algebra and its Applications},
    VOLUME = {527},
      YEAR = {2017},
     PAGES = {53--72},
      ISSN = {0024-3795,1873-1856},
   MRCLASS = {05C25 (05C50)},
  MRNUMBER = {3647469},
MRREVIEWER = {Gregory\ Tycho\ Markowsky},
       DOI = {10.1016/j.laa.2017.03.033},
       URL = {https://doi.org/10.1016/j.laa.2017.03.033},
}

@incollection {BMPRV18,
    AUTHOR = {Bondarenko, Andriy and Mellit, Anton and Prymak, Andriy and
              Radchenko, Danylo and Viazovska, Maryna},
     TITLE = {There is no strongly regular graph with parameters (460, 153,
              32, 60)},
 BOOKTITLE = {Contemporary computational mathematics---a celebration of the
              80th birthday of {I}an {S}loan. {V}ol. 1, 2},
     PAGES = {131--134},
 PUBLISHER = {Springer, Cham},
      YEAR = {2018},
      ISBN = {978-3-319-72455-3; 978-3-319-72456-0},
   MRCLASS = {05E30},
  MRNUMBER = {3822232},
}

@article {AB16-braid-like,
    AUTHOR = {Allcock, Daniel and Basak, Tathagata},
     TITLE = {Geometric generators for braid-like groups},
   JOURNAL = {Geom. Topol.},
  FJOURNAL = {Geometry \& Topology},
    VOLUME = {20},
      YEAR = {2016},
    NUMBER = {2},
     PAGES = {747--778},
      ISSN = {1465-3060,1364-0380},
   MRCLASS = {57M05 (20F36 32S22 52C35)},
  MRNUMBER = {3493096},
MRREVIEWER = {Jon\ McCammond},
       DOI = {10.2140/gt.2016.20.747},
       URL = {https://doi.org/10.2140/gt.2016.20.747},
}

@article {AB18-generators,
    AUTHOR = {Allcock, Daniel and Basak, Tathagata},
     TITLE = {Generators for a complex hyperbolic braid group},
   JOURNAL = {Geom. Topol.},
  FJOURNAL = {Geometry \& Topology},
    VOLUME = {22},
      YEAR = {2018},
    NUMBER = {6},
     PAGES = {3435--3500},
      ISSN = {1465-3060,1364-0380},
   MRCLASS = {57M05 (20F36 32S22 52C35)},
  MRNUMBER = {3858768},
MRREVIEWER = {Sang-hyun\ Kim},
       DOI = {10.2140/gt.2018.22.3435},
       URL = {https://doi.org/10.2140/gt.2018.22.3435},
}

@article {AB25-deflation,
    AUTHOR = {Allcock, Daniel and Basak, Tathagata},
     TITLE = {The {D}eligne-{M}ostow 9-ball, and the {M}onster},
   JOURNAL = {Geom. Topol.},
  FJOURNAL = {Geometry \& Topology},
    VOLUME = {29},
      YEAR = {2025},
    NUMBER = {2},
     PAGES = {791--828},
      ISSN = {1465-3060,1364-0380},
   MRCLASS = {20F36 (20D08 22E40)},
  MRNUMBER = {4879367},
MRREVIEWER = {A.\ S.\ Kondrat\cprime ev},
       DOI = {10.2140/gt.2025.29.791},
       URL = {https://doi.org/10.2140/gt.2025.29.791},
}

@misc {ABL,
  author = {Allcock, Daniel and Basak, Tathagata and Looijenga, Eduard},
   title = {The moduli space of cubic surfaces and its orbifold fundamental group},
  YEAR          = {2026},
  eprint        = {2610.16080},
  archivePrefix = {arXiv},
  primaryClass  = {math.AG},
  url           = {https://arxiv.org/abs/2610.16080},
  note   = {Preprint, arXiv:2610.16080}
}

@misc {B-monodromy,
  author = {Basak, Tathagata},
   title = {Petersen graph and monodromy of the 27 lines on the Clebsch surface},
  YEAR          = {2026},
  eprint        = {2607.01878},
  archivePrefix = {arXiv},
  primaryClass  = {math.AG},
  url           = {https://arxiv.org/abs/2607.01878},
  note   = {Preprint, arXiv:2607.01878}
}

@article {ACT02-cubic-surfaces,
    AUTHOR = {Allcock, Daniel and Carlson, James A. and Toledo, Domingo},
     TITLE = {The complex hyperbolic geometry of the moduli space of cubic
              surfaces},
   JOURNAL = {J. Algebraic Geom.},
  FJOURNAL = {Journal of Algebraic Geometry},
    VOLUME = {11},
      YEAR = {2002},
    NUMBER = {4},
     PAGES = {659--724},
      ISSN = {1056-3911,1534-7486},
   MRCLASS = {32G20 (11F23 14J10 14J15 14K30 32M15 32Q45)},
  MRNUMBER = {1910264},
MRREVIEWER = {I.\ Dolgachev},
       DOI = {10.1090/S1056-3911-02-00314-4},
       URL = {https://doi.org/10.1090/S1056-3911-02-00314-4},
}

@article {Allcock00-leech,
    AUTHOR = {Allcock, Daniel},
     TITLE = {The {L}eech lattice and complex hyperbolic reflections},
   JOURNAL = {Invent. Math.},
  FJOURNAL = {Inventiones Mathematicae},
    VOLUME = {140},
      YEAR = {2000},
    NUMBER = {2},
     PAGES = {283--301},
      ISSN = {0020-9910,1432-1297},
   MRCLASS = {11H06 (11F55)},
  MRNUMBER = {1756997},
       DOI = {10.1007/s002220050363},
       URL = {https://doi.org/10.1007/s002220050363},
}

@article {Allcock09-Y555,
    AUTHOR = {Allcock, Daniel},
     TITLE = {On the {$Y_{555}$} complex reflection group},
   JOURNAL = {J. Algebra},
  FJOURNAL = {Journal of Algebra},
    VOLUME = {322},
      YEAR = {2009},
    NUMBER = {5},
     PAGES = {1454--1465},
      ISSN = {0021-8693,1090-266X},
   MRCLASS = {20F55},
  MRNUMBER = {2543618},
MRREVIEWER = {Sait\ Hal\i c\i o\u glu},
       DOI = {10.1016/j.jalgebra.2009.05.027},
       URL = {https://doi.org/10.1016/j.jalgebra.2009.05.027},
}

@incollection {Allcock09-monstrous,
    AUTHOR = {Allcock, Daniel},
     TITLE = {A monstrous proposal},
 BOOKTITLE = {Groups and symmetries},
    SERIES = {CRM Proc. Lecture Notes},
    VOLUME = {47},
     PAGES = {17--24},
 PUBLISHER = {Amer. Math. Soc., Providence, RI},
      YEAR = {2009},
      ISBN = {978-0-8218-4481-6},
   MRCLASS = {20D08 (20F05)},
  MRNUMBER = {2500552},
MRREVIEWER = {Jamshid\ Moori},
       DOI = {10.1090/crmp/047/03},
       URL = {https://doi.org/10.1090/crmp/047/03},
}

@article {Basak07-bimonster-I,
    AUTHOR = {Basak, Tathagata},
     TITLE = {The complex {L}orentzian {L}eech lattice and the {B}imonster},
   JOURNAL = {J. Algebra},
  FJOURNAL = {Journal of Algebra},
    VOLUME = {309},
      YEAR = {2007},
    NUMBER = {1},
     PAGES = {32--56},
      ISSN = {0021-8693,1090-266X},
   MRCLASS = {20K30 (11H06 20D08)},
  MRNUMBER = {2301231},
MRREVIEWER = {Robert\ A.\ Wilson},
       DOI = {10.1016/j.jalgebra.2006.05.033},
       URL = {https://doi.org/10.1016/j.jalgebra.2006.05.033},
}

@article {Basak07-qLeech,
    AUTHOR = {Basak, Tathagata},
     TITLE = {Reflection group of the quaternionic {L}orentzian {L}eech
              lattice},
   JOURNAL = {J. Algebra},
  FJOURNAL = {Journal of Algebra},
    VOLUME = {309},
      YEAR = {2007},
    NUMBER = {1},
     PAGES = {57--68},
      ISSN = {0021-8693,1090-266X},
   MRCLASS = {11H06 (20F55)},
  MRNUMBER = {2301232},
MRREVIEWER = {Robert\ A.\ Wilson},
       DOI = {10.1016/j.jalgebra.2006.05.032},
       URL = {https://doi.org/10.1016/j.jalgebra.2006.05.032},
}

@article {Basak14-proj-planes,
    AUTHOR = {Basak, Tathagata},
     TITLE = {Modular lattices from finite projective planes},
   JOURNAL = {J. Th\'eor. Nombres Bordeaux},
  FJOURNAL = {Journal de Th\'eorie des Nombres de Bordeaux},
    VOLUME = {26},
      YEAR = {2014},
    NUMBER = {2},
     PAGES = {269--279},
      ISSN = {1246-7405,2118-8572},
   MRCLASS = {11H56 (11E12 11E39 51E15)},
  MRNUMBER = {3320480},
MRREVIEWER = {Gabriele\ Nebe},
       DOI = {10.5802/jtnb.867},
       URL = {https://doi.org/10.5802/jtnb.867},
}

@article {Basak16-bimonster-II,
    AUTHOR = {Basak, Tathagata},
     TITLE = {The complex {L}orentzian {L}eech lattice and the {B}imonster
              ({II})},
   JOURNAL = {Trans. Amer. Math. Soc.},
  FJOURNAL = {Transactions of the American Mathematical Society},
    VOLUME = {368},
      YEAR = {2016},
    NUMBER = {6},
     PAGES = {4171--4195},
      ISSN = {0002-9947,1088-6850},
   MRCLASS = {11H06 (11H56 20F36 20F55)},
  MRNUMBER = {3453368},
MRREVIEWER = {Gabriele\ Nebe},
       DOI = {10.1090/tran/6558},
       URL = {https://doi.org/10.1090/tran/6558},
}

@article {Basak20-barnes-wall,
    AUTHOR = {Basak, Tathagata},
     TITLE = {A new complex reflection group in {$PU(9,1)$} and the
              {B}arnes-{W}all lattice},
   JOURNAL = {Pacific J. Math.},
  FJOURNAL = {Pacific Journal of Mathematics},
    VOLUME = {307},
      YEAR = {2020},
    NUMBER = {1},
     PAGES = {13--51},
      ISSN = {0030-8730,1945-5844},
   MRCLASS = {11H56 (11F22 11F27 20F05 20F55 51M10)},
  MRNUMBER = {4131798},
MRREVIEWER = {Haowu\ Wang},
       DOI = {10.2140/pjm.2020.307.13},
       URL = {https://doi.org/10.2140/pjm.2020.307.13},
}

@article {Borcherds00-refl-LL,
    AUTHOR = {Borcherds, Richard E.},
     TITLE = {Reflection groups of {L}orentzian lattices},
   JOURNAL = {Duke Math. J.},
  FJOURNAL = {Duke Mathematical Journal},
    VOLUME = {104},
      YEAR = {2000},
    NUMBER = {2},
     PAGES = {319--366},
      ISSN = {0012-7094,1547-7398},
   MRCLASS = {11H56 (11F11 11F27 11M36 51F15)},
  MRNUMBER = {1773561},
MRREVIEWER = {G.\ J.\ Heckman},
       DOI = {10.1215/S0012-7094-00-10424-3},
       URL = {https://doi.org/10.1215/S0012-7094-00-10424-3},
}

@article {Borcherds85-leech,
    AUTHOR = {Borcherds, R. E.},
     TITLE = {The {L}eech lattice},
   JOURNAL = {Proc. Roy. Soc. London Ser. A},
  FJOURNAL = {Proceedings of the Royal Society. London. Series A.
              Mathematical, Physical and Engineering Sciences},
    VOLUME = {398},
      YEAR = {1985},
    NUMBER = {1815},
     PAGES = {365--376},
      ISSN = {0962-8444,2053-9169},
   MRCLASS = {11H06},
  MRNUMBER = {791880},
}

@article {Borcherds87-aut-LL,
    AUTHOR = {Borcherds, Richard},
     TITLE = {Automorphism groups of {L}orentzian lattices},
   JOURNAL = {J. Algebra},
  FJOURNAL = {Journal of Algebra},
    VOLUME = {111},
      YEAR = {1987},
    NUMBER = {1},
     PAGES = {133--153},
      ISSN = {0021-8693},
   MRCLASS = {20B27 (51F15)},
  MRNUMBER = {913200},
MRREVIEWER = {J.\ S.\ Hsia},
       DOI = {10.1016/0021-8693(87)90245-6},
       URL = {https://doi.org/10.1016/0021-8693(87)90245-6},
}

@book {ConwaySloane99,
    AUTHOR = {Conway, J. H. and Sloane, N. J. A.},
     TITLE = {Sphere packings, lattices and groups},
    SERIES = {Grundlehren der mathematischen Wissenschaften [Fundamental
              Principles of Mathematical Sciences]},
    VOLUME = {290},
   EDITION = {Third},
      NOTE = {With additional contributions by E. Bannai, R. E. Borcherds,
              J. Leech, S. P. Norton, A. M. Odlyzko, R. A. Parker, L. Queen
              and B. B. Venkov},
 PUBLISHER = {Springer-Verlag, New York},
      YEAR = {1999},
     PAGES = {lxxiv+703},
      ISBN = {0-387-98585-9},
   MRCLASS = {11H31 (05B40 11H06 20D08 52C07 52C17 94B75 94C30)},
  MRNUMBER = {1662447},
MRREVIEWER = {Renaud\ Coulangeon},
       DOI = {10.1007/978-1-4757-6568-7},
       URL = {https://doi.org/10.1007/978-1-4757-6568-7},
}

@article {DeligneMostow86,
    AUTHOR = {Deligne, P. and Mostow, G. D.},
     TITLE = {Monodromy of hypergeometric functions and nonlattice integral
              monodromy},
   JOURNAL = {Inst. Hautes \'Etudes Sci. Publ. Math.},
  FJOURNAL = {Institut des Hautes \'Etudes Scientifiques. Publications
              Math\'ematiques},
    NUMBER = {63},
      YEAR = {1986},
     PAGES = {5--89},
      ISSN = {0073-8301,1618-1913},
   MRCLASS = {22E40 (32G20 33A30)},
  MRNUMBER = {849651},
MRREVIEWER = {Shigeaki\ Tsuyumine},
       URL = {http://www.numdam.org/item?id=PMIHES_1986__63__5_0},
}

@article {NeumaierSeidel83,
    AUTHOR = {Neumaier, A. and Seidel, J. J.},
     TITLE = {Discrete hyperbolic geometry},
   JOURNAL = {Combinatorica},
  FJOURNAL = {Combinatorica. An International Journal of the J\'anos Bolyai
              Mathematical Society},
    VOLUME = {3},
      YEAR = {1983},
    NUMBER = {2},
     PAGES = {219--237},
      ISSN = {0209-9683},
   MRCLASS = {51F15 (05C50 11F06 11H55 52A43)},
  MRNUMBER = {726460},
MRREVIEWER = {\`E.\ Vinberg},
       DOI = {10.1007/BF02579296},
       URL = {https://doi.org/10.1007/BF02579296},
}

@incollection {Thurston98,
    AUTHOR = {Thurston, William P.},
     TITLE = {Shapes of polyhedra and triangulations of the sphere},
 BOOKTITLE = {The {E}pstein birthday schrift},
    SERIES = {Geom. Topol. Monogr.},
    VOLUME = {1},
     PAGES = {511--549},
 PUBLISHER = {Geom. Topol. Publ., Coventry},
      YEAR = {1998},
   MRCLASS = {57M50 (20F67 22E40 53C21)},
  MRNUMBER = {1668340},
MRREVIEWER = {Feng\ Luo},
       DOI = {10.2140/gtm.1998.1.511},
       URL = {https://doi.org/10.2140/gtm.1998.1.511},
}

@article{ConwaySimons2001,
  author  = {Conway, John H. and Simons, Christopher S.},
  title   = {26 Implies the Bimonster},
  journal = {Journal of Algebra},
  volume  = {235},
  number  = {2},
  pages   = {805--814},
  year    = {2001},
  doi     = {10.1006/jabr.2000.8494}
}

@book{BrouwerCohenNeumaier89,
  title     = {Distance-Regular Graphs},
  author    = {Brouwer, Andries E. and Cohen, Arjeh M. and Neumaier, Arnold},
  publisher = {Springer},
  address   = {Berlin, Heidelberg},
  year      = {1989},
  series    = {Ergebnisse der Mathematik und ihrer Grenzgebiete},
  volume    = {18},
  isbn      = {978-3-540-50619-5},
  doi       = {10.1007/978-3-642-74341-2}
}

\end{document}